\RequirePackage{ifluatex}
\documentclass[12pt]{article}
\usepackage[margin=1.15in]{geometry}
\usepackage{graphics}
\usepackage{ifpdf}
\usepackage{amsmath,amssymb}
\usepackage{hyperref}
\usepackage{amssymb}
\usepackage{epsfig}
\usepackage{epstopdf}
\usepackage{tikz}
\usepackage{upgreek}
\usepackage{authblk}
\usepackage{nomencl}
\usepackage[maxbibnames=10, citestyle=numeric-comp,bibstyle=numeric, giveninits=true, backend=bibtex, sorting=none]{biblatex}
\usepackage[font=small,labelfont=bf]{caption}
\usepackage{titlesec}
\usepackage[titletoc,toc,title]{appendix} 
\usepackage{abstract}  
\usepackage{setspace}
\usepackage{xcolor}
\usepackage[]{subfig}
\usepackage{bm}
\usepackage{float}
\usepackage{enumerate}

\DeclareMathAlphabet      {\mathbfit}{OML}{cmm}{b}{it}
\usepackage{titlesec}
\usepackage{overpic}
\titleformat{\paragraph}
{\normalfont\normalsize\bfseries}{\theparagraph}{1em}{}
\titlespacing*{\paragraph}
{0pt}{3.25ex plus 1ex minus .2ex}{1.5ex plus .2ex}

\makeatletter
\def\@seccntformat#1{\@ifundefined{#1@cntformat}%
	{\csname the#1\endcsname\quad}
	{\csname #1@cntformat\endcsname}
}
\makeatother
\AtEveryBibitem{%
	\clearfield{isbn}
	\clearfield{issn}
	\clearfield{month}
}
\renewbibmacro{volume+number+eid}{%
	\printfield{volume}%
	\setunit{\space}%
	\printfield{number}%
	\setunit{\space}%
	\printfield{eid}
}
\DeclareFieldFormat[article]{number}{\mkbibparens{#1}}
\DeclareFieldFormat[article]{volume}{\mkbibbold{#1}}
\bibliography{References}

\makeatletter
\def\@seccntformat#1{\@ifundefined{#1@cntformat}%
{\csname the#1\endcsname\quad}  
{\csname #1@cntformat\endcsname}
}
\let\oldappendix\appendix 
\renewcommand\appendix{%
\oldappendix
\newcommand{\section@cntformat}{\appendixname~\thesection\quad}
}
\makeatother
\usepackage{cleveref} 

\def\centerarc[#1](#2)(#3:#4:#5);%
{
\draw[#1]([shift=(#3:#5)]#2) arc (#3:#4:#5);
}

\definecolor{ao}{rgb}{0.0, 0.5, 0.0}
\definecolor{OliveGreen}{rgb}{0,0.6,0}

\title{\Large{\textbf{Three-dimensional varying-order NURBS discretization method for enhanced IGA of large deformation frictional contact problems}}\vspace{-0.5em}}
\author[$^{\dag}$]{{Vishal Agrawal}\vspace{-0.75em}}
\affil[$^{\dag}$]{\footnotesize{{Department of Engineering Mechanics,}}} 
\affil[$\,$]{\footnotesize{{KTH Royal Institute of Technology, Stockholm, Sweden}}}
\affil[$\,$]{\footnotesize{Corresponding author email: \href{vishala@kth.se}{vishala@kth.se} \href{agrawalvishal34@gmail.com}{agrawalvishal34@gmail.com}}}
\affil[$\,$]{\footnotesize{Published in CMAME on March 5, 2025 with DOI: \href{https://doi.org/10.1016/j.cma.2025.117853}{10.1016/j.cma.2025.117853}}}
\date{}

\begin{document}
\maketitle
\sloppy	
\vspace{-3.5em}
\renewenvironment{abstract}
{\begin{quote}
    \noindent \rule{\linewidth}{.5pt}\par{\bfseries \abstractname.}}
	{\medskip\noindent \rule{\linewidth}{.5pt}
	\end{quote}
}

\newenvironment{noverticalspace}
{%
	\par
	\offinterlineskip
}
{\par}
\begin{abstract}
In this contribution, we introduce a varying-order (VO) NURBS discretization method to enhance the performance of the isogeometric analysis (IGA) technique for solving three-dimensional (3D) large deformation frictional contact problems involving two deformable bodies. Building on the promising results obtained from the previous work on the 2D isogeometric contact analysis~\cite{Agrawal2020}, this work extends the method's capability for tri-variate NURBS-based discretization. The proposed method allows for independent,  user-defined application of higher-order NURBS functions to discretize the contact surface while employing the minimum order NURBS for the remaining volume of the elastic solid. This flexible strategy enables the possibility to refine a NURBS-constructed solid at a fixed mesh with the controllable order elevation-based approach while preserving the original volume parametrization. The advantages of the method are twofold. First, employing higher-order NURBS for contact integral evaluations considerably enhances the accuracy of the contact responses at a fixed mesh, fully exploiting the advantage of higher-order NURBS specifically for contact computations. Second, the minimum order NURBS for the computations in the remaining bulk volume substantially reduces the computational cost inherently associated with the standard uniform order NURBS-based isogeometric contact analyses.

The capabilities of the proposed method are demonstrated using various contact problems between elastic solids with or without considering friction. The results with the standard uniform order of tri-variate NURBS-based discretizations are also included to provide a comprehensive comparative assessment. We show that to attain results of similar accuracy, the varying-order NURBS discretization uses a much coarser mesh resolution than the standard uniform-order NURBS-based discretization, hence leading to a major gain in computational efficiency for isogeometric contact analysis. The convergence study demonstrates the consistent performance of the method for efficient IGA of 3D frictional contact problems. Furthermore, the simplicity of the method facilitates its direct integration into the existing 3D NURBS-based IGA framework with only a few minor modifications. 
\\
\end{abstract}

\section{Introduction}
The significance of accurately modeling contact interactions has driven the development of advanced numerical methods in various engineering fields. As is widely known, the accuracy and robustness of the numerical methods in computational contact mechanics depend primarily on two main factors. The first is the type of method employed for enforcing contact constraints into the variational formulations of the problem~\cite{wriggers2006, DeLorenzis2014, DeLorenzis2017_Encycp}. The most popular choices are the penalty method~\cite{Hallquist1985, Wriggers1990}, Lagrange multiplier method~\cite{Bathe1985}, and a combination of both these methods, known as augmented Lagrange multiplier method~\cite{Alart1991, Simo1992, Pietrzak1999}, among others. The second aspect is the type of contact algorithm utilized for expressing the contact contributions to the discretized weak form in a specific manner. The widely adopted algorithms are: Node-to-Segment (NTS)~\cite{Hallquist1985, Wriggers1990, Papadopoulos1992}, Gauss-Point-to-Surface (GPTS)~\cite{Fischer2005, DeLorenzis2011, Temizer2011, Dimitri2014Tspline}, and mortar contact formulation~\cite{Puso2004, Yang2005, Yang2008, Tur2009, Popp2010, Popp2013}. For the relative merits and demerits of the above-mentioned methods and algorithms, we refer to the above-cited references and the textbooks of Laursen \cite{laursen2003} and Wriggers \cite{wriggers2006}.

Another factor that adversely affects the robustness and accuracy of the contact algorithms is geometrical discontinuity. It arises from the traditional $C^0$-continuous finite element (FE)-based faceted representation of the contact surface~\cite{wriggers2006, DeLorenzis2014, DeLorenzis2017_Encycp}. In order to alleviate the issues emanating from this discontinuity, various surface smoothing strategies have been introduced in the context of FE. They employed different CAD polynomials in their designs such as: Hermite polynomials \cite{Pietrzak1999}, Bernstein polynomials \cite{Wriggers2001}, B\'{e}zier interpolation functions~\cite{Lengiewicz2010}, splines polynomials ~\cite{El-Abbasi2001, Stadler2003}, Gregory patches \cite{Puso2002}, triangular B\'{e}zier patches~\cite{Krstulovi2002}, and subdivision surface scheme~\cite{Stadler2004}. Additionally, layer-wise enrichment of shape functions approach by Konyukhov and Schweizerhof \cite{Konyukhov2009, Konyukhov2013} can also be regarded as a surface smoothing technique. While these methods ensured the smooth representation of the FE discretized contact surface, they introduced an additional smooth layer on top of an existing FE mesh. Such layering complicated the implementation and data management, especially in the 3D cases. Furthermore, these methods failed to enhance the spatial convergence and struggled with large deformation or sliding contacts scenarios.

As an alternative to the FE discretization, Hughes et al.~\cite{hughes2005} introduced the non-uniform rational B-Splines (NURBS) based IGA technique. Attributed to the advantageous key features of its underlying basis functions, i.e., NURBS, viz., the ability to represent complex-shaped geometry in its original smooth CAD form, tailorable inter-element continuity, and the non-negativeness, IGA has established itself as an advantageous computational technology for various classes of problems. Among others, contact mechanics belong to one of these classes that has significantly benefited from the intrinsic features of the IGA technique. This is because the NURBS-based description naturally provides a smooth representation of the contact surface, with the desirable inter-element continuity. This description leads to the unique definition of contact kinematics, i.e. surface normal and tangential vector fields across the contact surface even at a very coarse mesh, also for the geometries having a very large curvature. Consequently, the NURBS-based discretization inherently eliminates the computational issues originating from the traditional C$^0$-continuous FE-based discretization, leading to superior quality results, especially for the large deformation and finite frictional sliding contact problems. Moreover, by enabling contact search on a patch level rather than at the element level, IGA considerably reduces the bookkeeping task associated with the local contact search and frictional sliding contact~\cite{DeLorenzis2014}.

To this date, substantial research progress has been made towards developing various isogeometric methods for the robust and efficient treatment of a wide range of contact problems with high accuracy. In the following, we provide a brief introduction to some advancements. Temizer et al.~\cite{Temizer2011} introduced isogeometric mortar contact formulation for two-dimensional (2D) frictionless contact problems. De Lorenzis et al.~\cite{DeLorenzis2011} extended the formulation for large deformation frictional contact problems. Dittman et al.~\cite{Dittmann2014} further extended the isogeometric mortar contact formulation for thermo-mechanical contact and impact problems. Matzen et al.~\cite{Matzen2013} presented an isogeometric collocation contact formulation as an efficient alternative to the standard Galerkin formulation for frictionless contact problems. In the same direction, De Lorenzis et al.~\cite{DELORENZIS2015} introduced hybrid collocation-Galerkin and enhanced collocation approaches to circumvent the issues emanating from the strong imposition of Neumann boundary conditions and the enforcement of contact constraints in the isogeometric collocation framework. Later, Nguyen-Thanh et al.~\cite{Nguyen2019} introduced an adaptive isogeometric collocation method that enables the adaptive isogeometric analysis of 2D frictional contact problems. Seitz et al.~\cite{SEITZ2016} developed a dual mortar method for isogeometric analysis of mechanical and tied contact problems and conducted the first systematic investigation into the spatial converge rates of different order NURBS. We refer to the comprehensive review paper by De Lorenzis et al.~\cite{DeLorenzis2014} for a detailed overview of the advancement in the field of isogeometric contact analysis, and Das and Gautam~\cite{Das2024} for the detailed applications of isogeometric contact analysis.

Due to their unique intrinsic features, the NURBS-based discretizations provide significantly superior results in terms of accuracy, stability, and robustness for a wide range of contact problems as compared to its Lagrange polynomial-based counterpart. However, due to their rigid tensor product nature, the NURBS discretization lacks local mesh refinement. Additionally, they do not support non-uniform order elevation-based refinement~\cite{hughes2005, Cottrell2007, Temizer2014}. Since contact is mainly a surface-dominated phenomenon, employing uniform discretization in the region away from the contact surface is not desirable, particularly from the analysis viewpoint. 

In this direction, the various strategies enabling local mesh refinement in the vicinity of the contact surface have been introduced: T-spline for efficient IGA of frictionless contact and cohesive interface debonding problems by Dimitri et al.~\cite{Dimitri2014Tspline, Dimitri2017}, hierarchical NURBS-based local refinement approach for frictional contact problems Temizer and Hesch~\cite{Temizer20161}, and an adaptive local surface refinement based on the locally refined NURBS for frictional sliding contact by Zimmerman and Sauer~\cite{Zimmermann2017}. 

Moreover, various FE-IGA methods have been introduced as an alternative to uniform NURBS discretization. Corbett and Sauer~\cite{Corbett2014, Corbett2015} introduced the NURBS-enriched contact finite elements for frictional contact and mixed-mode cohesive debonding problems. Their approach used the NURBS functions to describe the bottom or contact surface of the linear-order finite element. Maleki-Jebeli et al.~\cite{MALEKIJEBELI2018} presented the hybrid isogeometric-FE discretization method for cohesive interface contact/debonding problems. They used NURBS in the contact/cohesive interface regions having high-stress gradients and linear- or quadratic-order Lagrange finite elements to discretize the remaining domain. Otto et al.~\cite{Otto2018} introduced a coupled FE-NURBS discretization approach to model the contact problems between the arbitrary order of finite element discretized bodies. In this, an auxiliary NURBS layer is added between the higher-order FE discretized contacting bodies. In Dias et al.~\cite{DIAS2019}, the geometric mapping of the curved surfaces of high-order FE contact elements is performed with NURBS. The central feature of these approaches has been to combine the accuracy and stability features of the NURBS for the contact computations with the efficiency characteristic of the FE method for the bulk computations. 

But, still, no such effort that fully exploits the key advantages of the IGA technique by directly allowing the non-uniform order elevation-based refinement of an elastic solid for isogeometric contact analysis has been made. We fill this research gap in the present work by introducing a novel three-dimensional varying-order NURBS discretization method. It enables direct usage of higher-order NURBS for the computations of quantities at the contact surface and lower-order NURBS for the computations of quantities in the remaining volume of an elastic solid. By retaining the intrinsic key features of the IGA technique, the proposed method aims to enhance the performance of the IGA technique for contact applications involving deformable solids. The proposed method extends the previous work focused on isogeometric contact analysis in two-dimensional settings~\cite {Agrawal2020}. In this contribution, we show the efficacy of the proposed method for the 3D large deformation frictional contact problems between deformable solids.

The remaining paper is structured as follows. Section~\ref{sec:contact_formulation} briefly describes the general formulation for the frictional contact problem, including NURBS and Gauss-point-to-surface based contact algorithm in the IGA setting. Next, Section~\ref{sec:VO_method} describes the proposed 3D VO-based NURBS discretization method and its implementation procedure into the existing IGA code. After that, Section~\ref{sec:numerical_example} demonstrates the performance and capabilities of the method using various numerical examples. Finally, Section~\ref{sec:conclusion} concludes the paper with the future directions.

\section{3D large deformation frictional contact problem} \label{sec:contact_formulation}
This section summarizes the mathematical and algorithm background for frictional contact between two bodies undergoing large deformations. The present formulation is based on Refs.~\cite{wriggers2006, laursen2003, unbiased2013, unbiased2015}; we refer to them for more details.

\subsection{Problem description and weak form} \label{subsec:problem_descritption}
\begin{figure}[!ht]
	\centering \subfloat{\includegraphics[width=0.85\linewidth]{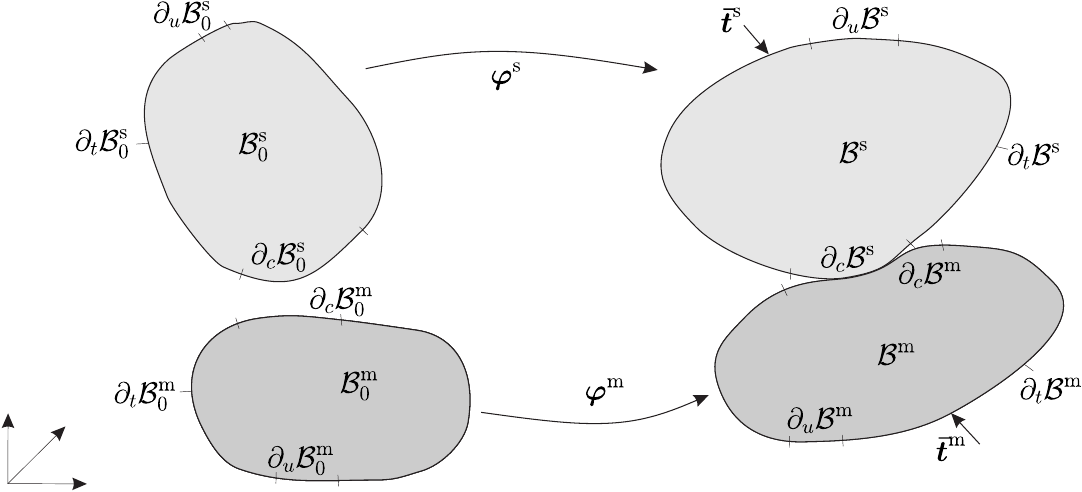}}
	\caption{Schematic illustration of the two-body contact in 3D.}
	\label{fig:Contact_deformation}
\end{figure}

It is assumed that two bodies come into contact and undergo large deformations, as shown in Fig.~\ref{fig:Contact_deformation}. In this work, according to the convention used in~\cite{wriggers2006}, one of them is denoted as the slave body $ \mathcal{B}^\textrm{s} $ and another as the master body $ \mathcal{B}^\textrm{m} $. Thus, in the following, a superscript $ k = \{\textrm{s}, \textrm{m}\} $ is used to specify them. It is assumed that the bodies come into contact by applying prescribed external surface traction and/or displacement. The boundary $ \partial \mathcal{B}^{k} $ of a body $ \mathcal{B}^{k} $ consists of the following three distinct parts: $\partial \mathcal{B}^{k} = \partial_t \mathcal{B}^{k}  \cup \partial_u \mathcal{B}^{k}  \cup \partial_c \mathcal{B}^{k} \,, ~~ \forall\, t \in  T$, where $ \partial_t \mathcal{B}^{k} $ and $ \partial_u \mathcal{B}^{k} $ are the boundaries where the traction and displacement are prescribed, and $ \partial_c \mathcal{B}^{k} $ is the potential surface where the contact occurs. Further, it is assumed that these boundaries satisfy $\partial_t \mathcal{B}^{k} \, \cap \, \partial_u \mathcal{B}^{k} \, \cap \partial_c \mathcal{B}^{k} = \emptyset\,,~~\forall\, t \in  T$. 

To identify the contact point on the master surface, it is parameterized using the convective coordinate vector $ \bm{\xi}_\textrm{m} = (\xi^1_{\mathrm{m}}, \xi^2_{\mathrm{m}}) $, as shown in Fig.~\ref{fig:friction_description}. Here, $ \xi_m $ denotes the parametric point computed at the master contact surface. With this, a contact point $ \bm{x}^\textrm{m} \in \Gamma^\textrm{m}_\textrm{c} $ on the master surface is defined by the parametric coordinates $ {\bm{\xi}}_{\textrm{m}} $ as $ \bm{x}^\textrm{m}(\bm{\xi_{\textrm{m}}}) $. The covariant tangent vectors at $ \bm{\xi_{\textrm{m}}} $ are defined by $\bm{\uptau}_{\alpha} = \frac{\partial \bm{x}^\textrm{m}}{\partial {\xi_\textrm{m}^{\alpha}}} =  \bm{x}^\textrm{m}_{,\alpha}\,, ~~\, \alpha = 1,2 \,$.  The contravariant tangent vectors can be obtained directly using the following relation $\bm{\uptau}^{\alpha} = m^{\alpha \beta}\, \bm{\uptau}_{\alpha}\,,$, where $ m^{\alpha \beta}  = [m_{\alpha \beta}] ^{-1}$ is the component of the inverse of the metric tensor defined by $ m_{\alpha \beta} = \bm{\uptau}_{\alpha} \cdot \bm{\uptau}_{\beta}$. The outward unit normal vector to the master surface at a parametric point $ \bm{x}^\textrm{m}(\bm{\xi_{\textrm{m}}}) $ is given by
\begin{equation}\label{eq:normal_vector}
	\bm{n} = \frac{\bm{\uptau}_{1} \times \bm{\uptau}_{2}}{||\bm{\uptau}_{1} \times \bm{\uptau}_{2}||}\,.
\end{equation}

Given the two reference configurations $ \mathcal{B}_{0}^k,~k = \{s,m\}$, the prescribed displacement $ \bar{\bm{u}}^k $ on $ \partial_u \mathcal{B}^{k} $, traction $ \bar{\bm{t}}^k $ on $ \partial_t \mathcal{B}^{k} $, and body forces $ \rho^k \, \bar{\bm{b}}^k $ in $ \mathcal{B}^k $, the objective is to find the displacement field $ \bm{u}^k \in \mathcal{B}^k $, such that at an instance it satisfies the two sets $ k= \{s,m\} $ of the following equilibrium equations
\begin{equation}\label{eq:CBPV}
	\begin{aligned}
	\mathrm{div}\,\bm{\sigma}^k + \rho^k\, \bar{\bm{b}}^k &= \bm{0}~~~~~ \mathrm{in} ~~ \mathcal{B}^k\,,  \\
	\bm{u}^k &= \bar{\bm{u}}^k ~~~~~ \mathrm{on} ~~\partial_u \mathcal{B}^{k} \,,\\
	\bm{t}^k = \bm{\sigma}^k\bm{n}^k &= \bar{\bm{t}}^k~~~~~~ \mathrm{on}~~ \partial_t \mathcal{B}^{k} \,. \\
    \bm{t}^k = \bm{t}_c & ~~~~~~ \mathrm{on}~~ \partial_c \mathcal{B}^{k} \,. \\
	\end{aligned}
\end{equation}
subjected to the following contact constraints for normal and frictional contact:
\begin{eqnarray}\label{eq:CBVP_cont_conditions}
		\textrm{g}^{}_{\textrm{N}}  \geq 0\,, & {t}_{\text{N}} \leq 0\,, &\textrm{g}^{}_{\textrm{N}} \, {t}^{}_{\textrm{N}} = 0\,, ~~~ \forall \bm{x}^{\mathrm{s}}_{c} \in \partial \mathcal{B}^\mathrm{s} \\
		\gamma \geq 0\,,&\Phi_{} \leq 0\,, &{\gamma}\, \Phi_{} = 0\,, ~~~ \forall \bm{x}^{\mathrm{s}}_{c} \in \partial \mathcal{B}^\mathrm{s}\,.
\end{eqnarray}
The first set of contact constraints are known as Karush-–Kuhn-–Tucker (KKT) conditions in the optimization theory~\cite{wriggers2006}. In the above equation, $\bm{\sigma}$ is the Cauchy stress tensor, and $\bar{\bm{b}}$ is body force. Moreover, $\bm{t}_c$ is contact traction, $\textrm{g}^{}_{\textrm{N}}$ is the normal gap, ${t}^{}_{\textrm{N}}$ is normal traction computed at the contact point, and $\gamma$ and $\Phi_{}$ are slip function and slip parameter, described in the next section. Since the partial differential equation in Eq.~\ref{eq:CBPV}$_1$ requires the displacement field to be at least twice differentiable, it imposes a strong $ C^1$-continuity requirement at every point of the domain $ \mathcal{B}^k $. This continuity condition can be relaxed by using the weak form of contact boundary problem, as outlined below
\begin{equation}\label{eq:weak_form_CBVP}
	\delta \mathcal{W} = \sum_{k}^{s,m} \left(  \delta \mathcal{W}^k_{\mathrm{int}} -  \delta \mathcal{W}^k_{\mathrm{ext}}\right) + \delta \mathcal{W}_{\mathrm{c}}  = 0 \,, ~~~~~\forall ~\delta \bm{\varphi}^k \in \mathcal{V}^{k}
\end{equation}
where 
\begin{eqnarray}
	\delta \mathcal{W}^k_{\mathrm{int}} &=& \int_{\mathcal{B}^k} \bm{\sigma}^k : \, \mathrm{grad} \,(\delta \bm{\varphi}^k) \, \mathrm{d}v\,,  \label{eq:int_potentials_CVBP} \\
	\delta \mathcal{W}^k_{\mathrm{ext}} &=&  \int_{\mathcal{B}^k} \rho^k\, \bar{\bm{b}}^k \cdot \delta \bm{\varphi}^k \, \mathrm{d}v + \int_{\partial_t \mathcal{B}^{k}} \bar{\bm{t}}^k \cdot \delta\bm{\varphi}^k \,\mathrm{d}a \label{eq:ext_potentials_CVBP}
\end{eqnarray}
and the virtual work due to contact traction is defined as
\begin{equation}\label{eq:contact_potential} 
	\delta \mathcal{W}_{\mathrm{c}} = - \int_{\partial_{c} \mathcal{B}^\mathrm{s}}  \boldsymbol{t}_{\mathrm{c}} \cdot \delta{\boldsymbol{\varphi}}^\textrm{s}  ~\textrm{d}a \, + \, \int_{\partial_{c} \mathcal{B}^\mathrm{s}} \boldsymbol{t}_{\mathrm{c}} \cdot \delta{\boldsymbol{\varphi}}^\textrm{m}  ~  \textrm{d}a\,.
\end{equation}

\subsection{Closest projection point procedure}
\begin{figure}[!b]
	\centering
	\subfloat{\includegraphics[width=0.55\linewidth]{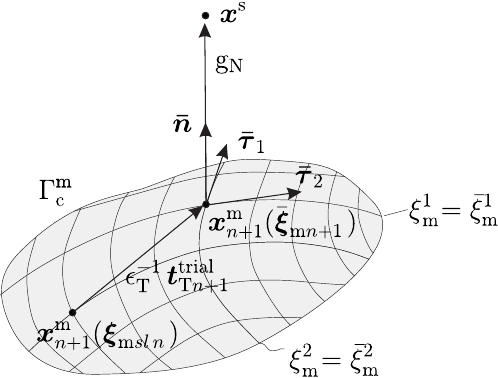}}
	\caption{A schematic illustration of the three-dimensional frictional contact quantities for a given slave point $ \bm{x}^{\mathrm{s}} $. The master contact surface $ \Gamma^{m}_{\mathrm{c}} $ is parametrized by the coordinates $ \xi^i $. The quantities computed at the projection point $ \bar{\bm{\xi}}_{\mathrm{m}} $ are denoted with bar.}
	\label{fig:friction_description}
\end{figure}
To distinguish whether the contact between the two surfaces is active or not, the closest normal gap function is determined. For the determination of a contact point $ \bm{x}^\textrm{m} $ on the master surface, a distance function is utilized: $\hat{d}(\bm{\xi}_\textrm{m}):= || \bm{x}^{\textrm{s}} - {\bm{x}}^\textrm{m}(\bm{\xi}_{\textrm{m}})||$. The necessary condition that yields the closest contact point ${\bm{x}}^\textrm{m}(\bm{\xi}_{\textrm{m}}) $  onto the master surface corresponding to a given slave point $\bm{x}^{\textrm{s}}$ is obtained by the minimum of the distance function $	\frac{d}{d{\xi}_{\textrm{m}}^{\alpha}} \, \hat{d}(\bm{\xi}_\textrm{m}) = [\bm{x}^{\textrm{s}} - {\bm{x}}^\textrm{m}(\bm{\xi}_{\textrm{m}})] \cdot  \bm{\uptau}_{\alpha}(\bm{\xi}_{\textrm{m}}) = 0\,.$ 

The corresponding parametric point $ \bar{\bm{\xi}}_{\textrm{m}} $ at which $ \frac{d}{d{\xi}_{\textrm{m}}^{\alpha}} \, \hat{d}({\bar{\bm{\xi}}}_\textrm{m}) = 0 $ is referred as {closest projection point}. Since the above equation is nonlinear, the Newton-Raphson (NR) method is used to solve it in this work. Note that in the subsequent description, all the quantities computed at the projection point will be denoted with a bar, e.g. $ {\bm{x}}^\textrm{m}(\bar{\bm{\xi}}_{\textrm{m}}) = \bar{\bm{x}}^{\mathrm{m}} $ and similarly unit normal vector $ \bm{n} $ at $ \bar{\bm{\xi}}_{\textrm{m}}  $ is denoted by $ \bar{\bm{n}} $, as shown in Fig.~\ref{fig:friction_description}.

With this, the normal gap $ \textrm{g}^{}_{\textrm{N}} $ between a slave point $ \bm{x}^{\textrm{s}} $ and projection point $ \bar{\bm{x}}^{\mathrm{m}} $ can be computed with
\begin{equation}\label{eq:gN}
	\textrm{g}_{\textrm{N}}  = (\bm{x}^\textrm{s} - \bar{\bm{x}}^\textrm{m})\cdot \bar{\bm{n}}\,,~~~~~~~~~~\forall\, \bm{x}^{\mathrm{s}} \in \partial \mathcal{B}_{c}^\mathrm{s}\,.
\end{equation}

\subsection{Penalty regularized contact constraints} \label{sec:contact_variables}
Since the contact constraints given by Eq.~\ref{eq:CBVP_cont_conditions} cannot be directly incorporated into the weak form, their regularization is needed. In the present work, we adopt the penalty method for this, and penalty regularized normal contact constraint is defined as
\begin{equation}\label{eq:normal_contact}
	\bm{t}_{\text{N}} = 
	\left\{\begin{array}{lr}
		-\epsilon_{\text{N}} \, \textrm{g}_{\textrm{N}}\, {{\bar{\bm{n}}}},  & \textrm{if}~ \textrm{g}_{\textrm{N}} < 0\\
		\bm{0},  & \textrm{otherwise}\,,
	\end{array}	\right.
\end{equation}
where $ \epsilon_{\text{N}} > 0 $ denotes the normal penalty parameter. 

Following~\cite{wriggers2006, unbiased2015}, the tangential traction is determined based on the changes in stick-slip status with
\begin{equation}\label{eq:tangential_traction}
	\bm{t}_{\textrm{T}} = 
	\left\{ \begin{array}{lr}
		\bm{t}_{\text{T}n+1}^{\textrm{trial}}, ~~ &  \textrm{if}~ \Phi_{n+1}^{\textrm{trial}} \leq 0\,,\\
		\mu_f\, {t}_{\textrm{N}n+1}\,  \bm{n}_{\textrm{T}n+1} ,~~ & \textrm{otherwise}\,. 
	\end{array}	\right.
\end{equation}
where $ \bm{n}_{\textrm{T}n+1} = {\bm{t}_{\text{T}n+1}^{\textrm{trial}}}/{||\bm{t}_{\text{T}n+1}^{\textrm{trial}}||} $, $\mu_f$ Coulomb's frictional coefficient, and $\bm{t}_{\text{T}n+1}^{\textrm{trial}}$ is the penalty-regularized trial contact traction defined as $ \bm{t}_{\text{T}n+1}^{\textrm{trial}} = \epsilon_{\textrm{T}} \left(\, \boldsymbol{x}^{\textrm{m}}_{n+1}(\bar{\bm{\xi}}_{\mathrm{m}\,{n+1}}) \,  - \, \boldsymbol{x}^{\textrm{m}}_{n+1}({\bm{\xi}}_{\textrm{m}\,{sl}\,n})\, \right) $, where $\epsilon_{\textrm{T}} $ denotes tangential penalty parameter. In the context of the predictor-corrector algorithm that is used to handle frictional sliding, $\bm{t}_{\text{T}n+1}^{\textrm{trial}}$ is utilized to compute the trial slip function $ \Phi_{n+1}^{\textrm{trial}} = || \bm{t}_{\text{T}n+1}^{\textrm{trial}} || -\mu_f \, {t}_{\text{N}n+1}$. With this, the total contact traction vector is defined as $\bm{t}_{\mathrm{c}} = {t}_{\mathrm{N}} \, \bar{{\bm{n}}} - {t}_{\mathrm{T}}^{\alpha} \, \bar{\bm{\uptau}}_{\alpha}\,$. 

\subsection{NURBS discretized weak form} \label{sec:NURBS_weak_form}
\subsubsection{Geometric modeling using NURBS}
In this section, the essential ingredients to construct a three-dimensional geometry using tri-variate NURBS functions are summarized. For a comprehensive description, the readers are referred to the monographs by Pigel and Tiller~\cite{nurbsbook} and Cottrell et al.~\cite{Cottrell2009}, as well as a tutorial paper by Agrawal and Gautam~\cite{Agrawal2018}, which provides stepwise implementation details of the IGA technique. Note that throughout this discussion, the interpolation order of splines - linear, quadratic, cubic, etc. is represented as $ p = 1,\, 2,\, 3,\,  $, and so on. This is as per the notation introduced by Hughes et al.~\cite{hughes2005}.

A \textit{knot vector} in one-dimension is a set of non-decreasing values of coordinates in the parametric space $ \xi^i~(i=1,\,2\,, \mathrm{~and }~  3) $, given by 
\begin{equation}\label{eq:knot_vector}
\Xi^i = \left[\xi_1^i, \, \xi_{2}^i, \, \dots\,, \, \xi^i_l, \, \dots\,, \, {\xi_{n_i+p_i+1}^i} \right] \,,
\end{equation}
where $ \xi_l^i \in \mathbb{R} ~~\forall \,i $ is the $ l^{\textrm{\scriptsize{th}}} $ knot, $ l = 1,2,\dots, n_i + p_i + 1$, $ p_i $ is the interpolation order, and $ n_i $ denotes the total number of basis functions. The non-zero interval between by any two consecutive knot values $ [\xi_l^i, \, \xi_{l+1}^i] > 0 $ is called as the \textit{knot span}. Each non-zero knot span partitions the parametric space $ \xi^i $ into elements, which, using the spline functions, are mapped to the physical space.

A NURBS-constructed solid is defined by 
\begin{equation}\label{eq:NURBS_solid}
\mathbf{V}(\xi^1,\xi^2, \xi^3) = \sum_{l=1}^{n_1} \sum_{m=1}^{n_2} \sum_{n=1}^{n_3}  R_{lmn}^{p_1,p_2,p_3} (\xi^1,\xi^2,\xi^3)\,\bm{X}_{lmn}\,, 
\end{equation}
where $n_i$ denotes the control points defined along the $i^{\textrm{th}}$ direction, and $\mathbf{X}_{A} $ is the control point coordinate vector. The tri-variate NURBS function $R_{lmn}^{p_1,p_2,p_3}$ is defined via the tensor product of the uni-variate B-spline functions as
\begin{equation}\label{eq:trivariate_NURBS_basis}
R_{lmn}^{p_1,p_2,p_3} (\xi^1,\xi^2, \xi^3)  =  \frac{w_{lmn}}{W(\xi^1,\xi^2,\xi^3) }\, N_{l,p_1}(\xi^1)\, N_{m,p_2}(\xi^2)\, N_{n,p_3}(\xi^3)\,.
\end{equation}
where $N_{l,p_i}(\xi^i)$ is a B-spline basis function of order $p_i$ along $\xi^i$ direction, and the normalizing weighting function is defined as
\begin{equation}\label{eq:3D_weight_function}
W(\xi^1,\xi^2,\xi^3)= {\sum_{l=1}^{n_1}\sum_{m=1}^{n_2} \sum_{n=1}^{n_3} \,w_{lmn}\, N_{l,p_1}(\xi^1)\, N_{m,p_2}(\xi^2)} \, N_{n,p_3}(\xi^3)\,.
\end{equation}

One of the key advantages of adopting CAD polynomials for IGA is that the constructed geometries can be refined in a straightforward fashion. The refinement strategies enrich the underlying basis functions, allowing the accurate capture of the sharp variations in the solution field~\cite{hughes2005, Cottrell2007}. The two basic methods that are used to refine the geometry are knot insertion, which introduces additional knots or elements in a given knot vector, and order elevation, which, as the name implies, uniformly elevates the interpolation order of the basis functions. A specific application of order elevation before the knot insertion leads to another strategy, namely $ k- $\textit{refinement}. These strategies provide control over the existing knot span or mesh size, inter-element continuity, and interpolation order of the underlying basis functions, while with their applications, the geometry and its parametrization remain intact. For further details on refining a solid constructed with NURBS, we refer to Pigel and Tiller~\cite{nurbsbook} and Cottrell et al.~\cite{Cottrell2009}.

In IGA, the NURBS-based description of the contact surface is directly inherited from the NURBS-described solid. The contact surface of NURBS discretized solid is defined by 
\begin{equation}\label{eq:phsical_point}
\begin{aligned}
\mathbf{S}(\xi^1,\xi^2) =&  \sum_{l=1}^{n_1} \sum_{m=1}^{n_2} R_{lm}^{p_1,p_2} (\xi^1,\xi^2) \, \mathbf{X}_{lm} 
 \end{aligned}
\end{equation}
where $\mathbf{X}_{lm} \subset \mathbf{X}_{lmn}$, and the bi-variate NURBS function is defined as
\begin{equation}\label{eq:bivariate_NURBS_basis}
R_{lm}^{p_1,p_2} (\xi^1,\xi^2) =  \frac{w_{lm}}{W(\xi^1,\xi^2) }\, N_{l,p_1}(\xi^1)\,N_{m,p_2}(\xi^2)\,.
\end{equation}

Following the isogeometric finite element notation, each domain $ \mathcal{B}^k,~k=\{s,m\} $ can be described in terms of $ n_{\mathrm{el}}^k $ number of NURBS elements. Consequently, the coordinates of an arbitrary point in the reference and current configurations are given as
\begin{equation}\label{eq:phsical_point2}
\begin{aligned}
\bm{X}^{ke} =& \sum_{A=1}^{n^{ke}_{\mathrm{cp}}} R^{e}_A(\bm{\xi}^k)\, \mathbf{X}^{k}_{A} =  \mathbf{R}^{ke}(\boldsymbol{\xi}^k)\,\mathbf{X}^{ke}\,,~~~~~ {\in}~~ \Omega^{ke}_0\,,  \\
 \bm{x}^{ke} =& \sum_{A=1}^{n^{ke}_{\mathrm{cp}}} R^{k}_A(\bm{\xi}^k)\, \mathbf{x}^{k}_{A}  = \mathbf{R}^{ke}(\boldsymbol{\xi}^k)\,\mathbf{x}^{ke}\,,~~~~~ {\in}~~ \Omega^{ke}\,,
 \end{aligned}
\end{equation}
where $ n_{\mathrm{cp}}^{ke} = \prod_{i=1,2,3} \, (p_i +1) $ denotes the control points whose corresponding NURBS basis functions $ R^k(\bm{\xi}^k) $ have the local support in an element. Moreover, $ \mathbf{R}^{ke} $ represents the array of NURBS basis functions, and $ \mathbf{X}^{ke} $ and $ \mathbf{x}^{ke} $ are coordinate vectors of the control points in the reference and current configurations. These arrays are defined as
\begin{equation}\label{eq:diff_arrays}
	\mathbf{R}^{ke} = 
	\begin{bmatrix} 
		R_1^k \, \mathbf{I}_d  \\	
		\\
		\vdots \\
		\\
		R_{n_{\mathrm{cp}}^{ke}}^k \, \mathbf{I}_d
	\end{bmatrix}\,, ~~~~~	
	\mathbf{X}^{ke} = 
	\begin{bmatrix} 
		\mathbf{X}_1^k \\	
		\\
		\vdots \\
		\\
		\mathbf{X}_{n_{\mathrm{cp}}^{ke}}^k
	\end{bmatrix}\,,	~~~~~\mathrm{and}~~
	\mathbf{x}^{ke} = 
	\begin{bmatrix} 
		\mathbf{x}_1^k \\
		\\
		\vdots \\
		\\
		\mathbf{x}_{n_{\mathrm{cp}}^{ke}}^k
	\end{bmatrix}\,.
\end{equation}
where $ \mathbf{I}_d $ is the identity matrix of size $3\times 3$, and the tri-variate NURBS functions. 

\subsubsection{Discretized weak form}\label{sec:final_weak_form}
Using the isoparametric concept, the NURBS-based approximation of the  displacement field $ \bm{u}^{ke} $ is given as
\begin{equation}\label{eq:disp_field}
 \boldsymbol{u}^{ke} =  \sum_{A=1}^{n^{ke}_{\mathrm{cp}}} R^{k}_A(\bm{\xi}^k)\, \mathbf{u}^{k}_{A} = \mathbf{R}^{ke}(\boldsymbol{\xi}^k)\,\mathbf{u}^{ke}\,, ~~~~~~  \textrm{where}~~~\mathbf{u}^{ke} = \mathbf{x}^{ke} - \mathbf{X}^{ke}\,.
\end{equation} 
Here, $ \mathbf{u}^{ke} $ denotes the displacement vector of control points and is defined as $ \mathbf{u}^{ke} = [\bm{u}_1^k,\, \dots, \bm{u}_{n_{\mathrm{cp}}^{ke}}^k]^{\mathrm{T}} $. Following the Galerkin approach, the same NURBS functions are utilized for the approximation of the test function, $ \delta \bm{\varphi}^{ke} $ as 
\begin{equation}\label{eq:disp_variation}
	\delta \boldsymbol{u}^{ke} = \mathbf{R}^{ke}(\boldsymbol{\xi}^k)\,\delta\mathbf{u}^{ke}.
\end{equation}

Using the definition of the displacement field, the NURBS discretized weak form for the contact boundary problem involving two deformable solids is given by
\begin{equation}\label{eq:global_force_vector}
	\sum_{k}^{s,m} 	(\delta \mathbf{u}^{k})^{\mathrm{T}} \left[  \mathbf{f}_{\textrm{int}}^k -  \mathbf{f}_{\textrm{ext}}^k +  \mathbf{f}_{\textrm{c}}\right]  = \mathbf{0} \,,~~~~\forall \delta \bm{u}^k \in \mathcal{V}^k
\end{equation}
where $ \mathbf{f}_{\mathrm{int}} $, $ \mathbf{f}_{\mathrm{ext}} $, and $ \mathbf{f}_{\mathrm{c}} $ are identified as the global internal force vector, global force vector due to external loads, and global contact force vector, respectively, corresponding to the two-body system. 

Since the variation of displacement is arbitrary, the Eq.~\ref{eq:global_force_vector} can be re-written as $\mathbf{f}_{\textrm{int}}^k -  \mathbf{f}_{\textrm{ext}}^k +  \mathbf{f}_{\textrm{c}}  = \mathbf{0}$, where these global force vectors are given by
\begin{equation}\label{eq:fint_element}
\begin{aligned}
	 \mathbf{f}_{\textrm{int}}^{ke} =& \int_{\Omega^{ke}} (\mathbf{B}^{ke})^{\mathrm{T}} \bm{\sigma}^{ke} \, \mathrm{d}\Omega\,, \\
  \mathbf{f}_{\textrm{ext}}^{ke} =& \int_{\Omega^{ke}} \rho^k\, (\mathbf{R}^{ke})^{\mathrm{T}} \, \bar{\bm{b}}^k \,\mathrm{d}\Omega + \int_{\Gamma_t^{ke}} (\mathbf{R}^{ke})^{\mathrm{T}}\, \bar{\bm{t}}^k \, \mathrm{d}\Gamma \\
  \mathbf{f}_{\textrm{c}}^{\mathrm{s}e} =& - \int_{\Gamma^{\textrm{s}e}_c}  (\mathbf{R}^{\textrm{s}e})^\textrm{T} \, \boldsymbol{t}_c ~\textrm{d}\Gamma\,, \\
  \mathbf{f}_{\textrm{c}}^{\mathrm{m}e} =& \int_{\Gamma^{\textrm{s}e}_c}  (\mathbf{R}^{\textrm{m}e})^\textrm{T} \, \boldsymbol{t}_c ~\textrm{d}\Gamma\,,
\end{aligned}  
\end{equation}
where $\mathbf{B}^{ke}$ is an array containing $n_{\mathrm{cp}}^{ke}$ tri-variate NURBS functions having local support in a NURBS element.

In this work, we use the Neo-Hookean material model to describe the large deformation of the bodies, and the corresponding constitutive relation is given by~\cite{bonet1997}
\begin{equation}\label{eq:Strain_energy_density_function}
\bar{W} = \frac{\lambda}{2}(\mathrm{ln}\,J)^2 + \frac{\mu}{2}(\mathrm{tr}(\bm{FF}^{\mathrm{T}}) - 3) - \mu(\mathrm{ln}J)\,,
\end{equation}
where $\bm{F}$ denotes the deformation gradient, $J = \mathrm{det}(\bm{F})$ the volume change, and  $ \lambda $ and $ \mu $ are the Lam\'{e} constants, i.e. bulk and shear moduli, and are given in terms of Young's modulus $  E $ and Poisson's ratio $ \nu $ as $\lambda = \frac{2\mu\nu}{ (1-2\nu)}\,, ~~\mathrm{and}~~ \mu =\frac{ E} {2(1+\nu)}\,.$ 

\subsection{Contact algorithm}
In this study, we adopt the classical full-pass version of the contact algorithm introduced by Sauer and De Lorenzis~\cite{unbiased2015}, owing to its simplicity and efficiency features over the existing approaches~\cite{wriggers2006, DeLorenzis2017_Encycp}. This algorithm was first proposed by Fischer and Wriggerss~\cite{Fischer2005} in the context of FE contact analysis and later extended by Temizer et al.~\cite{Temizer2011}, De Lorenzis et al.~\cite{DeLorenzis2011, DeLorenzis2012} and Dimitri et al.~\cite{Dimitri2014Tspline, Dimitri2014, Dimitri2017} for isogeometric contact analysis. Following~\cite{Dimitri2014Tspline, Dimitri2014, Dimitri2017}, we also denote this contact algorithm as ``Gauss-point-to-surface (GPTS)". In this, the contact constraints, impenetrability, and sticking given in Eq.~(\ref{eq:CBVP_cont_conditions}) are independently enforced at each quadrature point. It is noted that the GPTS-based contact algorithm passes the contact patch test within the quadrature error, which ensures the convergence of the solution upon mesh refinement~\cite{Fischer2005, Zavarise2009}. 

It is noted that the mathematically more sophisticated mortar contact algorithms, e.g. by Temizer et al.~\cite{Temizer2011, Temizer2012}, De Lorenzis et al.~\cite{DeLorenzis2011, DeLorenzis2012}, Seitz et al.~\cite{SEITZ2016}, and Duong et al.~\cite{Duong2018}, that have been introduced for the treatment of contact within the framework of IGA can also be employed. However, employing such contact algorithms lies beyond the scope of the present work. For an extensive reviews and recent developments on the treatment of contact problems using different mortar methods, the contributions by Popp and Wall~\cite{Popp2014} and Popp~\cite{Popp2018} can be referred.

The contact force vectors for the slave and master body's NURBS discretized contact element are given in Eq.~\ref{eq:fint_element}. For the sake of completeness, the contact tangent matrices that are obtained from the linearization of the contact force vectors are provided in Appendix~\ref{appendix:contact_matrices} for both the slave and master contact surfaces. 

\section{Varying order NURBS discretization method} \label{sec:VO_method}
For contact problems where the surface effects are predominant, enhancing the accuracy of contact integral evaluations is essential. This can be achieved by employing higher-order basis functions at a fixed mesh. However, in the current NURBS-based discretization approach, the interpolation order of the NURBS functions for both the contact boundary and the bulk of the contact body can only be elevated uniformly due to the rigid tensor product structure of NURBS~\cite{DeLorenzis2014}. As a result, higher-order NURBS, intended specifically for accurate contact computations, must also be applied to the bulk regions of the contact body that are far from the surface. This uniform elevation is inefficient from a computational standpoint, as it increases the cost significantly without providing added benefits in regions away from the contact. Additionally, it undermines one of the key strengths of isogeometric analysis (IGA), which is the ability to represent geometries with high accuracy, even with coarse meshes. Therefore, it is necessary to develop a method that overcomes the limitations of uniform order elevation in NURBS discretized geometries and enhances the performance of IGA in computational contact mechanics.

To address the above challenges, this work introduces a varying-order (VO) NURBS discretization method, building upon the two-dimensional VO approach in~\cite{Agrawal2020}. The proposed method allows user-defined higher-order NURBS basis functions to be applied specifically for contact integrals, while the bulk computations use the minimum-order NURBS that still represent the contact geometry accurately. This way, VO discretization enables the application of order elevation-based refinement in a controlled manner for isogeometric analysis. In the following, the theory of the proposed method is detailed for tri-variate NURBS constructed geometry.  Next, in Section~\ref{sec:implementation_aspect}, some remarks on integrating the proposed method into the standard NURBS-based IGA framework are given to summarize its implementation aspect.

\subsection{Three-dimensional VO NURBS discretization} \label{section:3D_VO}
\begin{figure}[!htbp]
	\centering
	\subfloat[]{\includegraphics[scale=0.81, trim={100 175 100 57}, clip]{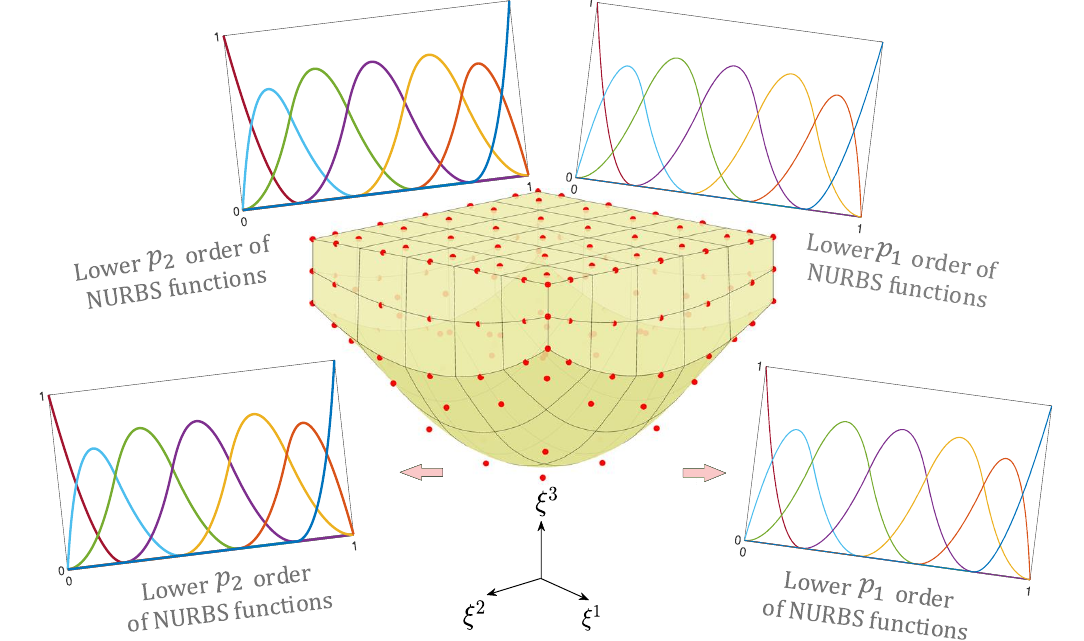}\label{fig:3D_VO_geometry_UO_body}} \\
	\subfloat[]{\includegraphics[scale=0.81, trim={100 175 100 57}, clip]{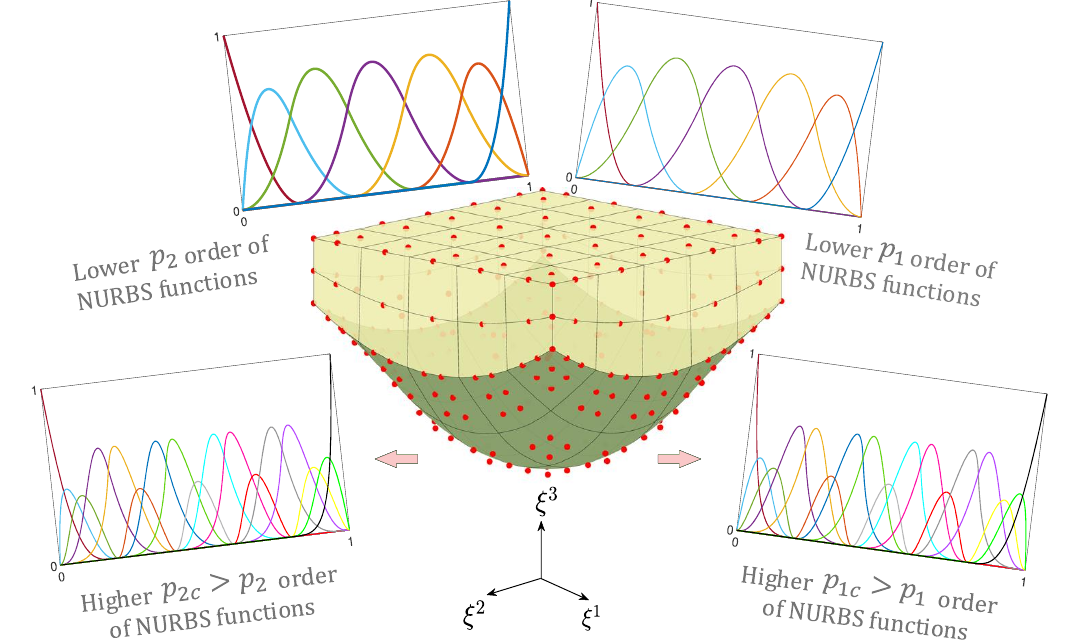}\label{fig:3D_VO_geometry_VO_body}}
	\caption{Schematic illustration of the standard and VO-based NURBS discretization of an example tri-variate geometry. (a) Standard NURBS discretization with minimum $ p_1 = p_2 = 2$ and $ p_3 = 1 $ order NURBS along the $ \xi^1 $, $ \xi^2 $, and $ \xi^3 $ parametric directions at a given mesh, where $ \Xi^1 = \Xi^2 = [0,0,0,1,2,3,4,5,5,5] $, and $ \Xi^3 = [0,0, 1, 2,2] $. (b) VO NURBS discretization, where the higher-order NURBS, i.e. N$ _{2 \cdot 1} $ are used for the description of the contact surface in both the $ \xi^1 $ and $ \xi^2 $ parametric direction, and the minimum N$_2 $ NURBS in its remaining volume.}\label{fig:3D_VO_geometry}
\end{figure}

The basic concept of the varying order NURBS discretization method for a contact geometry in three-dimension is illustrated in Fig.~\ref{fig:3D_VO_geometry}. Let $ p_1 $, $ p_2 $, and $ p_3 $ are the minimum interpolation order of the NURBS functions that are capable of accurately representing the given three-dimensional CAD geometry exactly. The coarsest mesh is given by the tensor product of the knot vectors defined along the $ \xi^1 $, $ \xi^2 $, and $ \xi^3 $ parametric directions, i.e. $ \Xi^1 \times \Xi^2 \times \Xi^3$. In the proposed method, the original minimum-order of NURBS discretized contact boundary surface, which is defined along the $ \xi^1 $ and $ \xi^2 $ parametric directions, is replaced with a new higher-order $ p_{1c} >p_1 $ and $ p_{2c} > p_2 $ NURBS surface such that it matches the volume parametrization as shown in Fig.~\ref{fig:3D_VO_geometry}. 

The new contact boundary surface of three-dimensional geometry can be constructed either using the $ k- $refinement or with a combination of $ k- $refinement and order elevation-based approaches, which, as a result, increase the inter-element continuity or interpolation order of NURBS on the contact surface. Thus, such a discretization enables the usage of higher-continuous and higher-order NURBS as basis functions for the computation of the quantities on the contact surface, while the minimum order of NURBS is utilized for the computation of the volume region that does not come into contact. 

The tri-variate NURBS basis functions for a VO discretized three-dimensional contact element are defined as 
\begin{equation}
\label{eq:Enriched_NURBS_element3D}
\begin{aligned}	
R_{1}^{p_{c1},p_{c2},p_3} (\xi^1, \xi^2, \xi^3) &:=  \frac{w_{111}}{W(\xi^1,\xi^2, \xi^3)} \, N_{1,p_{c1}}(\xi^1)\, N_{1,p_{c2}}(\xi^2)\, N_{1,p_3}(\xi^3)\,,  \\
\vdots~~~~~~~   & ~~~~~~~~~~~\vdots \\
R_{(p_{c1}+1)(p_{c2}+1)}^{p_{c1},p_{c2},p_3} (\xi^1, \xi^2, \xi^3) &:=  \frac{w_{(p_{c1}+1)(p_{c2}+1)1}}{W(\xi^1,\xi^2, \xi^3)} \, N_{p_{c1}+1,p_{c1}}(\xi^1)\, N_{p_{c2}+1,p_{c2}}(\xi^2)\, N_{1,p_3}(\xi^3)\,,  \\
R_{(p_{c1}+1)(p_{c2}+1)+1}^{p_{c1},p_{c2},p_3} (\xi^1, \xi^2, \xi^3) &:=  \frac{w_{112}}{W(\xi^1,\xi^2, \xi^3)} \, N_{1,p_{c1}}(\xi^1)\, N_{1,p_{c2}}(\xi^2) \, N_{2,p_3}(\xi^3)\,, \\
\vdots~~~~~~~   & ~~~~~~~~~~~\vdots \\
R_{n_{cp}^e}^{p_{c1},p_{c2},p_3} (\xi^1, \xi^2, \xi^3) &:= \frac{w_{(p_{c1}+1)(p_{c2}+1)2}}{W(\xi^1,\xi^2, \xi^3)} \, N_{p_{c1}+1,p_{c1}}(\xi^1)\, N_{p_{c2}+1,p_{c2}}(\xi^2)\, N_{2,p_3}(\xi^3)\,,
\end{aligned}	    
\end{equation} 
where $ n_{cp}^e = (p_{c1} +1)\times (p_{c2} +1) + (p_1 +1)\times (p_2 +1)\times p_3 $ denotes the total number of basis functions having local support in a VO based NURBS discretized contact element in 3D and the normalizing weight function is defined as
\begin{equation}\label{eq:3D_normalizing_weight}
\begin{aligned}	
W(\xi^1,\xi^2, \xi^3) = & \sum_{i=1}^{(p_{c1}+1)}\sum_{j=1}^{(p_{c2}+1)}  w_{ij1}  \,  
N_{i,p_{c1}}(\xi^1) N_{j,p_{c2}}(\xi^2) N_{1,p_3}(\xi^3) \\
	&+  \sum_{i=1}^{(p_{1}+1)}\sum_{j=1}^{(p_{2}+1)} w_{ij2}  \,  
	N_{i,p_{1}}(\xi^1) N_{j,p_{2}}(\xi^2) N_{2,p_3}(\xi^3)\,.
\end{aligned}
\end{equation}

The tri-variate NURBS functions defined in Eqs.~\ref{eq:Enriched_NURBS_element3D} exhibit the non-negativity property 
\begin{equation}\label{eq:non_negative}
R_{a}^{p_{c1},p_{c2},p_3} (\boldsymbol{\xi}) \geq 0 ~~~~\forall \,\boldsymbol{\xi} \in \Omega^e_c ~~~~~~\textrm{where}~ a = 1,2,\dots,n_{cp}^e\,,
\end{equation}
and also satisfy the partition of unity property
\begin{equation}\label{eq:partition_of_unity}
\sum_{a=1}^{n_{cp}^e} R_{a}^{p_{c1},p_{c2},p_3} (\boldsymbol{\xi}) = 1\,, ~~~~~~~ \forall \,\boldsymbol{\xi} \in \Omega^e_c\,.
\end{equation}

Since in the present work,$ p_{1c} = p_{2c} = p_c $ is taken in the considered numerical examples, the following notation is used for the VO-based NURBS discretization N$ _{p}-$N$ _{p_c \cdot p_s} $, where N$ _p $ is interpolation order of the NURBS used for the discretization of the volume, and $ N_{p_c} $ is for the contact surface. Moreover, the subscript $ p_s $ denotes the step number of order elevation applied to the N$ _{p_c} $ discretized contact surface.

Using the isoparametric concept, the same NURBS functions that are defined in Eq.~(\ref{eq:Enriched_NURBS_element3D}) for a VO discretized contact element are employed for the approximation of the unknown displacement field $ \mathbf{u}^e $, its variation $ \delta \mathbf{u}^e $, and to determine the current coordinates $ \bm{x}^e $ within each VO based NURBS discretized contact element $ \Omega_{\mathrm{c}}^e $ as 
\begin{equation}\label{eq:contact_surface}
\boldsymbol{u}^{e} = \mathbf{R}^e(\boldsymbol{\xi})\,\mathbf{u}^{e}\,,~~~\delta \boldsymbol{u}^{e} = \mathbf{R}^e(\boldsymbol{\xi})\,\delta\mathbf{u}^{e}\,,~~~\bm{x}^{e} = \mathbf{R}^e(\boldsymbol{\xi})\,\mathbf{x}^{e} ~~~~~ \forall \, \boldsymbol{\xi} \in \Omega_{\mathrm{c}}^e\,.
\end{equation} 
Here, the new basis function array $ \mathbf{R}(\boldsymbol{\xi}) $ is defined in a similar fashion to Eq.~(\ref{eq:diff_arrays}), except the modification that now it contains a total $n_{cp}^e = (p_{c1} +1)\times (p_{c2} +1) + (p_1 +1)\times (p_2 +1)\times p_3 $ number of trivariate NURBS functions $ R^{p_{c1}, p_{c2}, p_3}_a(\boldsymbol{\xi}) $ (\ref{eq:Enriched_NURBS_element3D}) having local support in a VO NURBS discretized contact element $ \Omega^e_c $. Moreover, the basis function array for the slave and master contact surface $ \mathbf{R}^{ke} $ in Eq.~(\ref{eq:fint_element}) makes use of bivariate $ R^{p_{c1}, p_{c2}}_a(\xi^1, \xi^2) $ NURBS functions for the evaluation of contact integrals in three-dimensions.

\subsection{Implementation details} \label{sec:implementation_aspect}
For integrating the VO NURBS discretization strategy into the standard NURBS-based isogeometric framework, only a few minor modifications are required. First of all, for a given mesh resolution, a $ p_{c} $ order of NURBS surface representing the contact surface of a NURBS-described solid is constructed. After that, the parametrization for an originally $ p $ order of NURBS discretized contact surface is replaced with that of the newly constructed $ p_{c} $ order of NURBS surface. For this, a certain number of conditions need to be fulfilled. The total number of control points defining the solid must be updated in such a manner that it allows the incorporation of the newly constructed $ p_{c} $ order of the new contact surface. This means that the connectivity array for contact elements must be adapted in a way that it contains the underlying control points of the VO-based discretized geometry. The derived contact element connectivity arrays can have a different length than the bulk element connectivity arrays. The trivariate NURBS basis functions defined in  Eq.~(\ref{eq:Enriched_NURBS_element3D}) are used for the evaluation of elemental quantities for VO-based NURBS described contact elements, while for the evaluation of the contact surface integrals, the standard bivariate $ p_c $ order of NURBS, functions are utilized. Except for these modifications, no other changes need to be made in the standard NURBS-based isogeometric framework. The local quantities, e.g. elemental tangent matrices and force vectors, are assembled to their global part in the same way as with the standard procedure. However, new control point connectivity arrays need to be utilized. In this work, a default $(p_c+1)$ number of Gauss-Legendre quadrature points in each parametric direction is employed to evaluate the contact integrals unless stated otherwise. For the bulk, $ (p+1) $ number of quadrature points along each parametric direction are utilized. Optimal quadrature rules~\cite{Hughes2010, AURICCHIO2012, FAHRENDORF2018}, which are well-suited for IGA, can also be opted for the reduced numerical evaluation. However, their application lies outside the scope of the present work.

\section{Numerical Examples} \label{sec:numerical_example}
In this section, we demonstrate the performance of the proposed varying-order-based NURBS discretization method for 3D large deformation contact problems, both with or without friction. The method's efficacy is assessed in terms of accuracy, robustness, and computational efficiency. For comparison, results obtained using the widely utilized standard N$_2$-based NURBS discretization in IGA are also included. It is noted that in this study, the interpolation order of NURBS is kept the same, i.e., $p=q$.
\begin{figure}[!ht]
	\centering
	\subfloat{\includegraphics[scale=0.26]{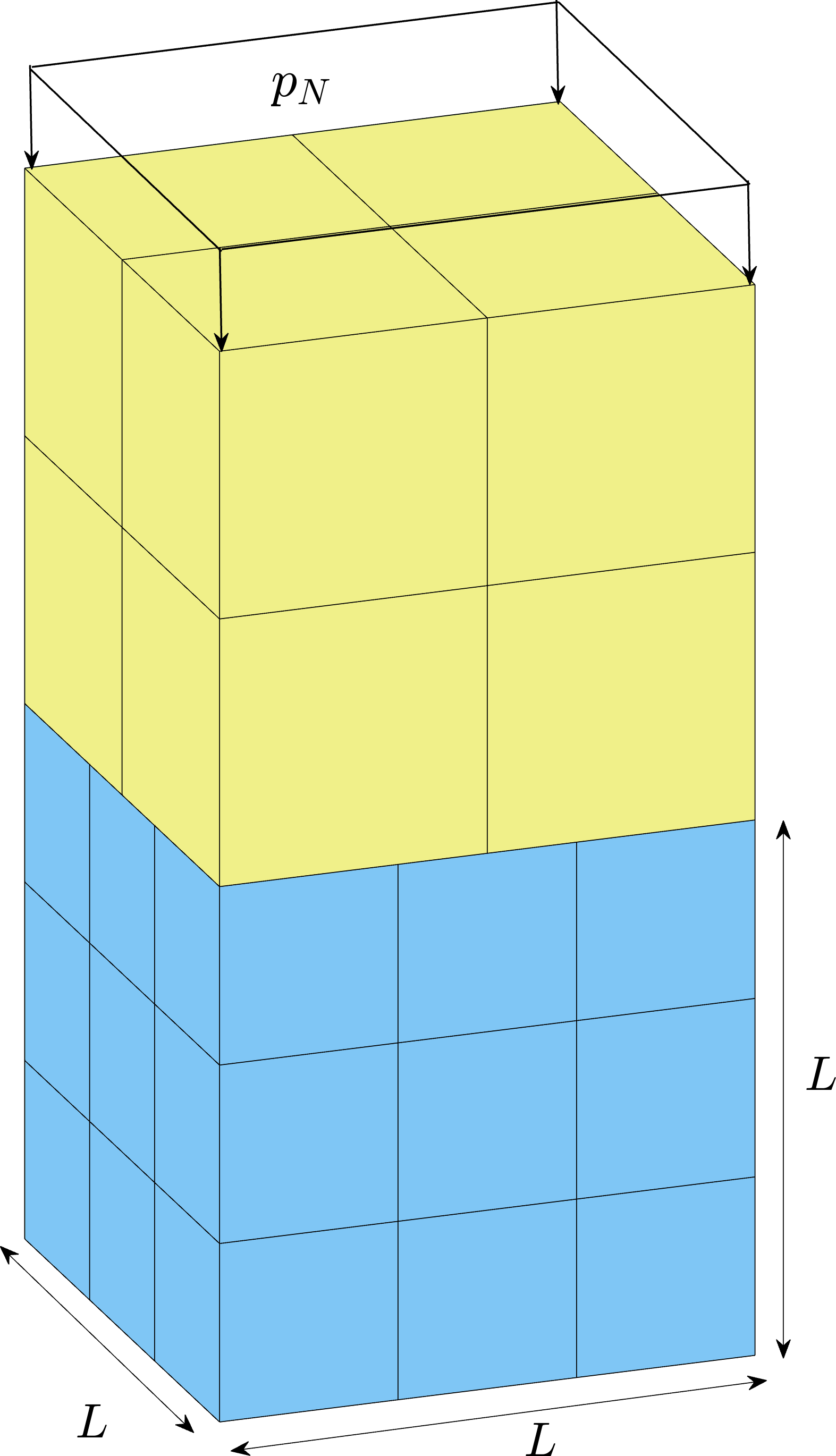}}
	\caption{Patch test: setup including the loading and coarsest mesh used for the cubes.} \label{fig:PathTest_setup_m1}
\end{figure}

    \subsection{Patch Test}
The contact patch test, originally introduced in~\cite{Papadopoulos1992}, is a benchmark used to evaluate whether a contact algorithm can exactly transmit a uniform pressure field across the contact interface.  Algorithms that fail this test introduce errors in the result at the interface, which may not diminish with the mesh refinement~\cite{Bathe2001}.

\begin{figure}[!htbp]
	\centering
    \begin{tikzpicture}
        \node[anchor=south west, inner sep=0] (img1) at (0,0) 
        {\subfloat{\includegraphics[scale=0.225]{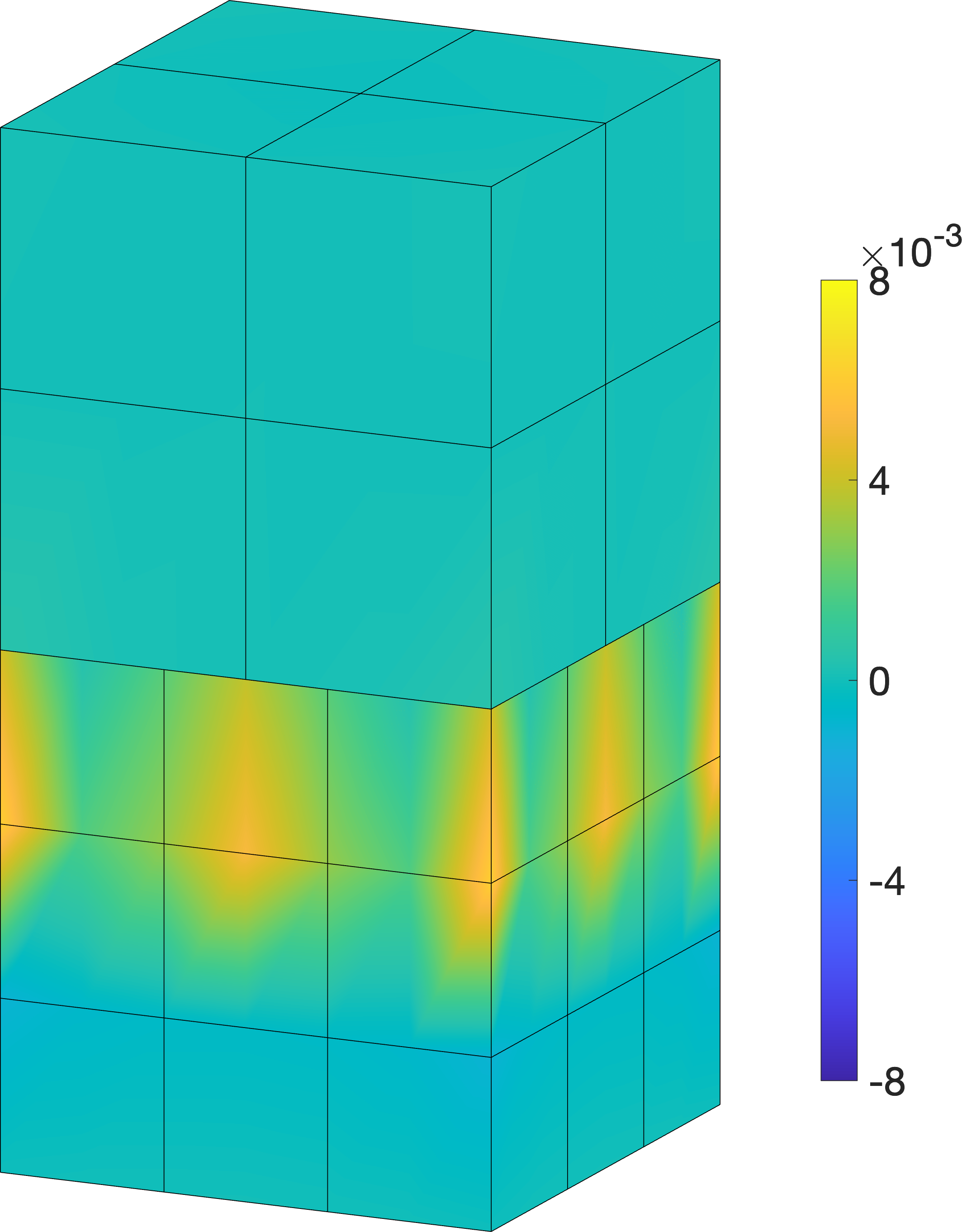}\label{fig:m1_ngps3}}};
        \node[font=\footnotesize, anchor=south] at (1.65, 6.215) { $n_\mathrm{gp}= (p+1)$ };
    \end{tikzpicture}~~
    \begin{tikzpicture}
        \node[anchor=south west, inner sep=0] (img2) at (0,0) {\subfloat{\includegraphics[scale=0.225]{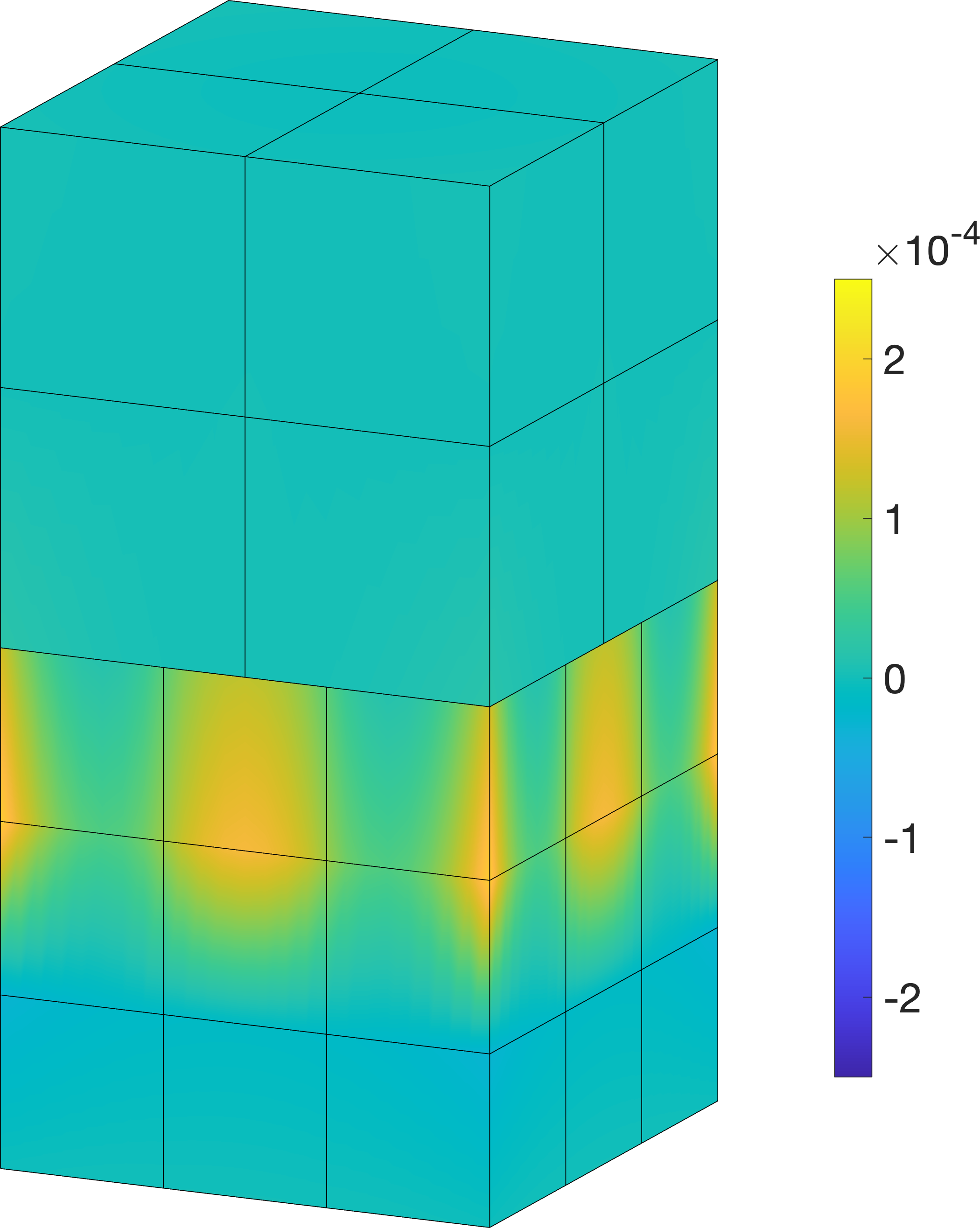}\label{fig:m1_ngps15}}};
        \node[font=\footnotesize, anchor=south] at (1.65, 6.215) {$n_\mathrm{gp}= 4(p+1)$ };
    \end{tikzpicture}~~
    \begin{tikzpicture}
        \node[anchor=south west, inner sep=0] (img3) at (0,0) {\subfloat{\includegraphics[scale=0.225]{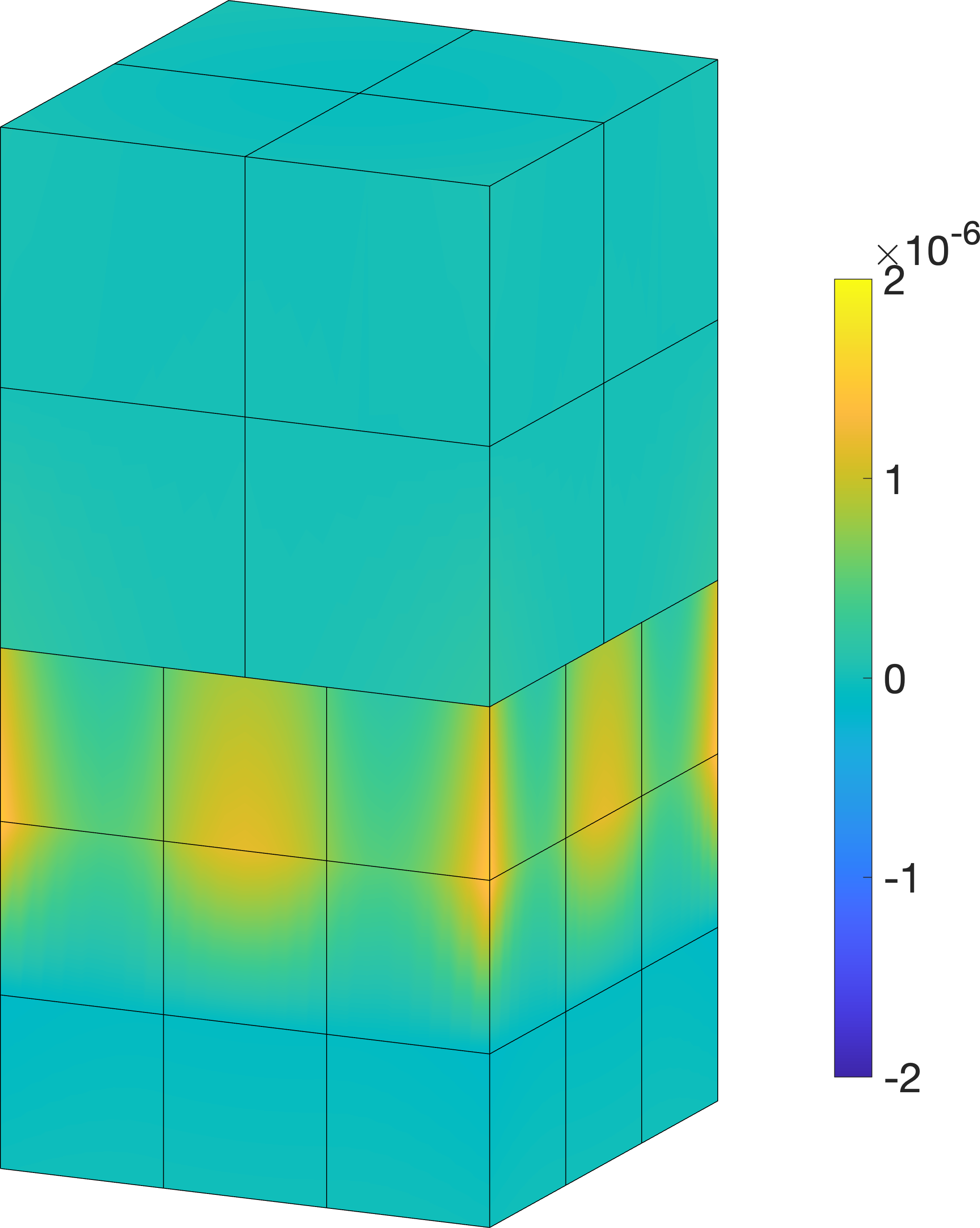}\label{fig:m1_ngps200}}};
        \node[font=\footnotesize, anchor=south] at (1.65, 6.215) { $n_\mathrm{gp}=200$ };
    \end{tikzpicture} \\
	\subfloat{\includegraphics[scale=0.225]{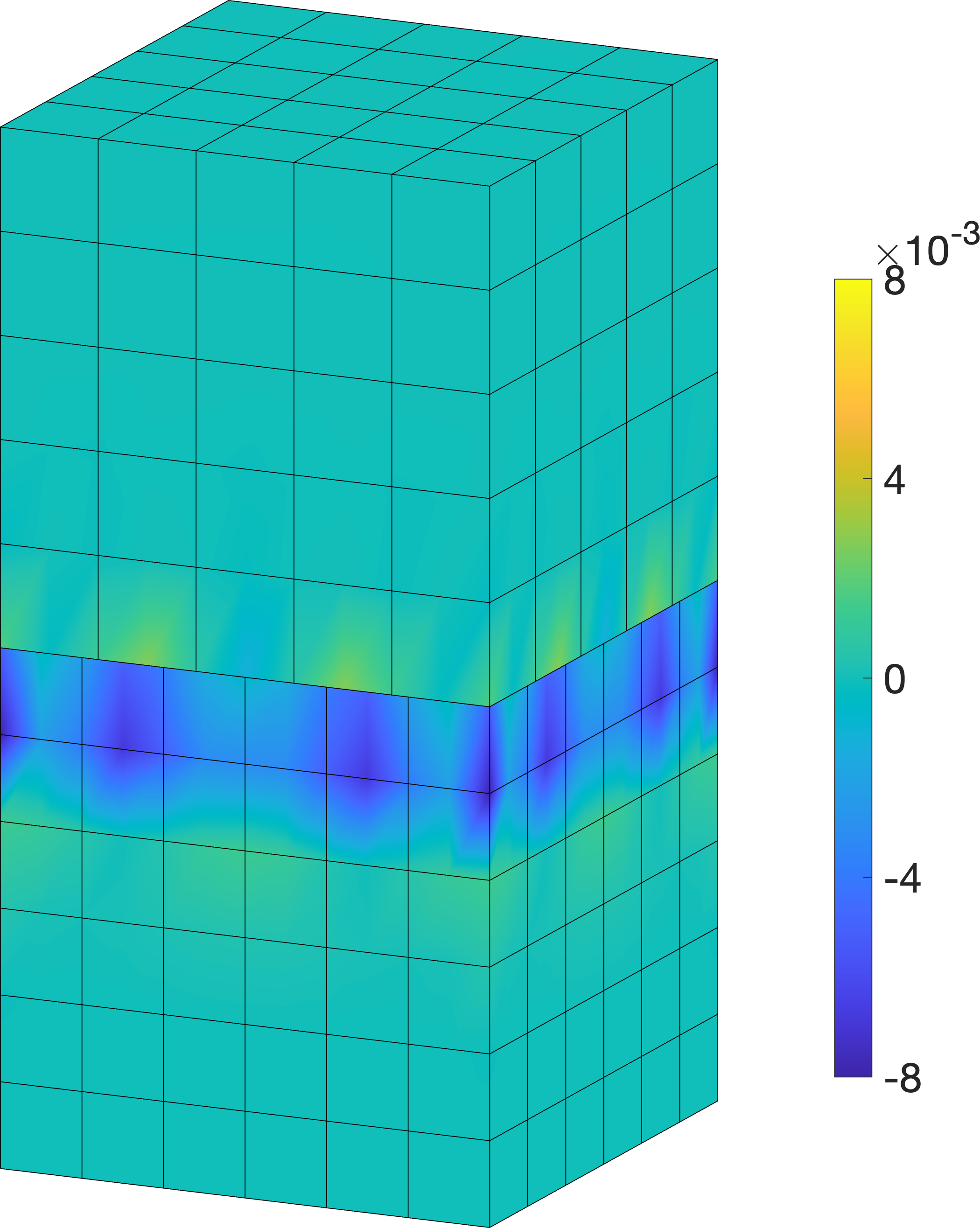}\label{fig:m2_ngps3}} ~~
    \subfloat{\includegraphics[scale=0.225]{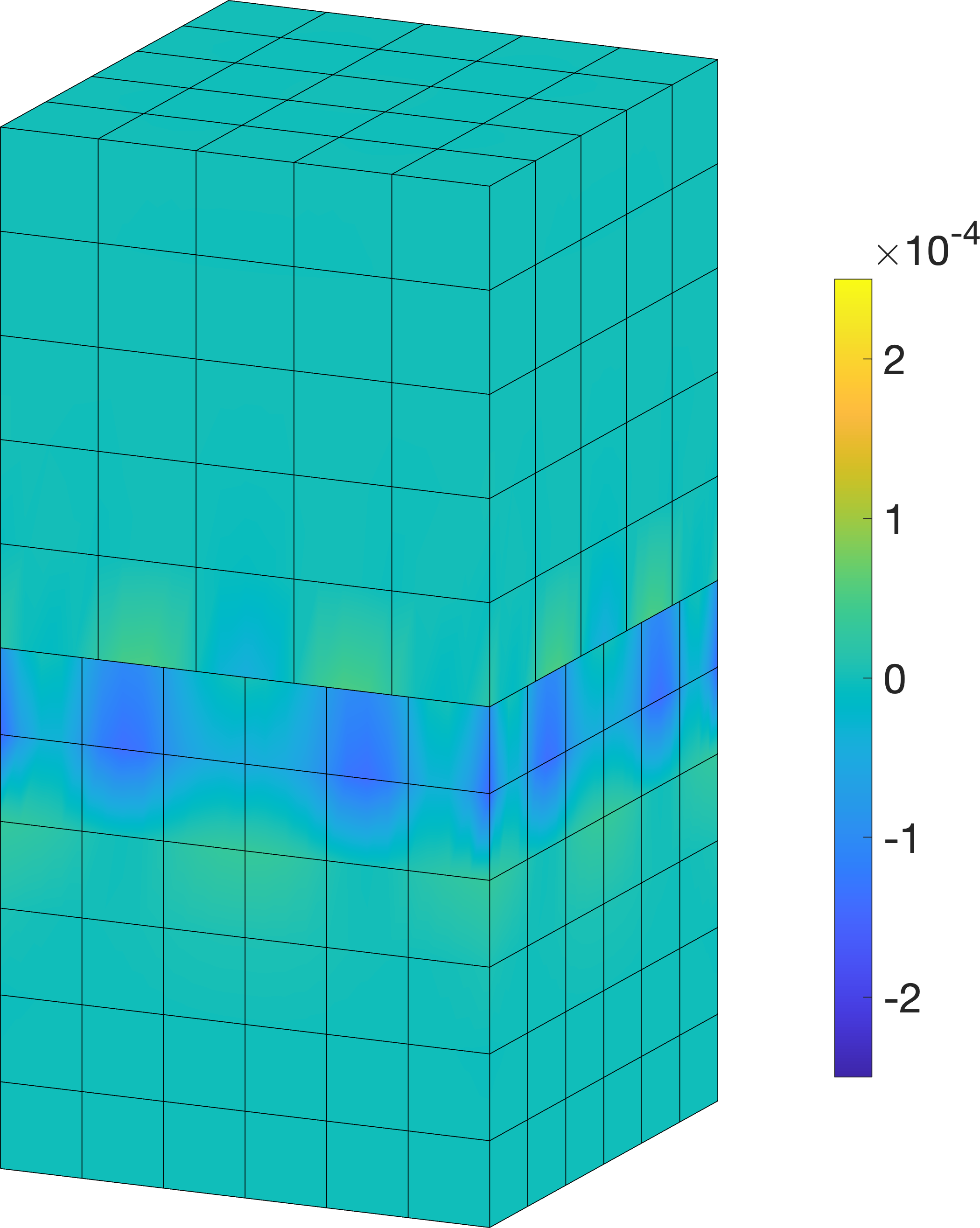}\label{fig:m2_ngps15}} ~~
    \subfloat{\includegraphics[scale=0.225]{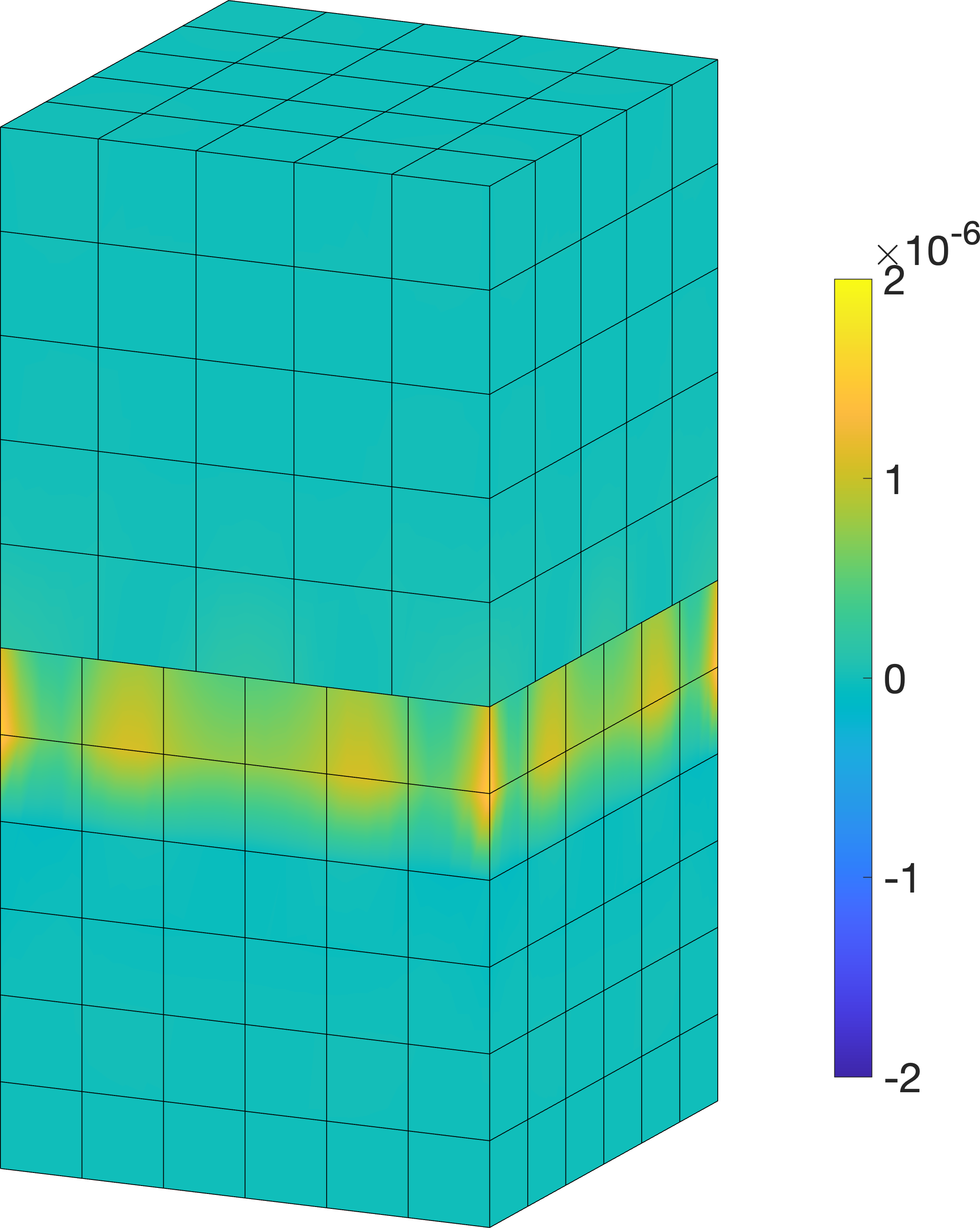}\label{fig:m2_ngps22}} \\
    \subfloat{\includegraphics[scale=0.225]{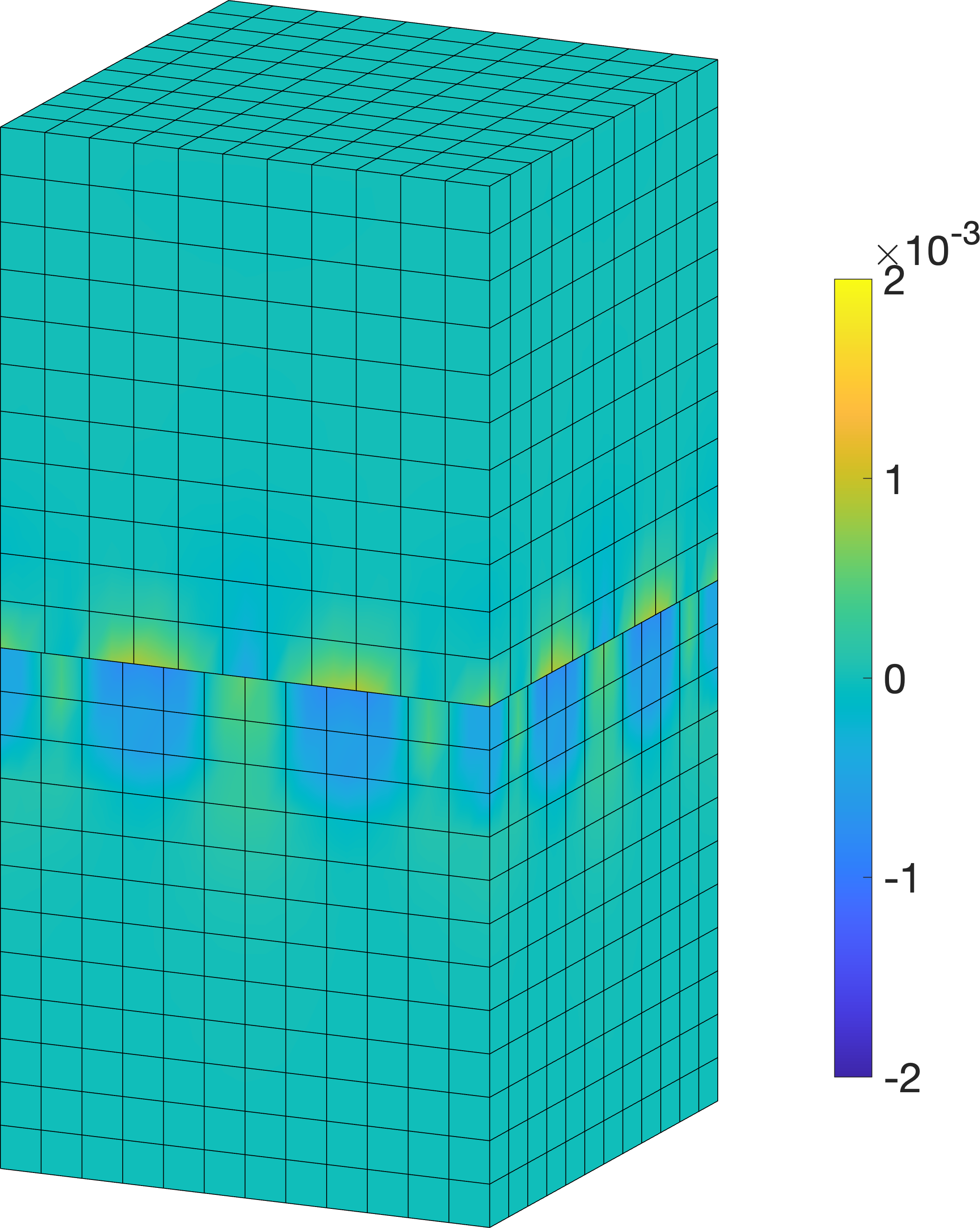}\label{fig:m4_ngps3}}~~
    \subfloat{\includegraphics[scale=0.225]{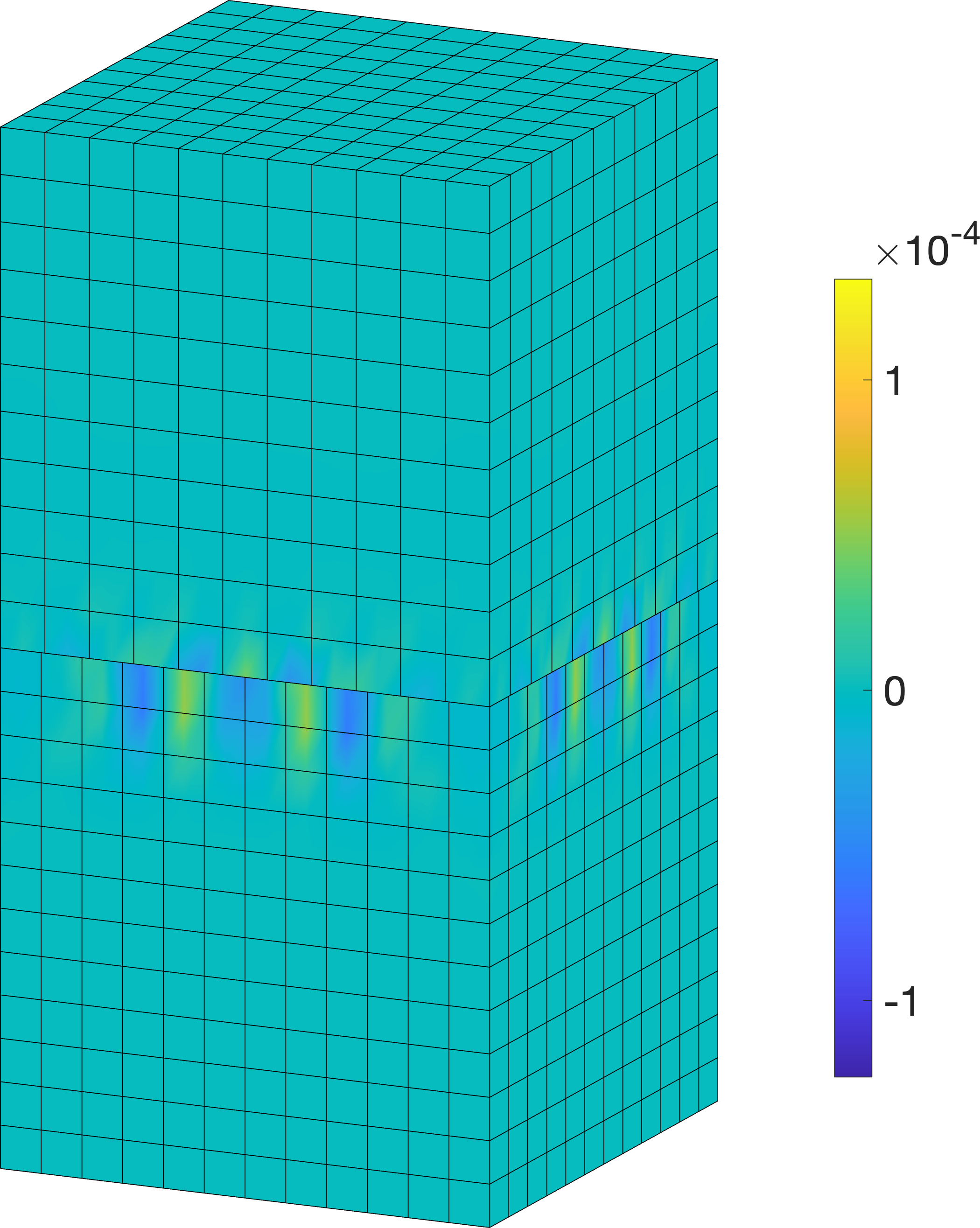}\label{fig:m4_ngps15}} ~~
    {\subfloat{\includegraphics[scale=0.225]{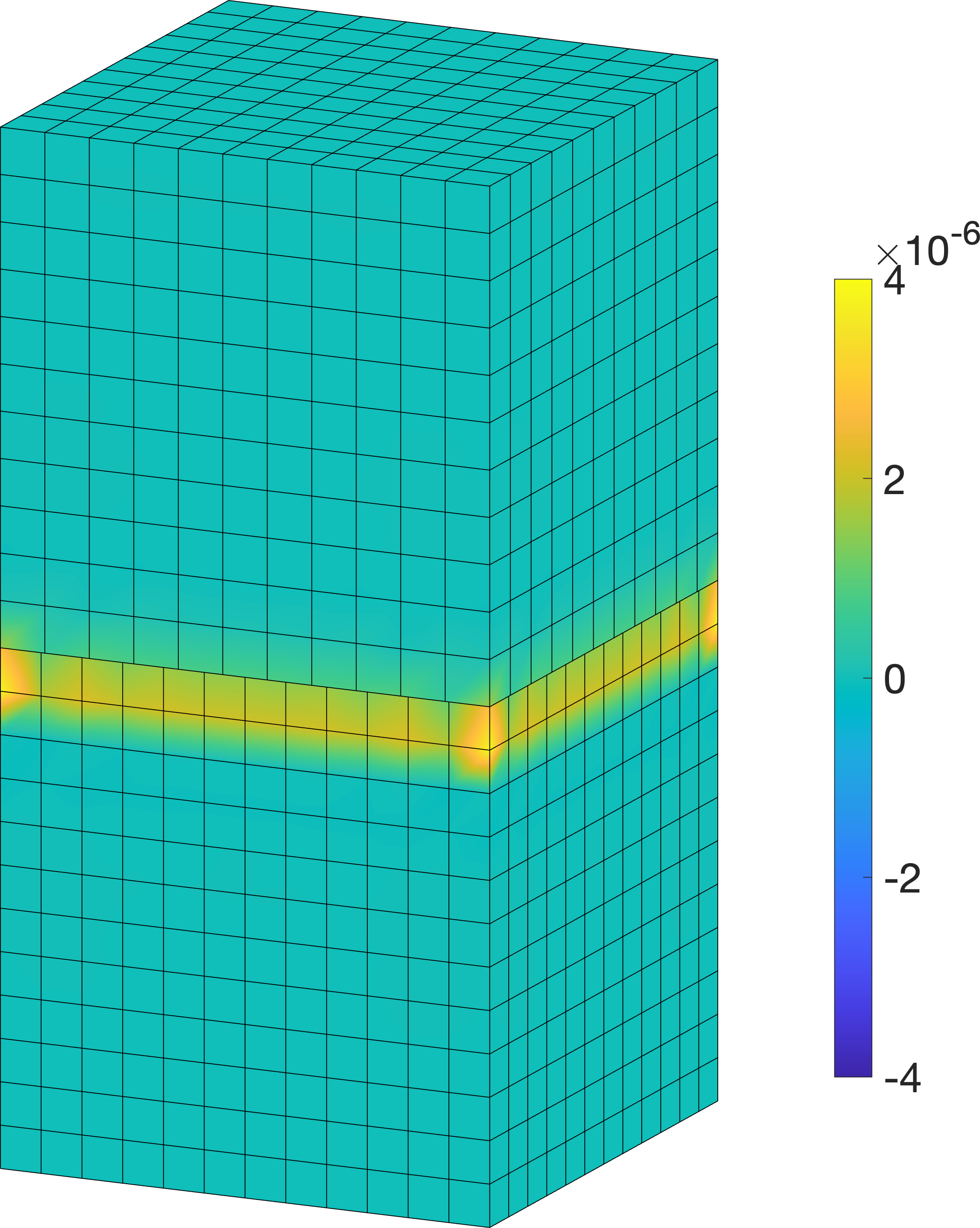}\label{fig:m4_ngps200}}};
    \caption{Patch test: Error in the vertical stress for GPTS contact algorithm over different meshes but at a fix number of quadrature points (column-wise), and with varying the number of quadrature points at a fixed mesh (row-wise).} \label{fig:PT_results}
\end{figure}

In this example, we examine the relative error in the vertical contact stress for the considered GPTS isogeometric contact algorithm while varying both the mesh resolution and the number of quadrature points at the interface. The setup involves a frictionless contact between two cubes of identical dimensions, $[L_x, L_y, L_z] = [1, 1, 1] $ mm with non-conforming meshes at their interface, as shown in Fig.~\ref{fig:PathTest_setup_m1}. A uniform contact pressure, $p_N = 1$ N/mm$^2$ is applied on the top surface of the upper cube, while the bottom surface of the lower cube is constrained against the vertical direction. Both bodies are modeled as isotropic, and linear elastic material is considered, with $E=1$ N/mm$^2$ and $\nu = 0.3$. The penalty parameter is set to $\epsilon_N = 100 E$. Three successively refined meshes, containing $ [3, 6, 12]$ and $ [2, 5, 11]$ number of quadratic NURBS elements N$_2$ in each direction for the lower and upper cubes, respectively, as illustrated in Fig.~\ref{fig:PT_results}. In all the cases, the upper body is always set as the slave body, while the lower one is the master body. For the evaluation of the contact integral, we use $n_{\mathrm{gp}} = (p+1)$, $4*(p+1)$, and $200$ Gauss-Legendre integration points in each direction in each element.

Figure~\ref{fig:PT_results} illustrates the error in the vertical contact stress over the different meshes and quadrature points. The results show that the error decreases only with increasing the number of quadrature points while remaining nearly unaffected by mesh resolution or interpolation order (not shown here) at a fixed number of quadrature points. Thus, the considered GPTS algorithm passes the contact patch test within the quadrature accuracy. Based on the patch test results, we set $ 4(p+1) $ as the default number of quadrature points for the evaluation of contact integral along each parametric direction of the $N_p$ discretized NURBS contact layer in the subsequent analysis to ensure numerical efficiency.

\subsection{Hertzian contact problem}
\begin{figure}[!ht]
	\centering
	\subfloat[]{\includegraphics[scale=0.32]{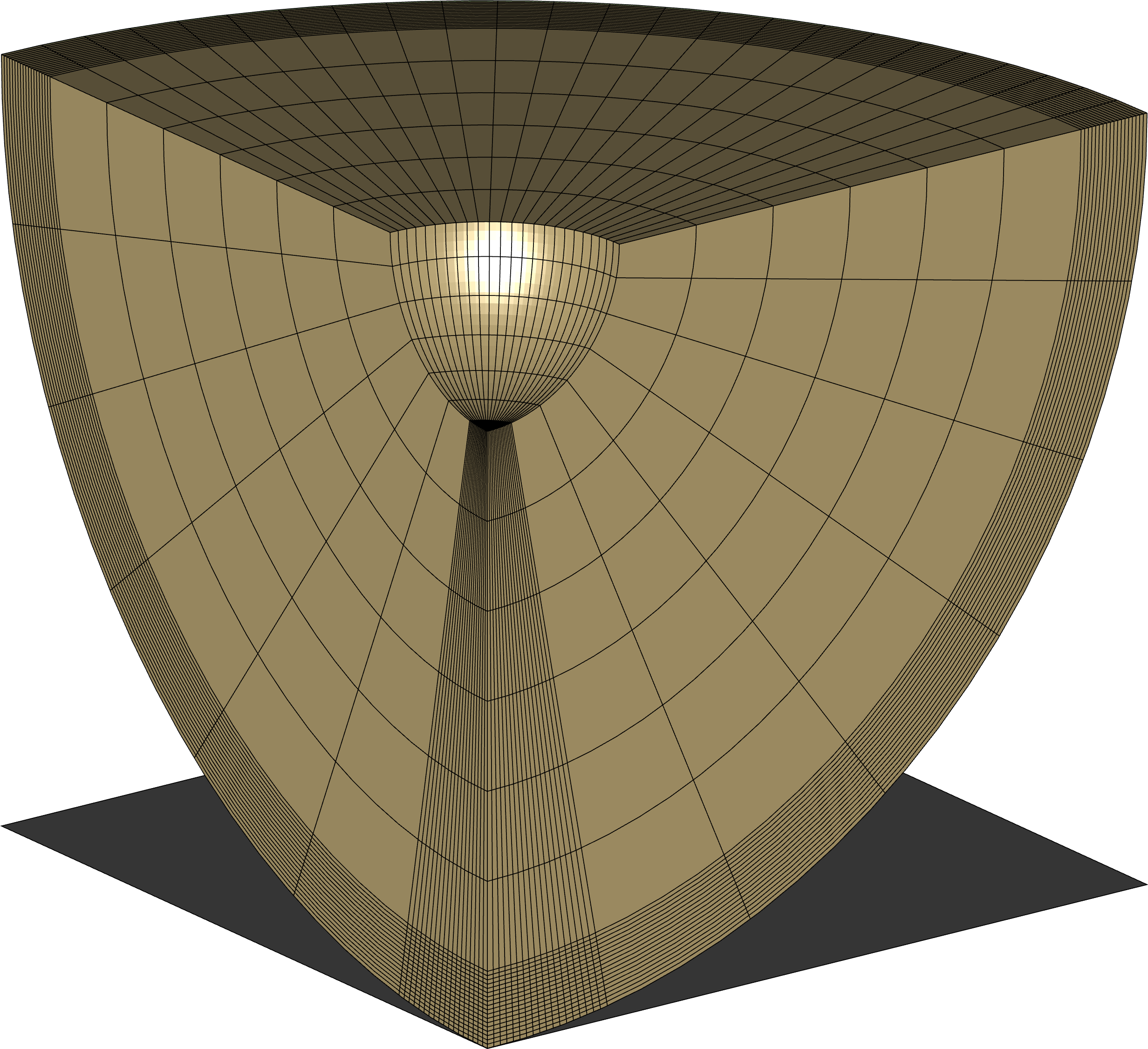}\label{fig:3DHertizian_setup}}~~
	\subfloat[]{\includegraphics[scale=0.12]{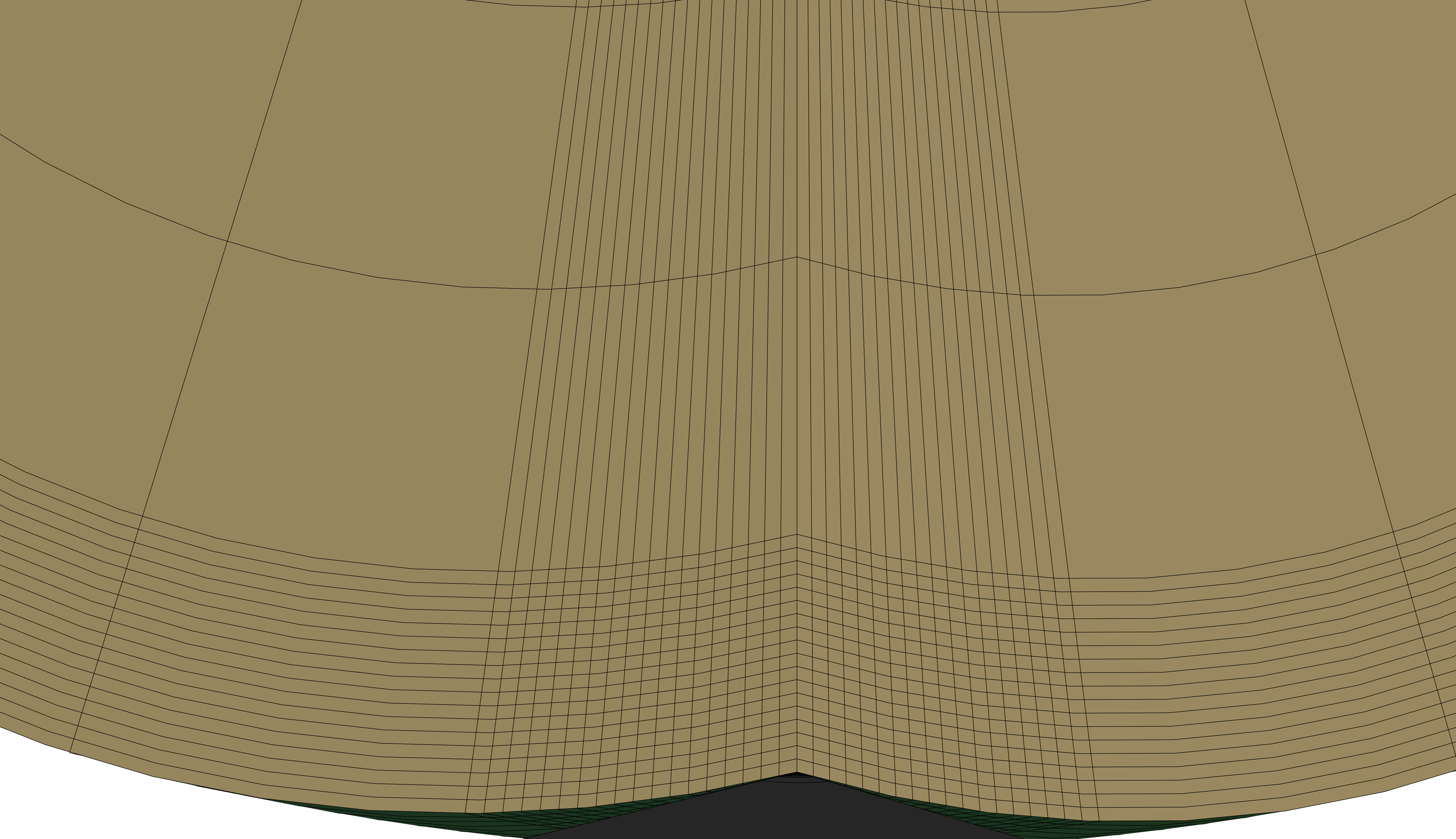}\label{fig:3DHertizian_mesh}} 
	\caption{Hertzian contact: (a) $1/8^{\mathrm{th}}$ of the NURBS constructed sphere with the redistributed coarsest mesh, and (b) enlarged view of the mesh in the contact zone.} \label{fig:3DHertzian_setup_mesh}
\end{figure}
The second examples considers the Hertzian contact problem between an elastic sphere and a rigid plane, as illustrated in Fig.~\ref{fig:3DHertizian_setup}. This example is used to evaluate the performance of the 3D varying order NURBS discretization method in capturing the contact pressure distribution at a fixed mesh, compared to the standard N$_2$ NURBS discretization. The geometry of the problem, boundary conditions, and material parameters are taken from De Lorenzis et al.~\cite{DeLorenzis2012}.

The sphere, having the outer radius $R_o=1$ mm, is made of a linear elastic material with Young's modulus $E=1$ N/mm$^2$ and Poisson's ratio $\nu = 0.3$. The penalty parameter $\epsilon_{\mathrm{N}} = 5\times 10^3$ is set as the default value. It is noted that the natural NURBS-based construction of the $1/8$ of the sphere includes an internal spherical surface with an inner radius $R_i = 0.8$ mm. Numerical experimentations verified that this internal radius does not influence the accuracy of the results. In the classical Hertz problem, a uniform normal pressure distribution is applied on the equatorial plane of the hemisphere. However, in this case, due to the presence of an internal spherical surface, a uniform vertical displacement is applied to its top surface, as in~\cite{Agrawal2020}. For the comparative convergence analysis, five progressive refined meshes,  driven by $ 12n \times 24$ number of elements along the angular and radial directions, where $n= \left[1, 2, 4, 8, 16 \right]$, denoted as $m_1$ to $m_5$. Numerical tests confirmed that the $24$ number of redistributed elements along the radial direction is sufficient to obtain accurate results for this example. For numerical efficiency purposes, the mesh is redistributed such that approximately $75 \% $ number of elements are concentrated within the $10\%$ of the total length of knot vectors defined in both directions, as illustrated for mesh $m_1$ in Fig.~\ref{fig:3DHertzian_setup_mesh}. 

It is worth noting that due to the axial symmetry of the problem, the contact pressure distribution can readily be obtained from a simplified 2D setup. However, this study primarily evaluates the performance of the VO NURBS discretization in capturing the contact pressure distribution in the 3D setting compared to the standard NURBS discretization. The same analysis was presented by De Lorenzis et al.~\cite{DeLorenzis2012}, which used the standard NURBS discretization with mortar contact elements. 

\begin{figure}[!ht]
	\centering
	\subfloat{\includegraphics[scale=0.26]{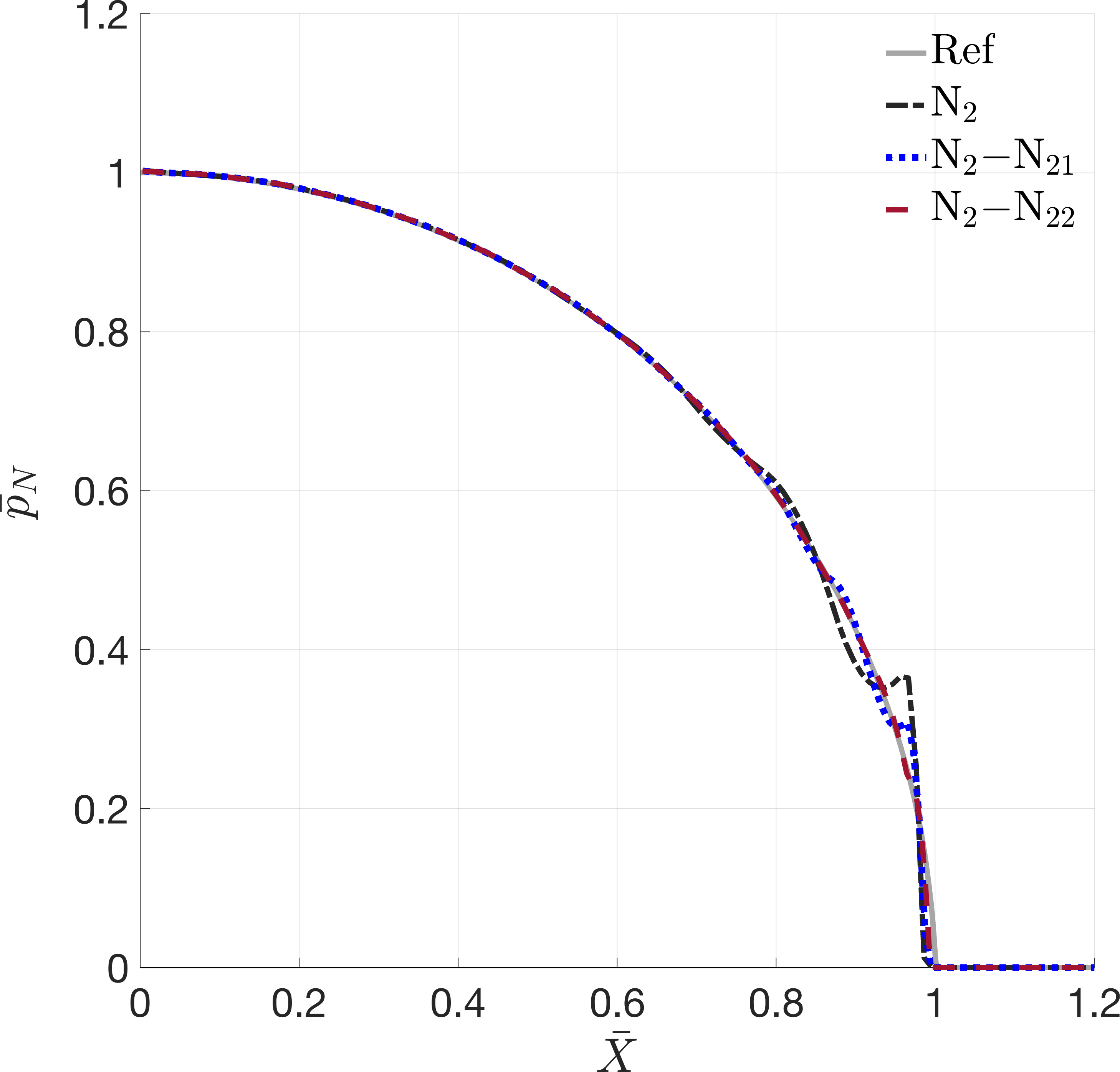}}
	\caption{Hertzian contact: Distribution of the contact pressure for different VO and standard NURBS discretization at mesh $m_3$. The result with N$_2$ at a fine mesh is used as a reference.} \label{fig:Hertzian_pressure}
\end{figure}
Figure~\ref{fig:Hertzian_pressure} qualitatively compares the dimensionless contact pressure distribution for VO-based N$_2$-N$_{2\cdot 1}$ and N$_2$-N$_{2 \cdot 2}$ with the standard N$_2$ discretization at an intermediate mesh $m_3$. As can be observed, at a fixed mesh, VO-based discretizations deliver superior quality results compared to the standard N$_2$ discretization. This is because the pressure curves are closer to the reference solution. This improvement stems from the increased number of DOFs available with VO discretization at the contact surface than N$_2$, as detailed in  Table~\ref{table:Hertzian_table}. A closer examination reveals a slight improvement in the quality of the result when the interpolation order of NURBS is increased at the contact surface within the VO setting, i.e., from N$_2$-N$_{2 \cdot 1}$ to N$_2$-N$_{2 \cdot 2}$.

\begin{table}[!h]
\begin{center}
\begin{tabular}{|c| c |c| c| c|}
\hline
\textbf{Discretization Type} & \textbf{Mesh} & \textbf{Interface} & \textbf{Bulk} & \textbf{Total} \\[1ex]
\hline
N$_2$ & & 588 & 14112 & 14700 \\
N$_2$-N$_{2\cdot 1}$ & m$_1$ & 2028 & 14112 & 16140 \\
N$_2$-N$_{2\cdot 2}$ & & 4332 & 14112 & 18444 \\
\hline
N$_2$ & & 2028 & 48672 & 50700 \\
N$_2$-N$_{2\cdot 1}$ & m$_2$ & 7500 & 48672 & 56172 \\
N$_2$-N$_{2\cdot 2}$ & & 16428 & 48672 & 65100 \\
\hline
N$_2$ & & 7500 & 180000 & 187500 \\
N$_2$-N$_{2\cdot 1}$ & m$_3$ & 28812 & 180000 & 208812 \\
N$_2$-N$_{2\cdot 2}$ & & 63948 & 180000 & 243948 \\
\hline
N$_2$ & & 28812 & 691488 & 720300 \\
N$_2$-N$_{2\cdot 1}$ & m$_4$ & 112908 & 691488 & 804396 \\
N$_2$-N$_{2\cdot 2}$ & & 252300 & 691488 & 943788 \\
\hline
N$_2$ &  & 112908 & 2709792 & 2822700  \\
N$_2$-N$_{2\cdot 1}$ & m$_5$ & 446988 & 2709792 & 3156780 \\
N$_2$-N$_{2\cdot 2}$ & & 1002252 & 2709792 & 3712044 \\
\hline
\end{tabular}
\caption{Hertzian contact: DOF distribution across the contact interface and bulk with the different VO and standard NURBS discretizations.}
\label{table:Hertzian_table}
\end{center}
\end{table}

\begin{figure}[!b]
	\centering
	\subfloat[]{\includegraphics[scale=0.25]{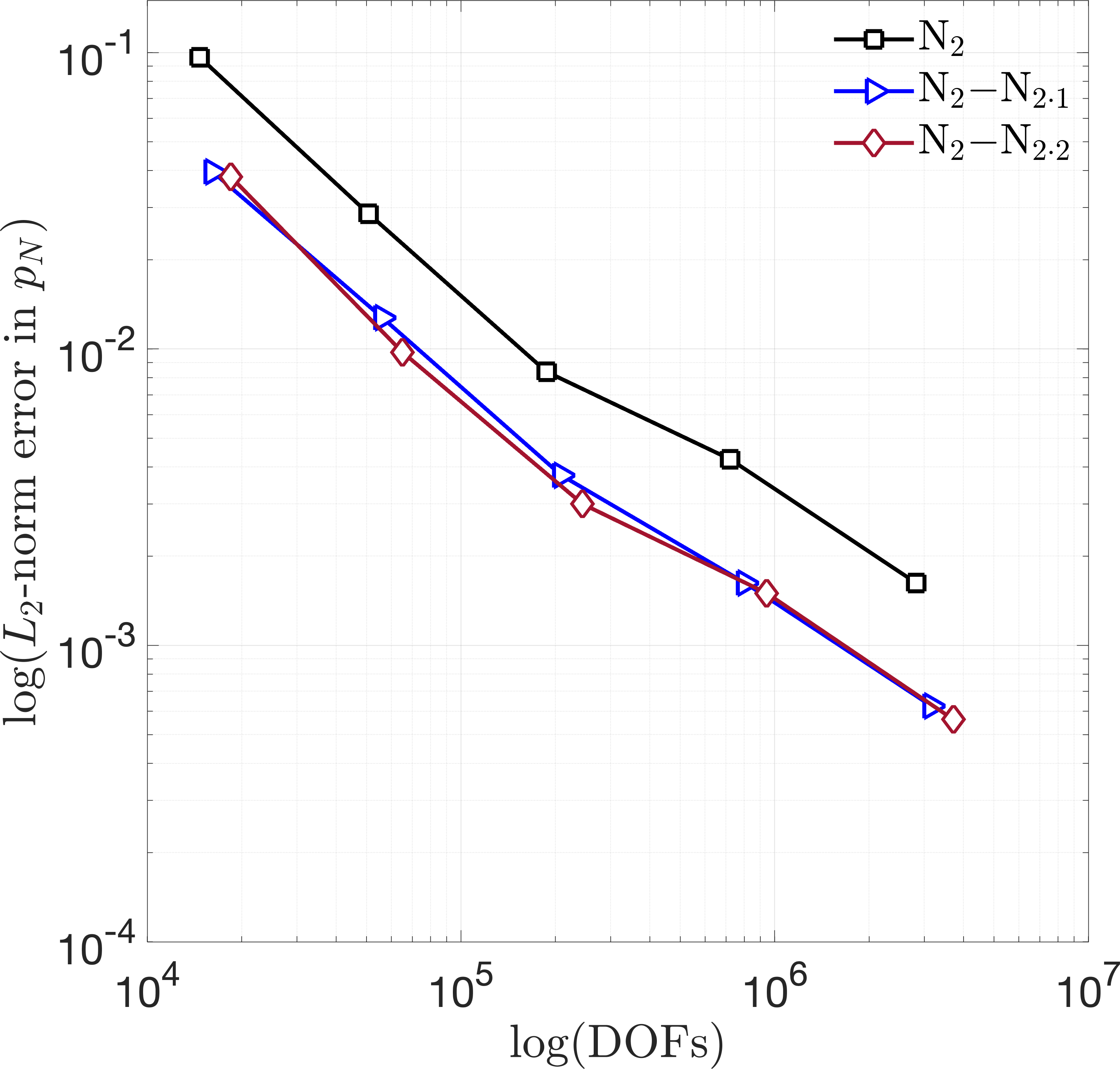}} ~~~
    \subfloat[]{\includegraphics[scale=0.25]{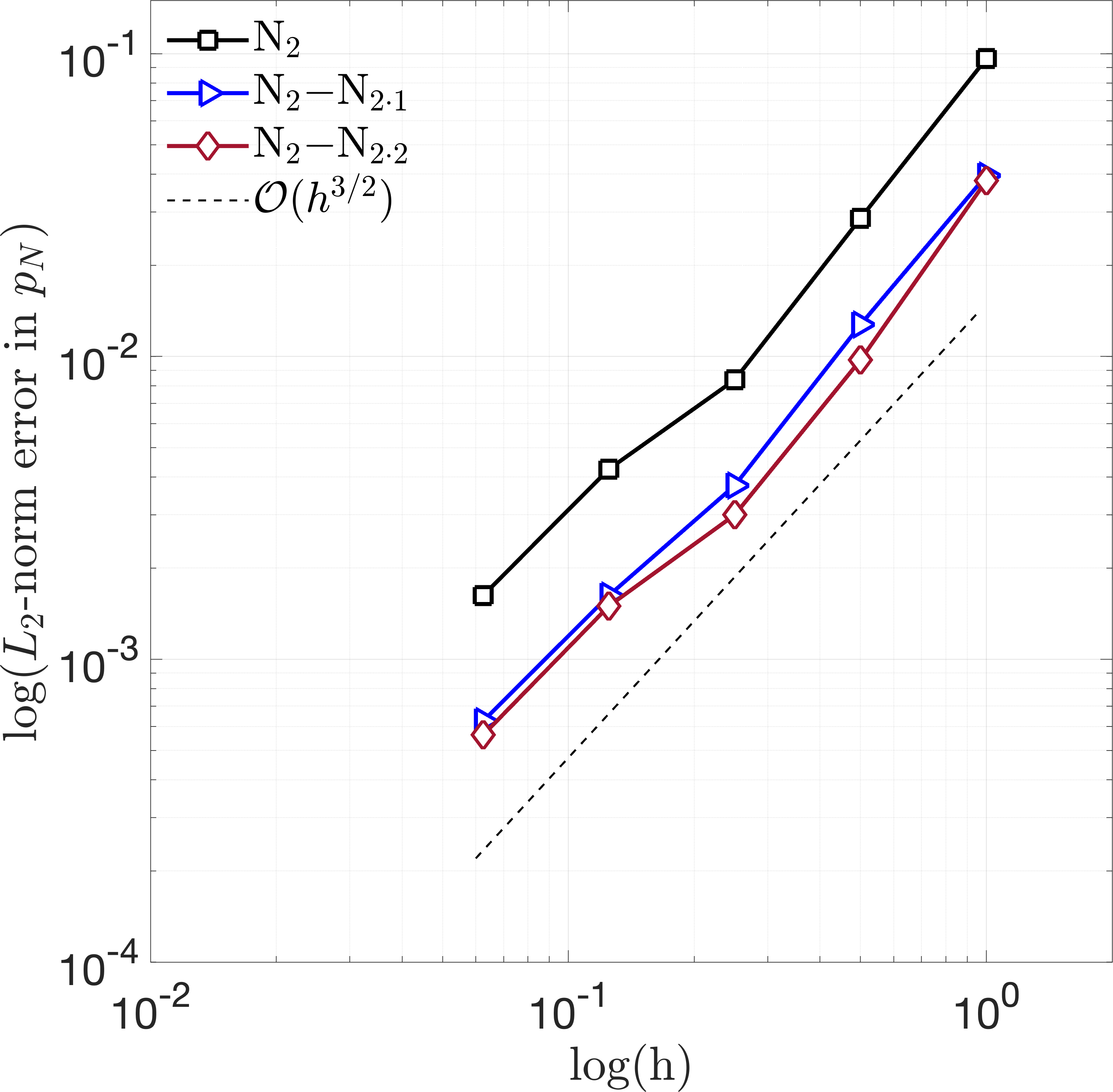} \label{fig:3DHertz_conv_eleSize}}
	\caption{Hertzian contact: convergence of the $L_2$-norm based error in the normalized contact pressure distribution with (a) total DOFs and (b) element size for both the standard and VO-based NURBS discretizations.}\label{fig:3DHertzian_mconv}
\end{figure}

Next, fig.~\ref{fig:3DHertzian_mconv} illustrates the convergence of the $L_2$-norm error in normalized contact pressure for these discretizations. The $L_2$-norm error at each mesh level is computed using $ ||e_{p_N}||_{{L}_2} = \sqrt{ \int_{\Gamma_c} \left[\bar{p}_{N_{\mathrm{ref}}} - \bar{p}_{N}\right]^2\,\textrm{d}\Gamma}\, $, where the reference solution is obtained at sufficiently fine mesh, containing approx. $6.4$ million of DOFs with N$_2$ discretization.

For unilateral contact problems, the spatial convergence rate of the displacement field is theoretically limited to $\mathcal{O}(h^{3/2})$ in the energy norm upon uniform mesh refinement when using second-order contact elements. This limitation arises due to the reduced regularity of the contact solution, as discussed in detail in Popp et al.~\cite{Popp2008} and Seitz et al.~\cite{SEITZ2016}. In the current work, as mentioned earlier, the loading on the top surface of the sphere is prescribed via displacement field rather than pressure due to the presence of an internal spherical surface. Consequently, we analyze the convergence of the $L_2$-norm of the contact pressure error, which differs from displacement-based convergence estimates such as $H^1$-norm. Our findings indicate that the convergence rate of this error is slightly below than $\mathcal{O}(h^{3/2})$ for both the standard and VO-based NURBS discretizations. This deviation may stem from the large mesh size associated with the numerical reference solution, i.e. ${p}_{N_{\mathrm{ref}}}$, which is limited by our computational resources.

From fig.~\ref{fig:3DHertzian_mconv}, it is evident that the error curve with the VO discretization consistently lies below that of the standard N$_2$ discretization, demonstrating improved accuracy at each mesh level. Notably, the accuracy with N$_2$-N$_{2\cdot 1}$ at mesh $m_2$, $m_3$, and $m_4$ is comparable to that with N$_2$ at next finer meshes, i.e. $m_3$, $m_4$, and $m_5$, respectively. This indicates that VO-based N$_2$-N$_{2\cdot 1}$ achieves similar accuracy while requiring approx. $3.5$ times fewer DOFs than N$_2$, see  Table~\ref{table:Hertzian_table} for DOF details. This significant reduction in DOFs translates to a considerable gain in computational efficiency.

Additionally, the error values obtained with N$_2$-N$_{2 \cdot 2}$ are only marginally better that those with N$_2$-N$_{2 \cdot 1}$. This indicates that increasing the number of DOFs at the contact surface, as in N$_2$-N$_{2 \cdot 2}$, may not be useful, as it provides only marginal improvement in the accuracy, making the use of additional DOFs than N$_2$-N$_{2 \cdot 1}$ difficult to justify.

\subsection{Frictional Ironing with finite deformations}
\subsubsection{Setup}
\begin{figure}[!hb]
\centering
\subfloat{\includegraphics[scale=0.4]{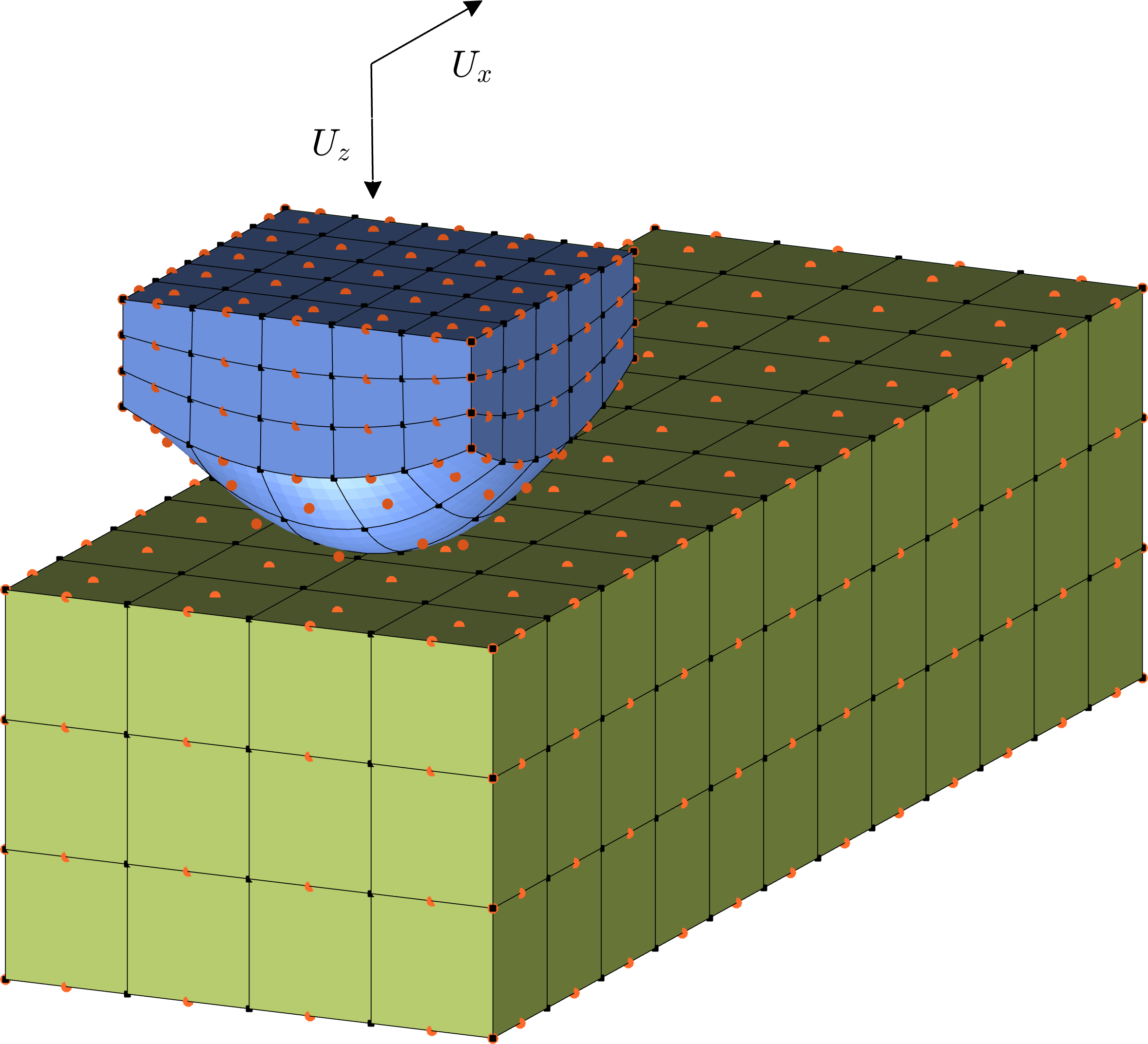}}
\caption{Frictional ironing: Setup and the coarsest mesh considered for the indentor and slab. The control points and knot entries are shown with red dots and black squares for N$_2$ with mesh m$_1$.}\label{fig:ironing_setup}
\end{figure}
This example considers a 3D frictional sliding contact between an indentor-slab pair undergoing large deformations. The setup of the problem, which is adapted from Temizer et al.~\cite{Temizer2012}, is shown in Fig.~\ref{fig:ironing_setup}. We simulate this example in two overall stages. In the first one, the elastic indentor having a spherical-shaped contact surface and geometrical dimension $[20, 20, 16]$ mm is pressed against an elastic slab of size $ [L_x, L_y, L_z] = [80, 28, 21]$ mm by applying the vertical displacement $U_z = -7$ mm on its top surface with $70$ uniform load steps. In the second stage, the indentor is dragged across the slab by applying horizontal displacement $U_x = 54$ mm on its top surface with uniform $540$ load steps. Note that during the dragging process, the vertical displacement is kept constant. To avoid the physical fold-in effect at the leading contact surface during the dragging process, Coulomb's friction coefficient $\mu = 0.1$ is used as in~\cite{Temizer2012}. The foundation of the slab is fixed in all the displacement directions. An isotropic Neo-Hookean material behavior is considered for both bodies, where Young's modulus is $E = 100$ N/mm$^2$ for the indentor and $E = 1$ N/mm$^2$ for the slab, and the same Poisson's ratio $\nu = 0.3$ for both. The penalty parameters are taken as $\epsilon_N = \epsilon_T = 100$ N/mm and are scaled with the minimum element size. The coarsest mesh used for the indentor and slab is shown in Fig.~\ref{fig:ironing_setup}. Other nested meshes that are obtained through the uniform knot insertion and considered for the convergence analysis are listed in Table~\ref{table:ironing_element}. 
\begin{table}[!ht]
\begin{center}
\begin{tabular}{|c c c| c c c | c c c|}
\hline
&\textbf{Mesh} &&  \multicolumn{6}{c|}{\textbf{Elements}} \\ [0.4ex]
\cline{4-9}		
&& && {Indentor} &&& {Slab} &\\ 
\hline
&m$ _1 $ &&&  $ 5 \times 5 \times 3 $ &&&  $ 12 \times 4 \times 3 $  &\\
&m$ _2 $ &&&  $ 10 \times 10 \times 4 $ &&&  $ 24 \times 8 \times 6 $  &\\
&m$ _3 $ &&&  $ 20 \times 20 \times 5 $ &&&  $ 48 \times 16 \times 12 $ &\\
\hline
\end{tabular} \caption{Frictional ironing: Elements with different meshes.} \label{table:ironing_element}
\end{center}
\end{table}

\begin{figure}[!p]
	\centering
	\subfloat[]{\includegraphics[scale=0.235]{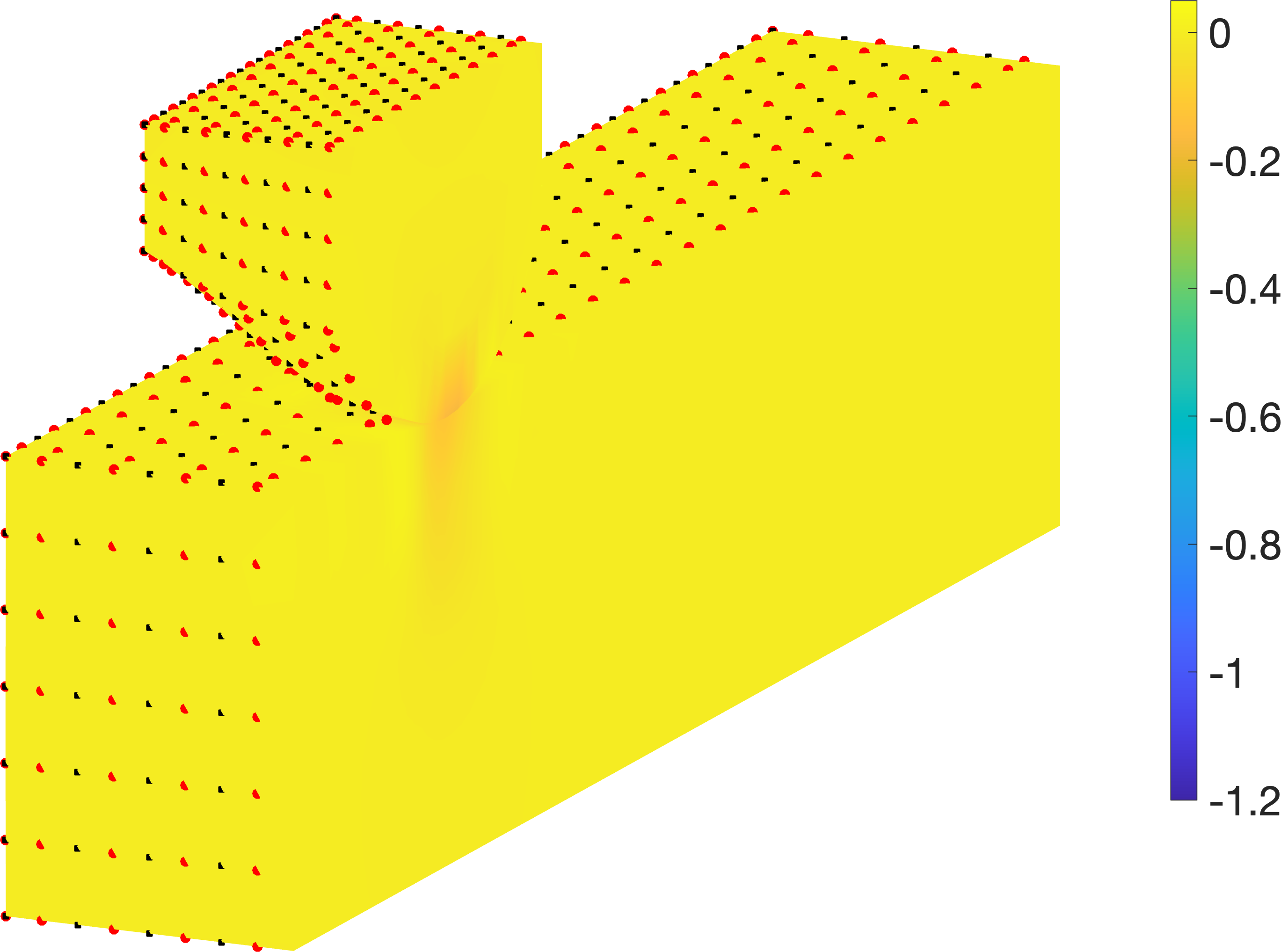}\label{fig:comp_s33_a}} ~~~~
    \subfloat[]{\includegraphics[scale=0.235]{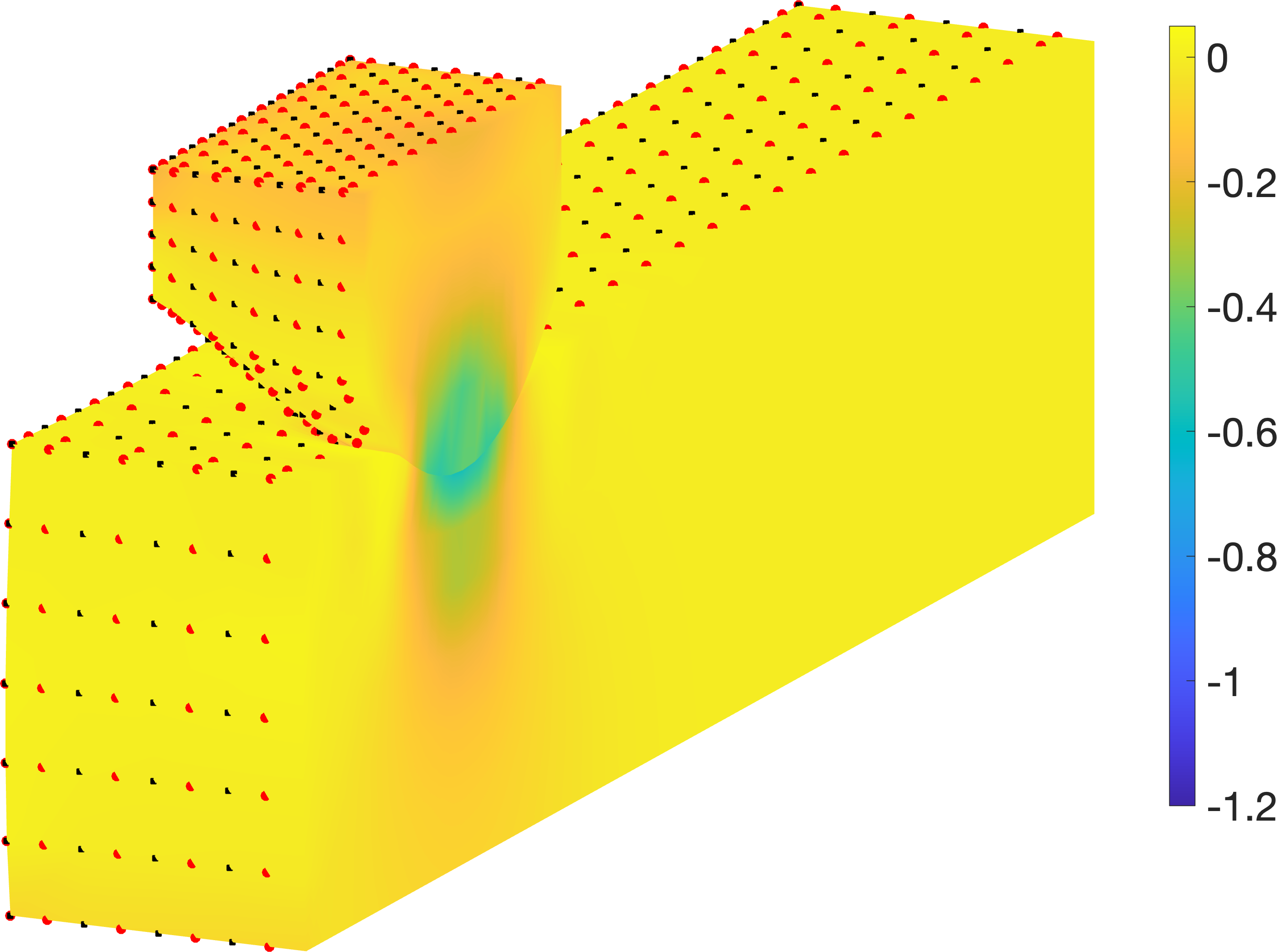}\label{fig:comp_s33_b}} \\
    \subfloat[]{\includegraphics[scale=0.235]{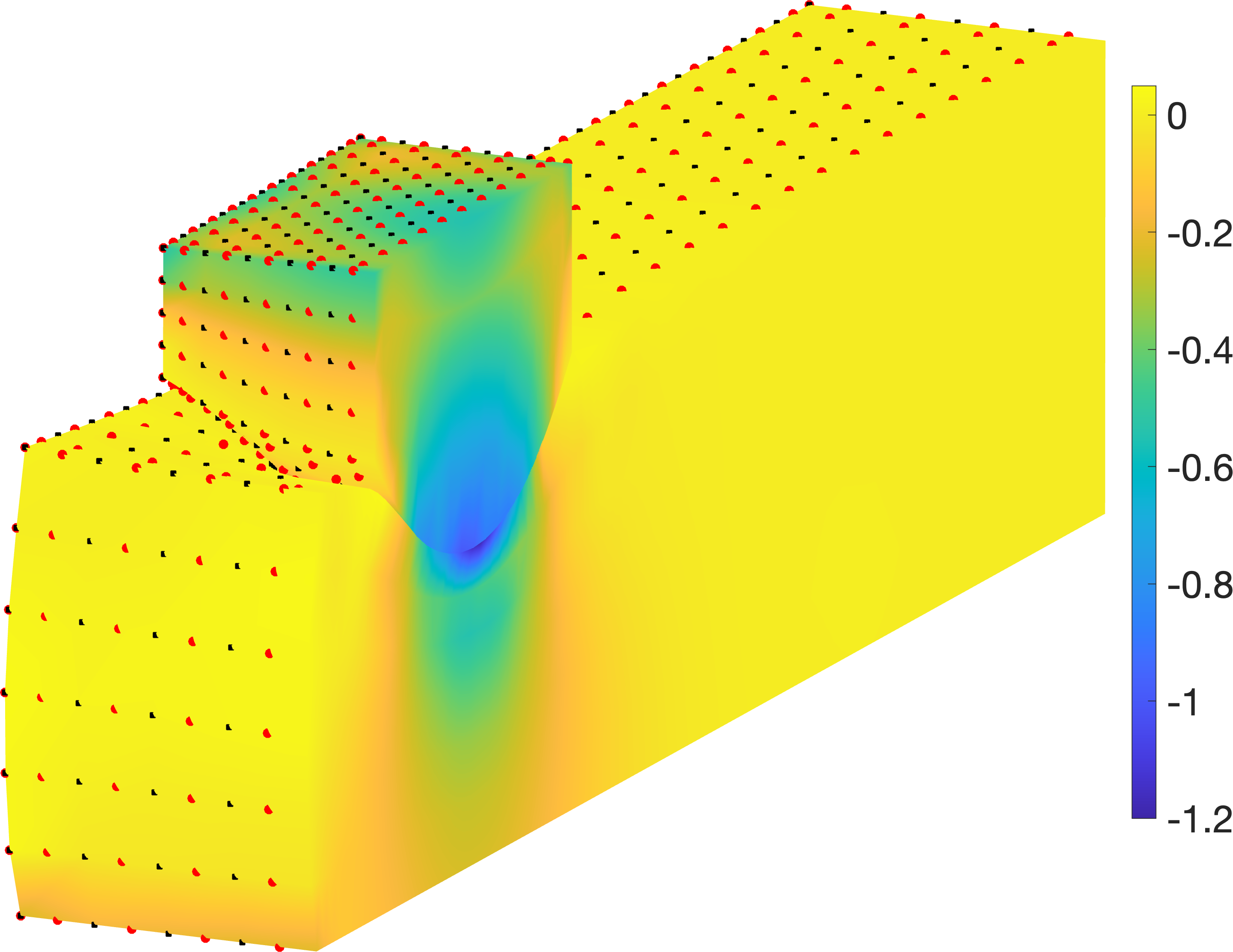}\label{fig:comp_s33_c}} ~~~~
    \subfloat[]{\includegraphics[scale=0.235]{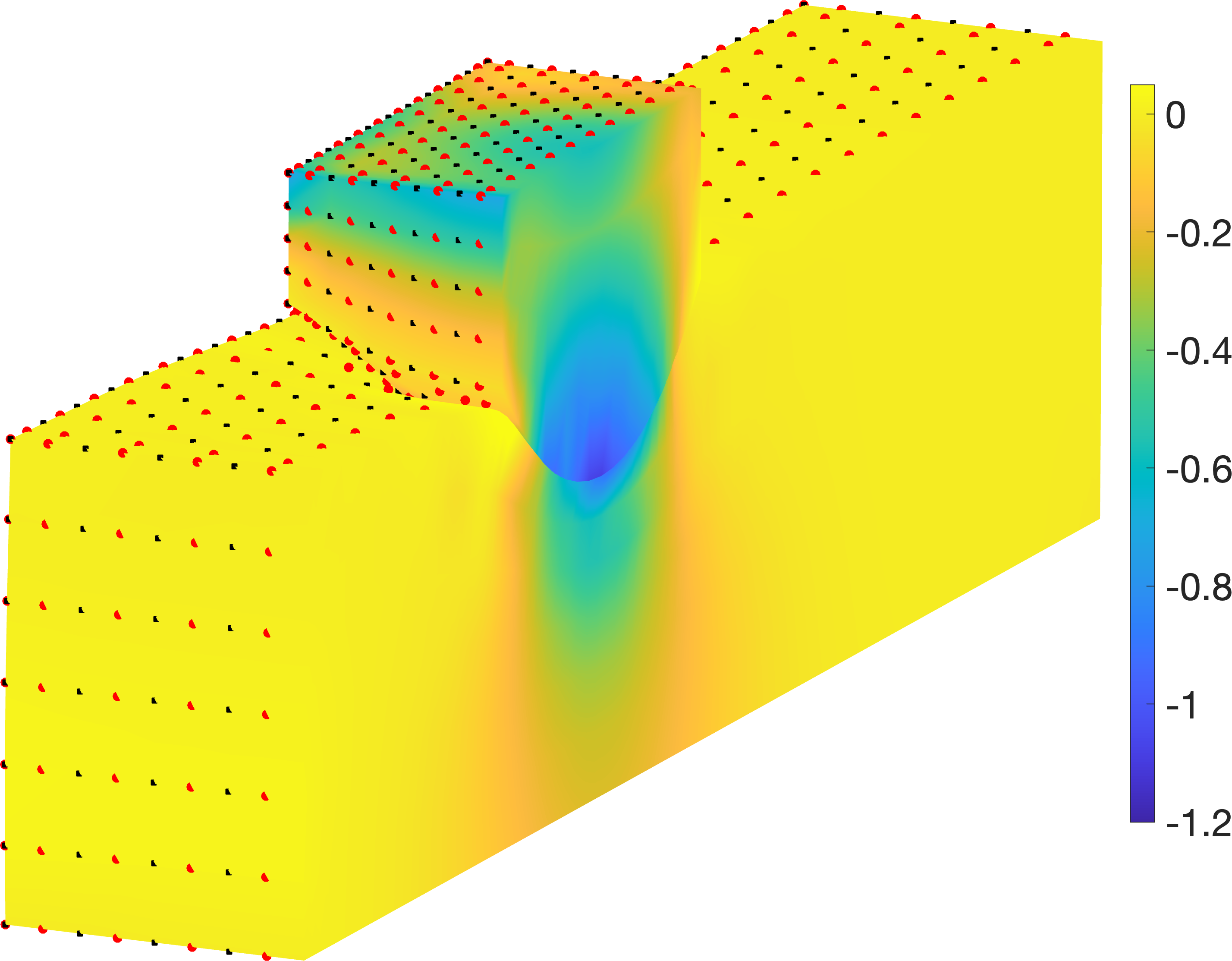}\label{fig:slid_s33_a}} \\
    \subfloat[]{\includegraphics[scale=0.235]{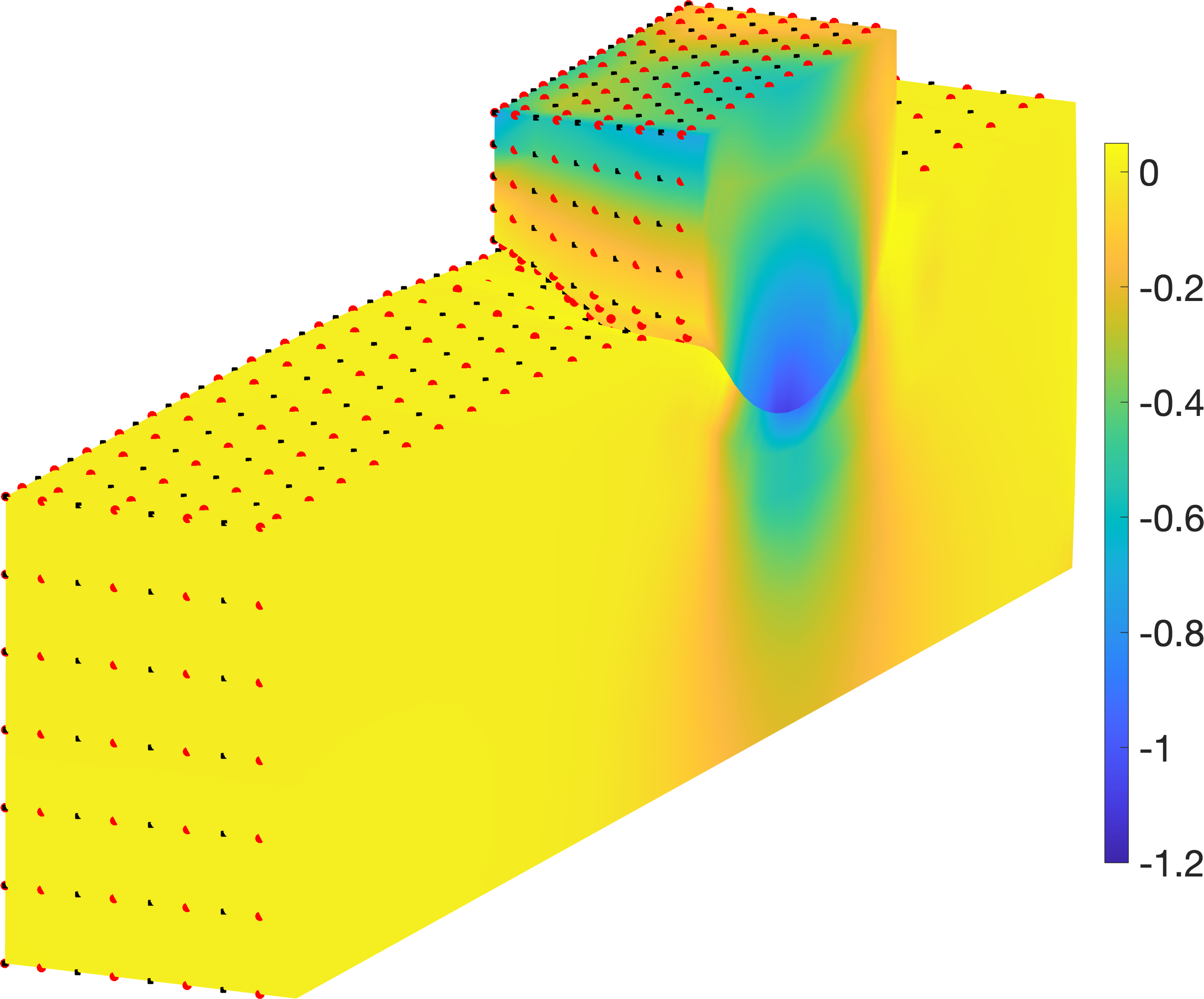}\label{fig:slid_s33_b}} ~~~~
    \subfloat[]{\includegraphics[scale=0.235]{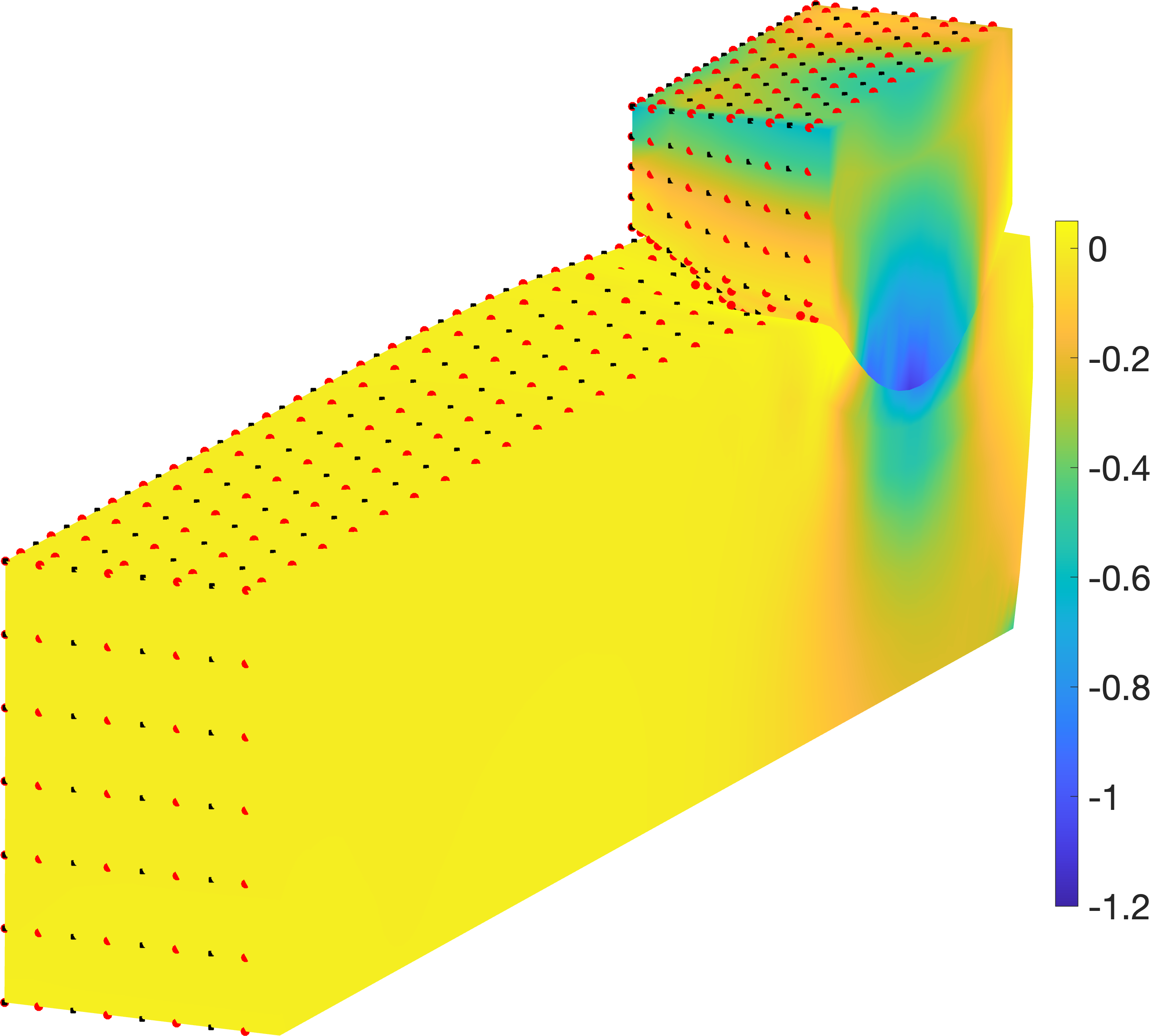}\label{fig:slid_s33_c}}
	\caption{Frictional ironing: Deformed configuration of the setup at the vertical indentation: (a) 0.05 mm, (b) 3.5 mm, and (c) 7 mm during the compression stage, and at (d) 18 mm, (e) 36 mm, and (f) 54 mm during the dragging stage. The color shows the distribution of stress $ \sigma_{33} $. The control points and knot entries with N$_2$ at mesh m$_2$ are shown with red dots and black squares.} \label{fig:ironing_deformed}
\end{figure}

The deformed configurations of the setup at different simulation instances are shown in Fig.~\ref{fig:ironing_deformed} with quadratic NURBS discretizations at mesh m$_2$. Here, the half-view of the setup is shown to provide better visuals of the contact surface deformation and the stress developed within the interior of the bulk domain.

\subsubsection{Performance at a fixed mesh}
\begin{figure}[!hb]
	\centering
	\subfloat{\includegraphics[scale=0.345]{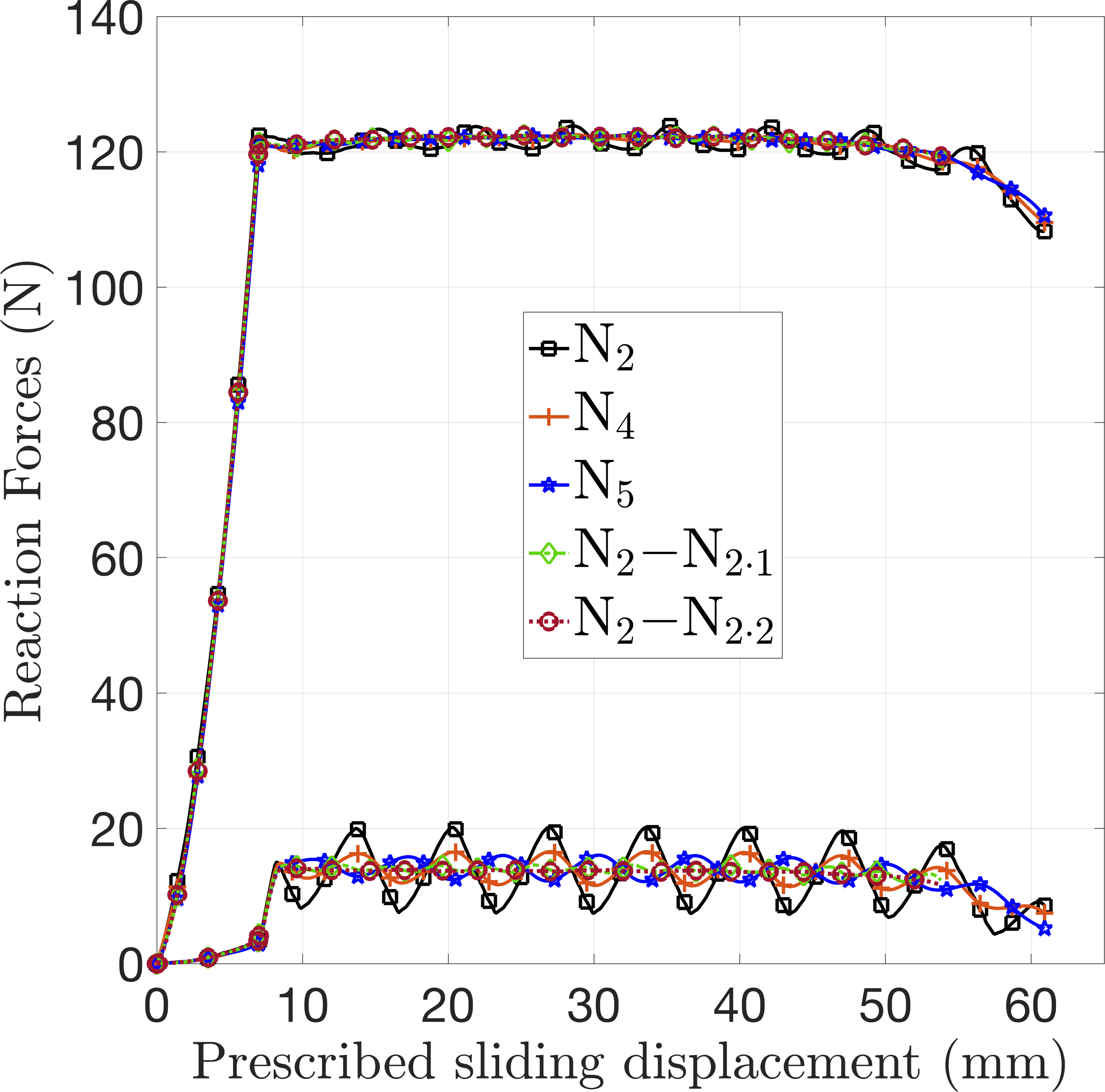}} 
    \caption{Frictional ironing: Evolution of total vertical and frictional contact forces over prescribed displacement with different NURBS discretizations.} \label{fig:Ironing_full_m1}
\end{figure}

First, we show the performance of the proposed method at the coarsest mesh level m$_1$ for the considered example. Figure~\ref{fig:Ironing_full_m1} illustrates the evolution of the total normal and frictional contact forces, computed at the top surface of the indentor as a function of prescribed displacement with different standard and VO NURBS-based discretizations. The results indicate that the evolution of the contact forces appears similar to that of different NURBS discretizations. Even though the contact surface is discretized using the $C^1$-continuous NURBS, the oscillation error in contact forces is noticeable, particularly during the dragging phase. The origin of these oscillations lies in the inadequate discretization of the contact layers, e.g. when using a coarse mesh and/or lower interpolation order of basis functions. This limitation hinders the contact layer from adequately conforming to the finite deformations of the bodies involved. However, at a fixed mesh, IGA provides the possibility to enhance the conforming ability of contact layers through the use of higher-order NURBS functions. However, using a fine mesh or higher-order NURBS in the region away from the contact surface may not be desirable as a substantial computational cost is associated with bulk computations. 
\begin{figure}[!hb]
	\centering
	\subfloat[]{\includegraphics[scale=0.245]{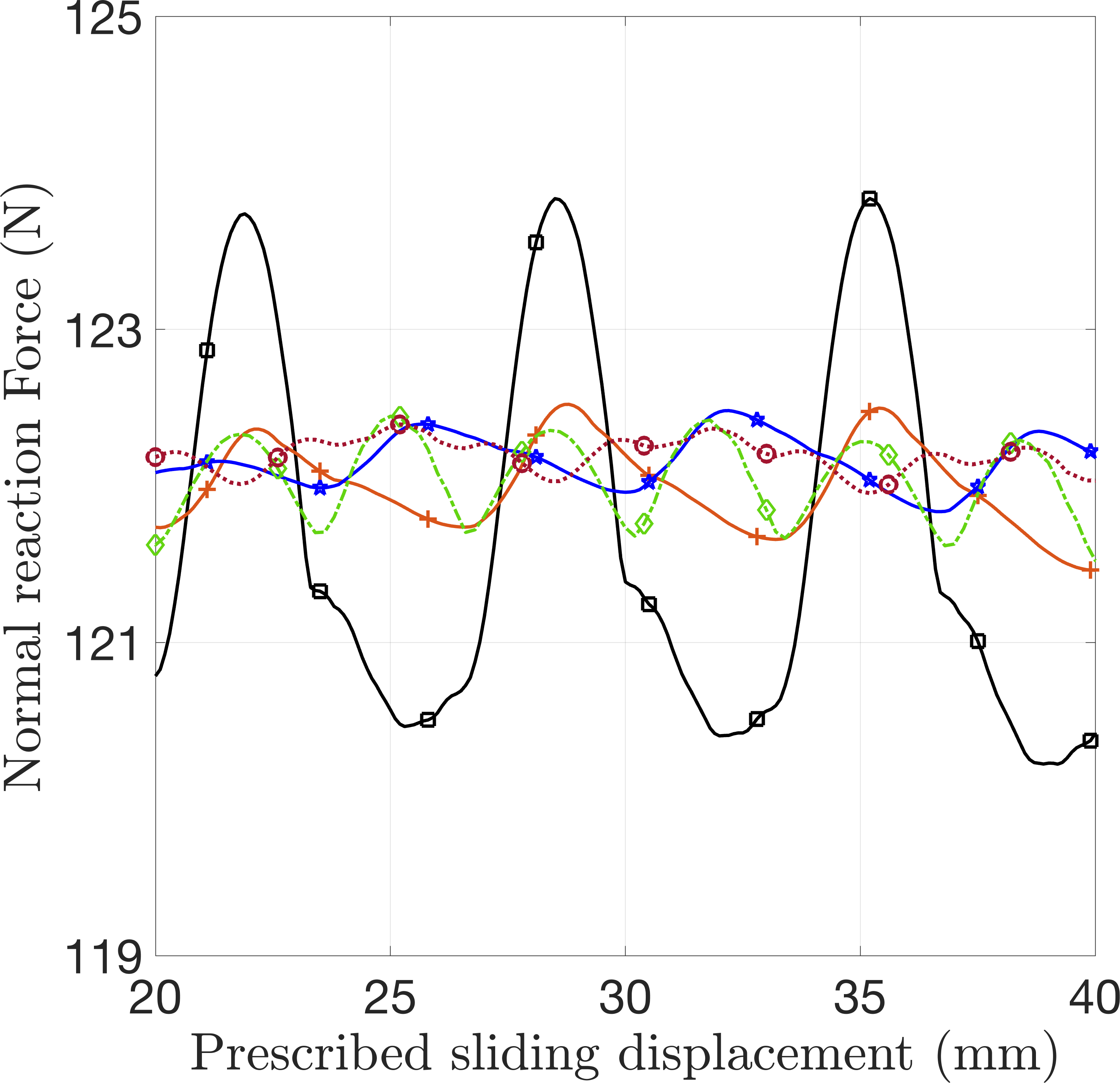}\label{fig:Ironing_Normal_m1}} ~~~
    \subfloat[]{\includegraphics[scale=0.245]{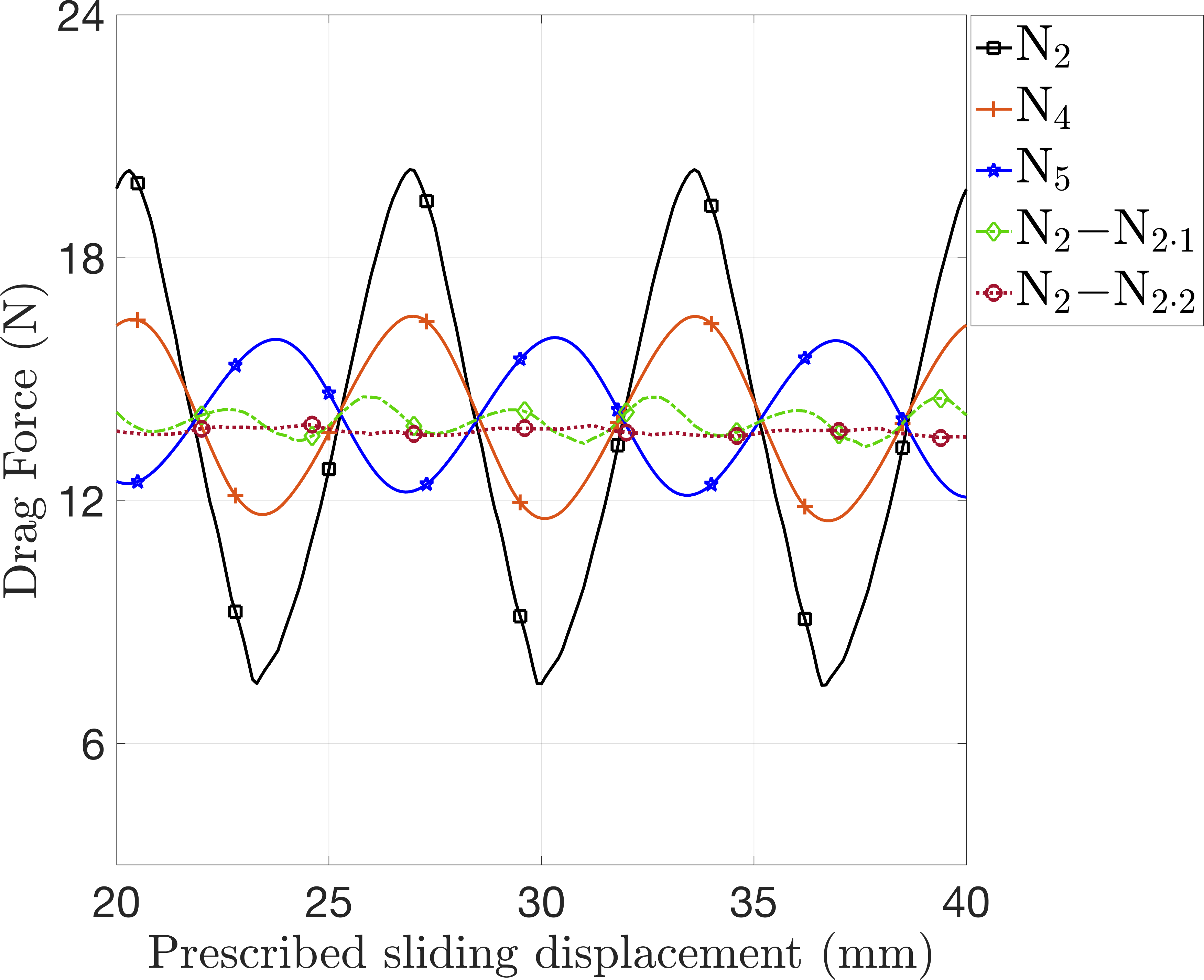}\label{fig:Ironing_Frictional_m1}}
	\caption{Frictional ironing: Enlarged view of (a) normal and (b) frictional contact forces with the VO and standard NURBS discretizations at mesh m$_1$.} \label{fig:ironing_result_m1}
\end{figure}

To illustrate the performance of our VO discretization over the standard NURBS discretizations, the enlarged views of the normal and horizontal contact forces are provided in Figs.~\ref{fig:Ironing_Normal_m1} and ~\ref{fig:Ironing_Frictional_m1}. For quantitative analysis, the reduction in the oscillation error with VO discretizations w.r.t. to popular standard N$_2$ is provided in Table~\ref{table:amplitude_error3D} for both the reaction forces. Here, we use: $ \Delta P_j := \text{max}(P_j) - \text{min}(P_j),~\text{where }j = z,\, \mathrm{and}\, x $ to compute the oscillation amplitude at the fixed mesh. Moreover, the total DOFs associated with the contact surface and remaining bulk domain of each solid with different discretizations are listed in Table~\ref{table:DOF_density_Irnoning3D}.
\begin{table}[!h]
\begin{center}
\begin{tabular}{| ccc | ccc | ccc |}
    \hline
    &\textbf{Discretization} &&& $ \Delta \textrm{P}_z\, (\%)$ &&& $ \Delta \textrm{P}_x \,(\%)$ &\\ 
    & \textbf{Type} &&& &&& & \\ \hline 
    &N$ _2 $  &&& 100   &&& 100   &\\ 
    &N$ _4 $  &&& 74.64 &&& 68.75 &\\ 
    &N$ _5 $  &&& 84.54 &&& 83.95 &\\ 
    &N$ _2- $N$ _{2\cdot1} $ &&& 79.97  &&& 93.16 &\\ 
    &N$ _2- $N$ _{2\cdot 2} $&&& 95.40  &&& 98.22 &\\ 
    \hline
\end{tabular}
\caption{Frictional ironing: Reduction in the oscillation amplitude of the normal and frictional contact forces with different standard and VO NURBS discretizations at m$_1$. The result with N$ _2 $ is used as a reference.}\label{table:amplitude_error3D}
\end{center}
\end{table}

\begin{table}[!h]
\begin{center}
\begin{tabular}{|c| c |c| c| c| c| c | c |}
\hline
\textbf{Discretization} & 
\multicolumn{3}{c|}{\textbf{DOFs for indentor}} & \multicolumn{3}{c|}{\textbf{DOFs for slab}} & \textbf{Total} \\[1ex]
\cline{2-7} 		
\textbf{Type} & {Interface} & {Bulk} & {Total} & {Interface} & {Bulk} & {Total} & \textbf{DOF}\\  
\hline
N$_2$ & 147 & 441 & 588 & 252 & 756 & 1008 & 1596 \\
N$_4$ & 243 & 729 & 972 & 384 & 1152 & 1536 & 2508 \\
N$_5$ & 300 & 900 & 1200 & 459 & 1377 & 1836 & 3036 \\
N$_2$-N$_{2\cdot 1}$ & 432 & 441 & 873 & 780 & 756 & 1536 & 2409 \\
N$_2$-N$_{2\cdot 2}$ & 867 & 441 & 1308 & 1596 & 756 & 2352 & 3660 \\
\hline
\end{tabular} \caption{Frictional ironing: DOFs details for the indentor and slab for different VO and standard NURBS-based discretizations at mesh m$_1$.} \label{table:DOF_density_Irnoning3D}
\end{center}
\end{table}

Based on the obtained results, the following three key observations are drawn. First, the VO NURBS-based discretization, specifically N$_2- $N$ _{p_c} $ ($ p_c = 2\cdot 1,$ and $ 2 \cdot 2 $) achieves a substantial improvement in accuracy by reducing the oscillation amplitude of forces compared to the N$_2$ based NURBS discretization. Specifically, N$_2- $N$ _{2\cdot 1} $ achieves approx. $80\%$ and $93\%$ reduction in the oscillation amplitude of the normal and frictional contact forces relative to N$_2$. Second, while the employment of higher-continuous NURBS discretization N$_5$ improves the accuracy over N$_2$ at the fixed mesh, N$_2- $N$ _{2\cdot 1} $ still delivers the results with superior accuracy. This improvement with the VO NURBS is attributed to the large number of the DOFs that are additionally present on the contact surface and moderate multi-knot span support of the higher-order basis function, as compared to the standard discretizations. The third observation is that further increasing the interpolation order of the NURBS for the contact layer, as in  N$_2- $N$ _{2\cdot 2} $ results in only a marginal improvement in accuracy over N$_2- $N$ _{2\cdot 1} $. In other words, this suggests that adding more number of the DOFs for the contact surface than N$_2- $N$ _{2\cdot 1} $ offers only a slight improvement in the accuracy. 

\subsubsection{Convergence analysis}
\begin{figure}[!ht]
	\centering
	\subfloat[]{\includegraphics[scale=0.25]{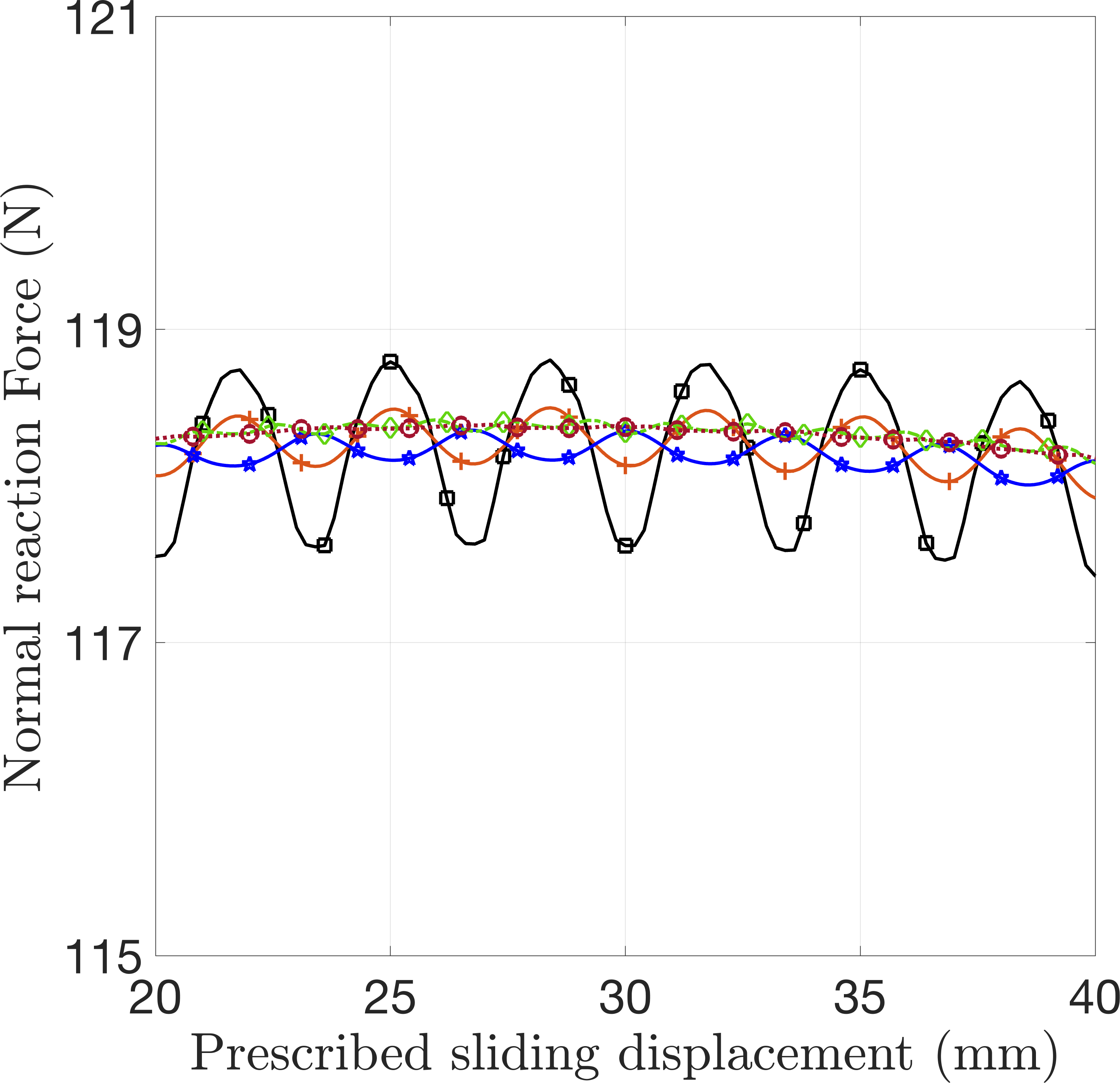}\label{fig:Ironing_Normal_m2}} ~~~
    \subfloat[]{\includegraphics[scale=0.25]{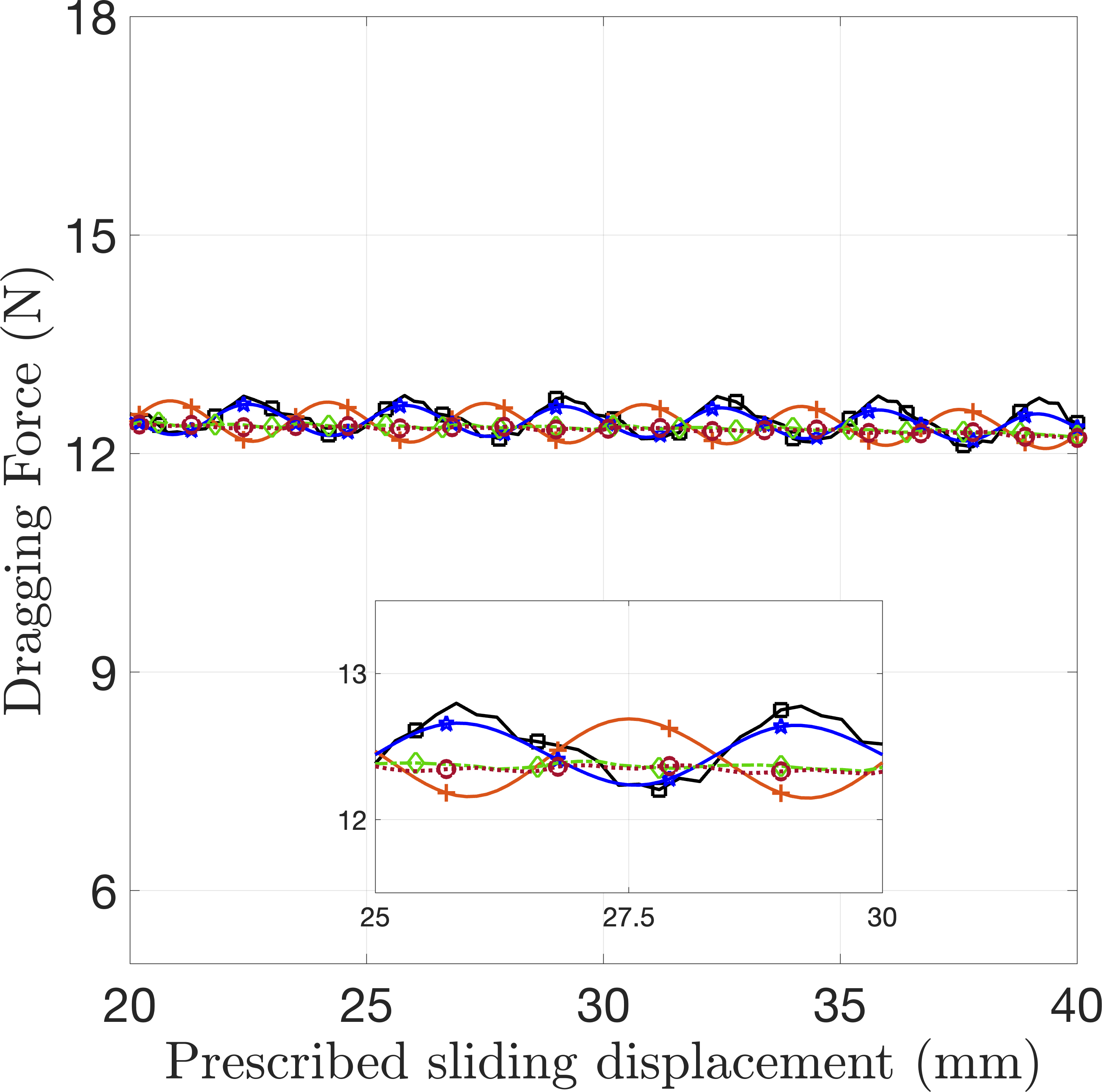}\label{fig:Ironing_Frictional_m2}}
	\caption{Frictional ironing: Enlarged view of (1) normal and (b) frictional contact forces with the VO and standard NURBS discretizations at mesh m$_2$..} \label{fig:ironing_result_m2}
\end{figure} 
Next, we analyze the mesh convergence of the normal and horizontal contact forces for both the VO and standard NURBS discretizations. For this, the resolution of the coarsest mesh for the indentor and slab is uniformly increased, as given in Table~\ref{table:ironing_element} with the uniform knot insertion along each parametric direction. The reduction in the oscillation amplitude of the normal and frictional contact forces with different discretizations at mesh level m$_2$ is illustrated in Figs.~\ref{fig:Ironing_Normal_m2} and Figs.~\ref{fig:Ironing_Frictional_m2}. As shown, the obtained results demonstrate the consistent performance of the VO NURBS over the standard uniform NURBS discretizations. Such a consistent gain in accuracy aligns with the observations made at mesh level m$_1$ with the VO and standard discretizations.

\begin{figure}[H]
\centering
\subfloat[]{\includegraphics[scale=0.25]{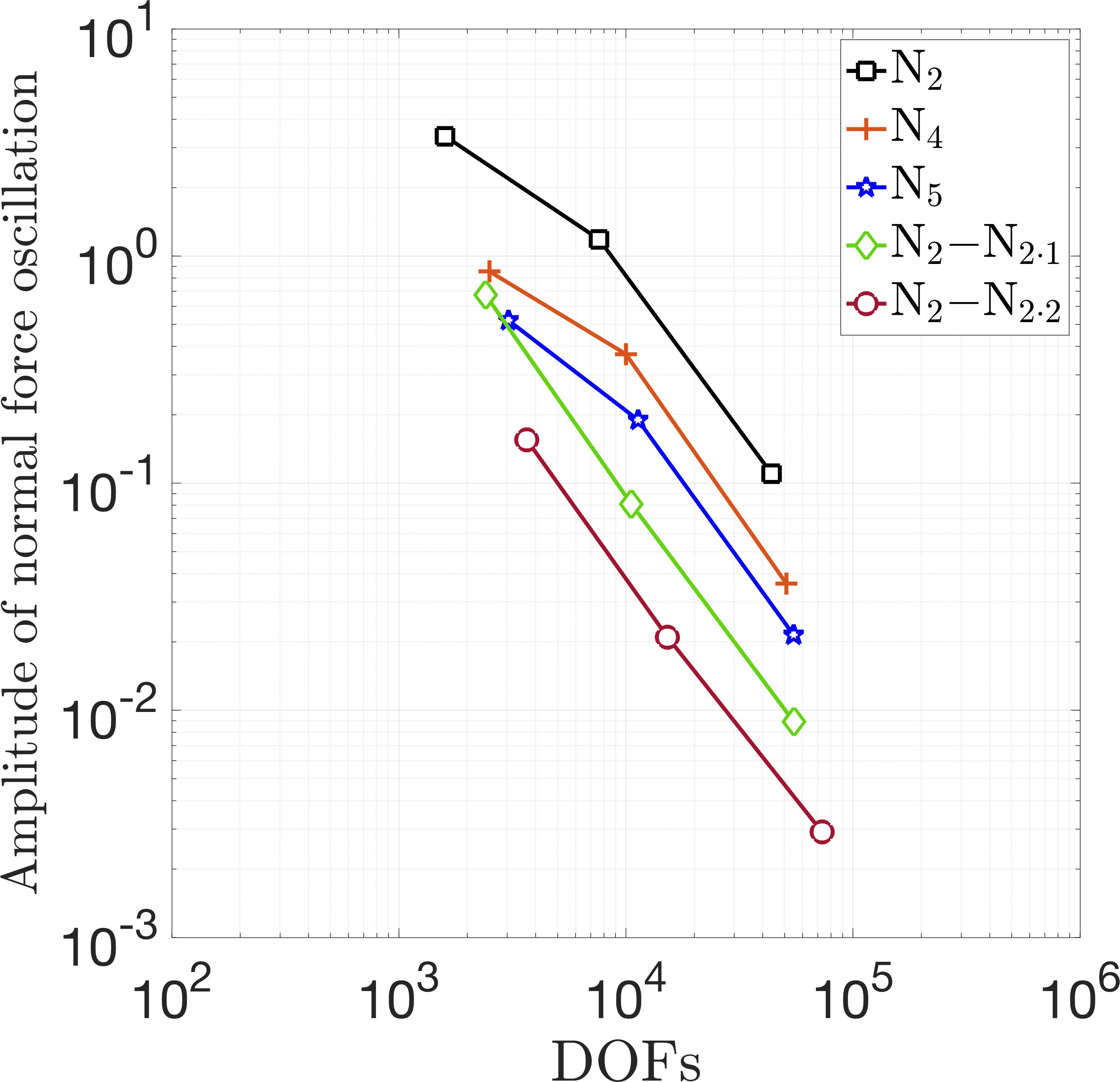}\label{fig:Ironing_normal_conv}} ~~
\subfloat[]{\includegraphics[scale=0.25]{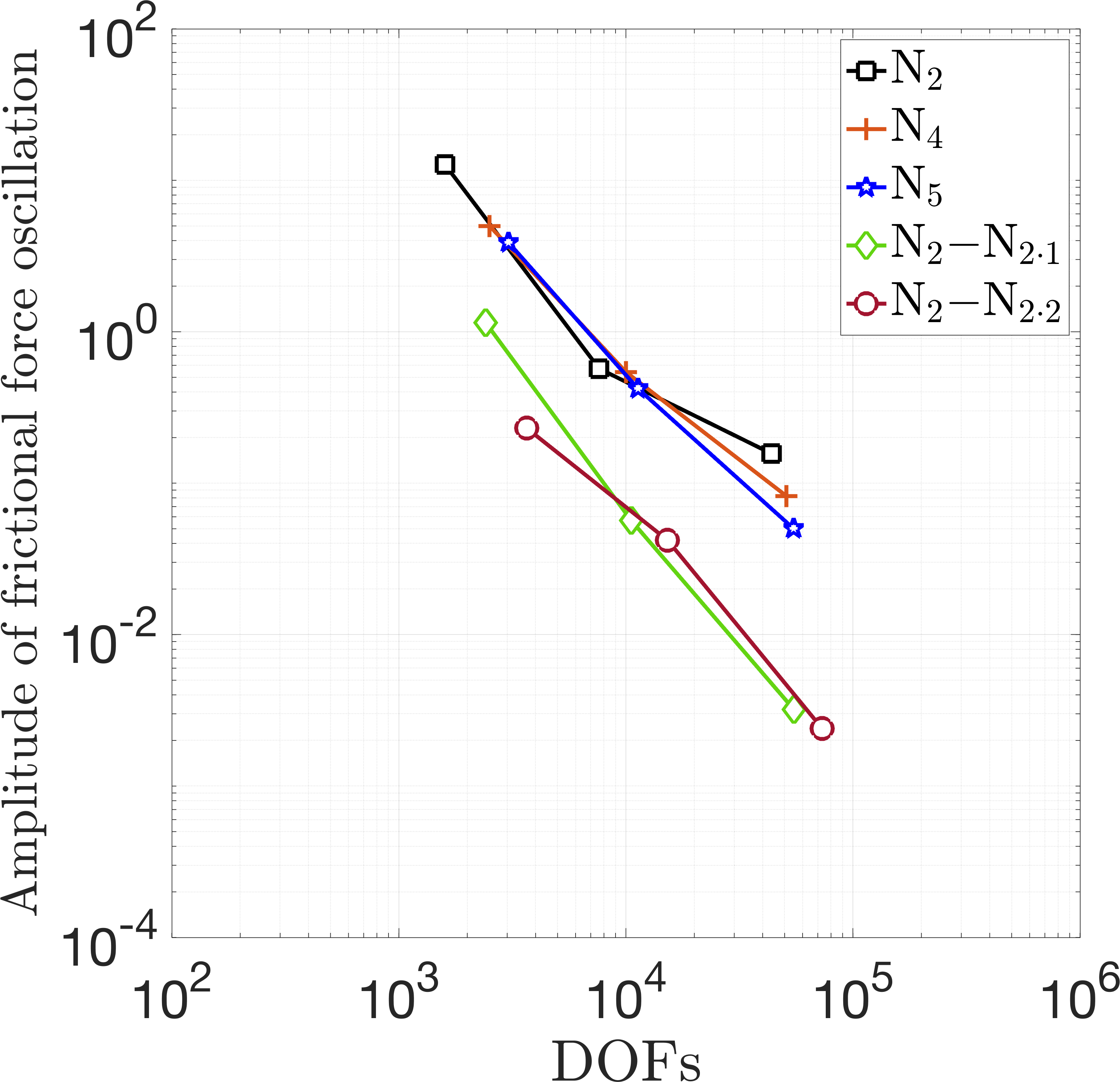}\label{fig:Ironing_Frictional_conv}}
\caption{Frictional ironing: Convergence of the oscillation amplitude of the (a) normal and (b) frictional contact forces for different VO and standard NURBS discretizations upon mesh refinement.} \label{fig:ironing_conv}
\end{figure}

The convergence plots for both the normal and frictional contact force oscillation amplitude over total DOFs with mesh refinement are presented in Figs.~\ref{fig:Ironing_normal_conv} and~\ref{fig:Ironing_Frictional_conv}. It can be seen that VO discretization offers significant accuracy improvement compared to the standard NURBS discretization at the same mesh level. A closer look reveals that the accuracy obtained with N$_2- $N$ _{2\cdot 1} $ at mesh levels m$_1$ and m$_2$ is comparable to results with N$_2$ at subsequent finer meshes m$_2$ and m$_3$. For this, N$_2- $N$ _{2\cdot 1} $ takes approx. $3.2$ and $4.1$ times fewer DOFs, leading to a considerable gain in computational efficiency. Furthermore, the higher-order VO discretization, N$_2- $N$ _{2\cdot 2} $ yields a slight improvement in the accuracy compared to N$_2- $N$ _{2\cdot 1} $ at the same mesh level.

\subsection{Twisting contact between a hemisphere and a cube}
In this section, we show the performance of the VO discretization method using a challenging 3D benchmark: twisting contact between a deformable hemisphere and a cube, as introduced in~\cite{Corbett2014, Corbett2015}. We first analyze the frictionless twisting scenario to showcase the capability of the method to handle pure normal contact in 3D. In the subsequent example, friction is included, resulting in a more complex interaction involving frictional sticking and sliding, in addition to normal contact during twisting. This progression allows us to assess the robustness of the method in handling complex contact scenarios.

\subsubsection{Frictionless Twisting}
\begin{figure}[!ht]
	\centering
	\subfloat[]{\includegraphics[scale=0.325]{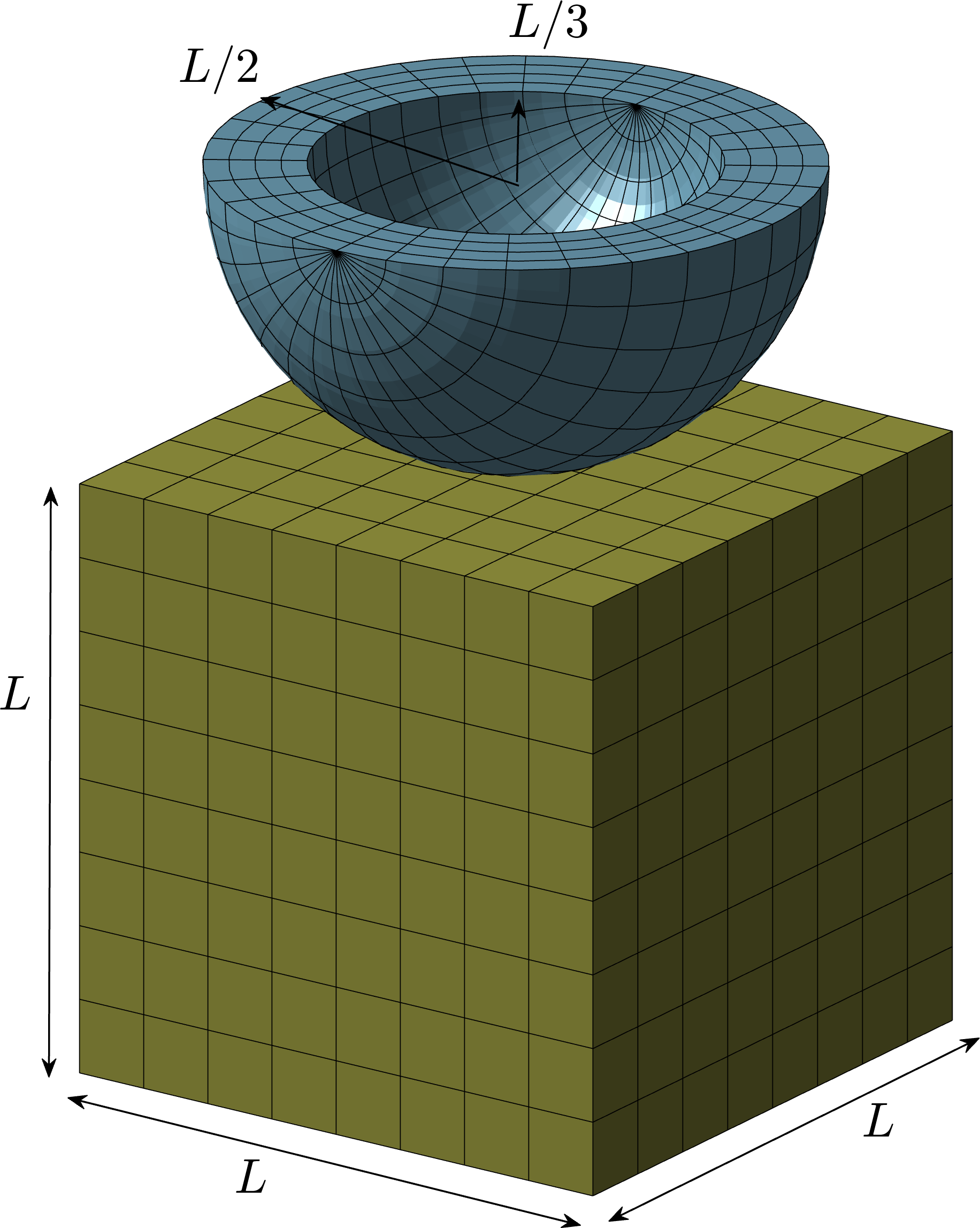}\label{fig:twist_setup}} ~~~
	\subfloat[]{\includegraphics[scale=0.325]{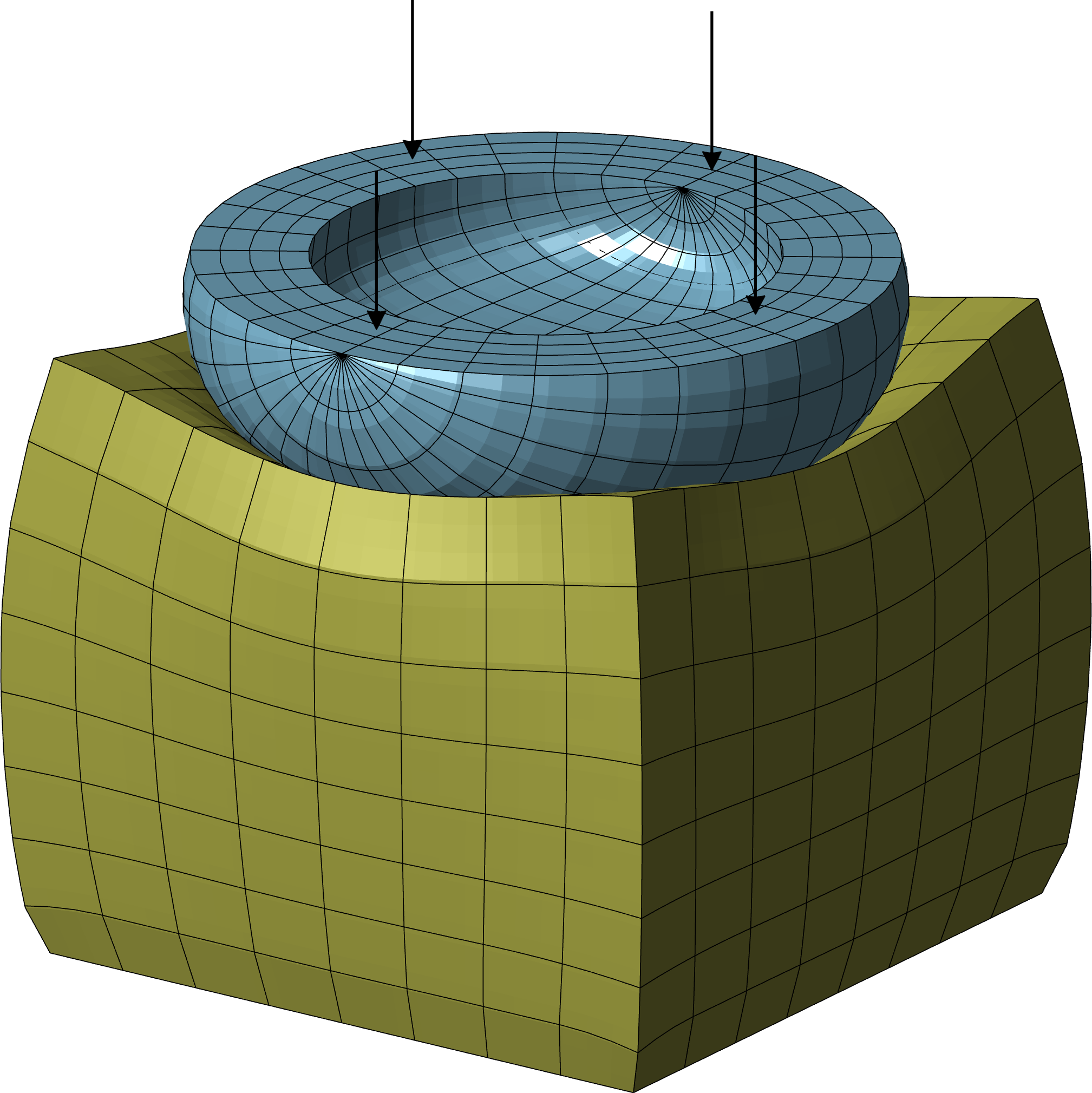}\label{fig:twist_compression}} 
	\caption{Setup for both the frictionless and frictional twisting contact. (a) Geometry and finest mesh for the hemisphere and cube. (b) Deformed configuration of the setup at the end of compression.} \label{fig:twisting_setup_compressed}
\end{figure}
The setup of the problem, which involves a hollow hemisphere having an outer radius $L/2$ mm and an inner radius $L/3$ mm, and a cube of size $[L \times L \times L ],$ where $L = 2$ mm is shown in Fig.~\ref{fig:twist_setup}. In this benchmark problem, the hemisphere is first pressed into a cube by applying the prescribed vertical downward displacement of $U_z = L/2 $ mm to all the control points on the ring-shaped top surface of the hemisphere. The vertical displacement is then held constant in the subsequent frictionless rotation of $180^{o}$ applied around the vertical axis in uniform load steps. The foundation of the cube is held stationary by fixing displacement in all directions. A Neo-Hookean hyper-elastic material model is considered with material properties: $E= 5$ N/mm$^2$ for the hemisphere and $E= 1$ N/mm$^2$ for the cube, and Poisson's ration $\nu=0.3$ for both. Three nested meshes that are obtained with the uniform knot insertion along each parametric direction are summarised in Table~\ref{table:twisting_element}. The penalty parameter used for enforcing the normal contact constraint is also included in Table~\ref{table:twisting_element} with each mesh level. Fig.~\ref{fig:twist_setup} depicted the finest mesh used for the setup. Further, as can be seen, the exact NURBS-based construction of the hemisphere consists of degenerated triangle elements at its two polar points. However, it noted that those degenerated elements do not negatively affect the accuracy of the result as the displacement of the degenerated points is obtained with the prescribed Dirichlet boundary condition. Moreover, the NURBS-based construction of the hemisphere enables the achievement of the required C$^1$-continuity across the contact surface with N$_2$ based discretization. The DOFs over the contact surface and bulk domain of the hemisphere and cube with each NURBS discretization at different meshes are provided in Table~\ref{table:twisting_DOFs}. Furthermore, the deformed configuration of the setup at the end of the compression of the hemisphere against the cube is shown in Fig.~\ref{fig:twist_compression}.

\begin{table}[!ht]
\begin{center}
\begin{tabular}{|c c c| c c c | c c c| c c c|}
\hline
&\textbf{Mesh} && \multicolumn{6}{c|}{\textbf{Elements}} & \multicolumn{3}{c|}{$\epsilon_N$} \\ [0.4ex]
\cline{4-9}		
&& && {Hemisphere} &&& {Cube} &&& &\\ 
\hline
&m$ _1 $ &&&  $ 4 \times 4 \times 2 $ &&&  $ 2 \times 2 \times 2 $  &&& 100 &\\
&m$ _2 $ &&&  $ 8 \times 8 \times 3 $ &&&  $ 4 \times 4 \times 4 $  &&& 200 &\\
&m$ _3 $ &&&  $ 16 \times 16 \times 4 $ &&&  $ 8 \times  8 \times 8 $ &&& 400 &\\
\hline
\end{tabular} \caption{Elements with three meshes used for hemisphere and cube along with the penalty parameter at each mesh for both frictionless and frictional twisting.} \label{table:twisting_element}
\end{center}
\end{table}
\begin{table}[!hb]
\begin{center}
\begin{tabular}{|c| c |c| c| c| c| c | c | c|}
\hline
\textbf{Discretization} & \textbf{Mesh} & \multicolumn{3}{c|}{\textbf{DOFs for hemisphere}} & \multicolumn{3}{c|}{\textbf{DOFs for cube}} & \textbf{Total} \\[1ex]
\cline{3-8}
\textbf{Type} & & {Interface} & {Bulk} & {Total} & {Interface} & {Bulk} & {Total} & \textbf{DOF}\\  
\hline
N$_2$ & m$_1$ & 147 & 294 & 441 & 48 & 96 & 144 & 585 \\
N$_2$-N$_{2\cdot 1}$ & m$_1$ & 363 & 294 & 657 & 108 & 96 & 204 & 861 \\
N$_2$-N$_{2\cdot 2}$ & m$_1$ & 675 & 294 & 969 & 192 & 96 & 288 & 1257 \\
\hline
N$_2$ & m$_2$ & 363 & 1089 & 1452 & 108 & 432 & 540 & 1992 \\
N$_2$-N$_{2\cdot 1}$ & m$_2$ & 1083 & 1089 & 2172 & 300 & 432 & 732 & 2904 \\
N$_2$-N$_{2\cdot 2}$ & m$_2$ & 2187 & 1089 & 3276 & 588 & 432 & 1020 & 4296 \\
\hline
N$_2$ & m$_3$ & 1083 & 4332 & 5415 & 300 & 2400 & 2700 & 8115 \\
N$_2$-N$_{2\cdot 1}$ & m$_3$ & 3675 & 4332 & 8007 & 972 & 2400 & 3372 & 11379 \\
N$_2$-N$_{2\cdot 2}$ & m$_3$ & 7803 & 4332 & 12135 & 2028 & 2400 & 4428 & 16563 \\
\hline
\end{tabular}
\caption{DOFs details over the contact surface and in volume of the hemisphere and cube for different VO and standard NURBS-based discretizations at three nested meshes for both frictionless and frictional twisting contact.}
\label{table:twisting_DOFs}
\end{center}
\end{table}

\begin{figure}[!ht]
	\centering
	\subfloat[]{\includegraphics[scale=0.255]{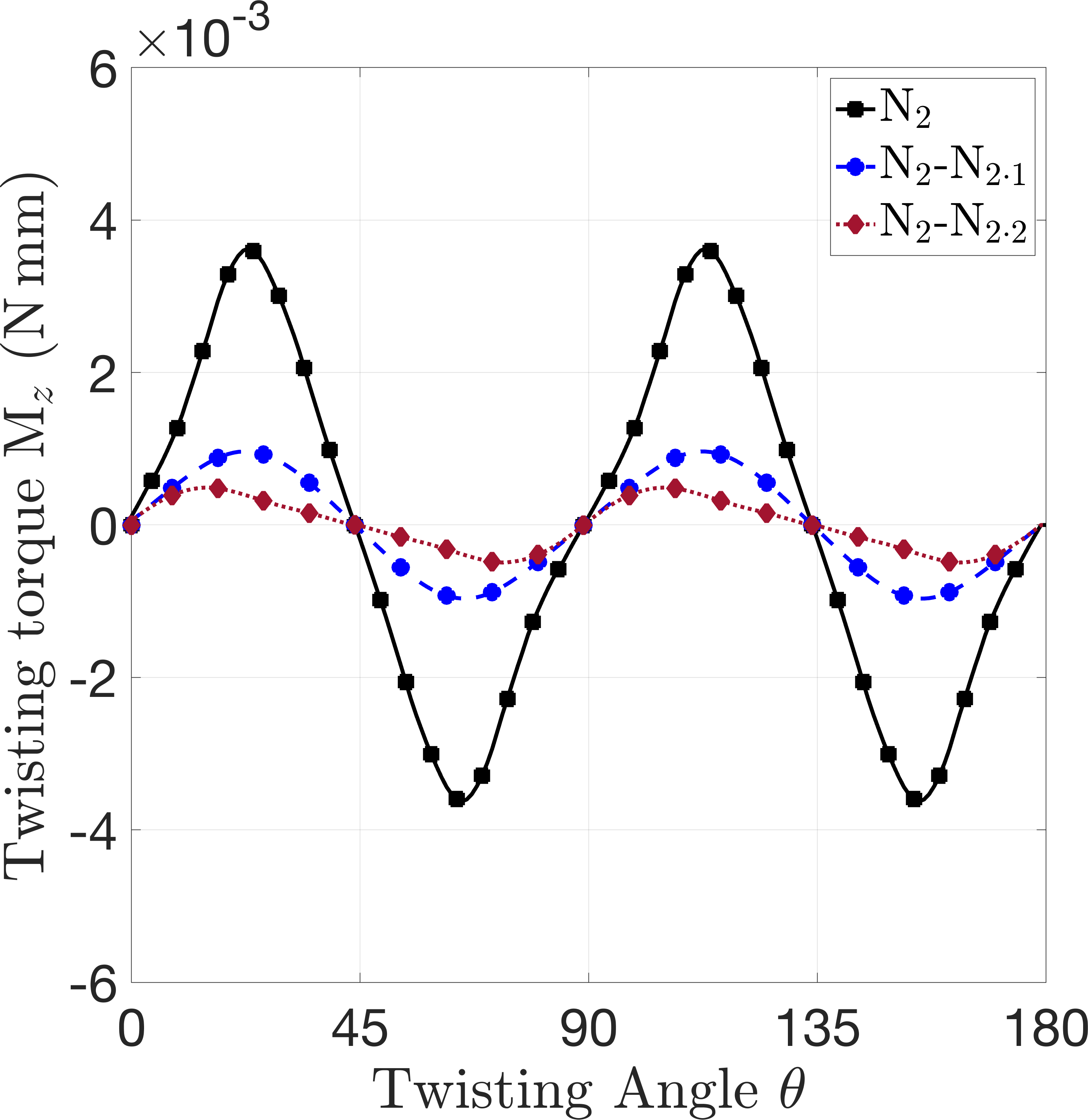}\label{fig:twist_fless_m1}} ~~
	\subfloat[]{\includegraphics[scale=0.255]{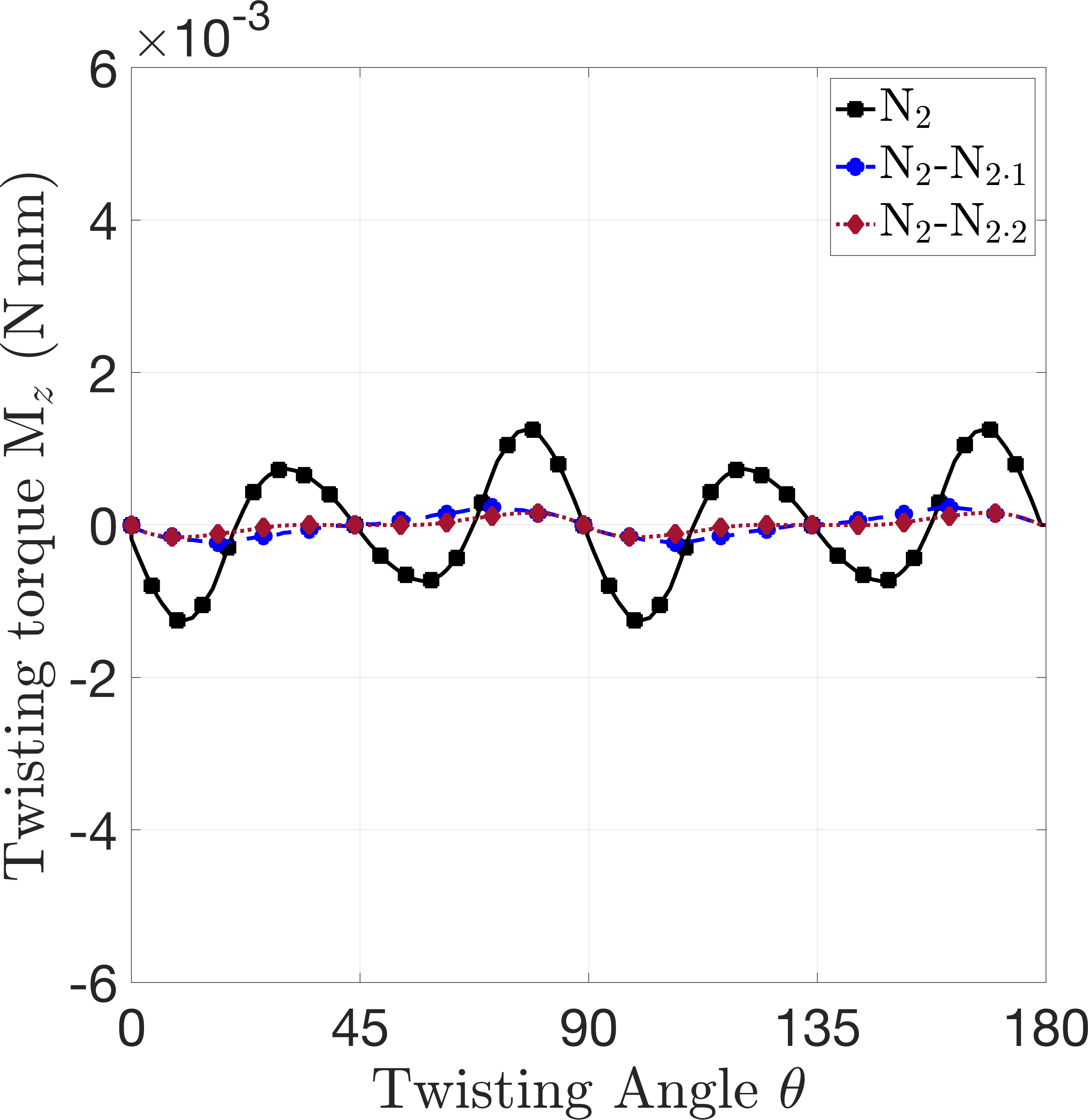}\label{fig:twist_fless_m2}} 
	\caption{Frictionaless twisting contact: Variation of torque around the vertical axis of the setup over prescribed twisting angle for different VO and standard NUBRS discretization at mesh (a) m$_1$ and (b) m$_2$.} \label{fig:twisting_fless_result}
\end{figure}
First, we show the performance of the proposed method at a very coarse mesh. The evolution of torque computed around the rotation axis over the prescribed rotation angle with VO discretizations is shown in Fig.~\ref{fig:twist_fless_m1} at the mesh level m$_1$. The result with the most popular standard N$_2$ based discretization is also included for the comparative assessment. Considering the frictionless twisting, theoretically, the torque should be zero for the prescribed rotation angles. However, as noted above, such non-physical oscillation of torque, as can be seen from Fig.~\ref{fig:twist_fless_m1}, stems from using insufficient resolution for the contact surface, even if it is constructed with the required C$^1$-continuous NURBS. One of the most common practices to alleviate such a numerical error in the contact responses is to use a finer mesh. At a fixed mesh, the employment of higher-order basis functions offered by IGA can also lead to the reduction of such errors. However, using the finer mesh or standard higher-order NURBS functions in the region away from the contact surface may not be required, specifically the higher-order NURBS, due to their higher computational cost.

From Fig.~\ref{fig:twist_fless_m1}, it can be observed that the VO-based N$_2- $N$ _{2\cdot 1} $ discretization is able to capture the torque more accurately as compared to N$_2$ at the very coarse mesh. It reduces the amplitude error in torque by approx. $73\%$ over N$_2$ at mesh m$_1$. Such a gain in accuracy with VO is attributed to the additional DOFs present at the contact surface. Another observation is that although the application of further higher-order NURBS, i.e. N$_2- $N$ _{2\cdot 2} $ at the contact surface further improves the accuracy of the torque, the gain is marginal as it reduces amplitude error by only approx. $11\%$ over N$_2- $N$ _{2\cdot 1} $. Hence, as also observed in the frictional sliding example, using more number of DOFs for contact surface than with N$_2- $N$ _{2\cdot 1} $ is not useful, as it yields only a marginal gain in accuracy at a fixed mesh.

Next, we illustrate the performance of the VO discretizations at a finer mesh, m$_2$ in Fig.~\ref{fig:twist_fless_m2}. It is evident that the VO discretization, specifically N$_2- $N$ _{2\cdot 1} $, showcases superior performance over N$_2$ at mesh m$_2$. It is consistent with the observation made at mesh m$_1$ for VO compared to N$_2$. Another observation is that, at this mesh, the accuracy of the result obtained with N$_2- $N$ _{2\cdot 2} $ becomes comparable to N$_2- $N$ _{2\cdot 1} $, hence supporting the argument that more DOFs than obtained with N$_2- $N$ _{2\cdot 1} $ are not useful, as N$_2- $N$ _{2\cdot 2} $ only marginally improves the accuracy of the contact response at a fixed mesh.

\begin{figure}[!ht]
	\centering
	{\includegraphics[scale=0.275]{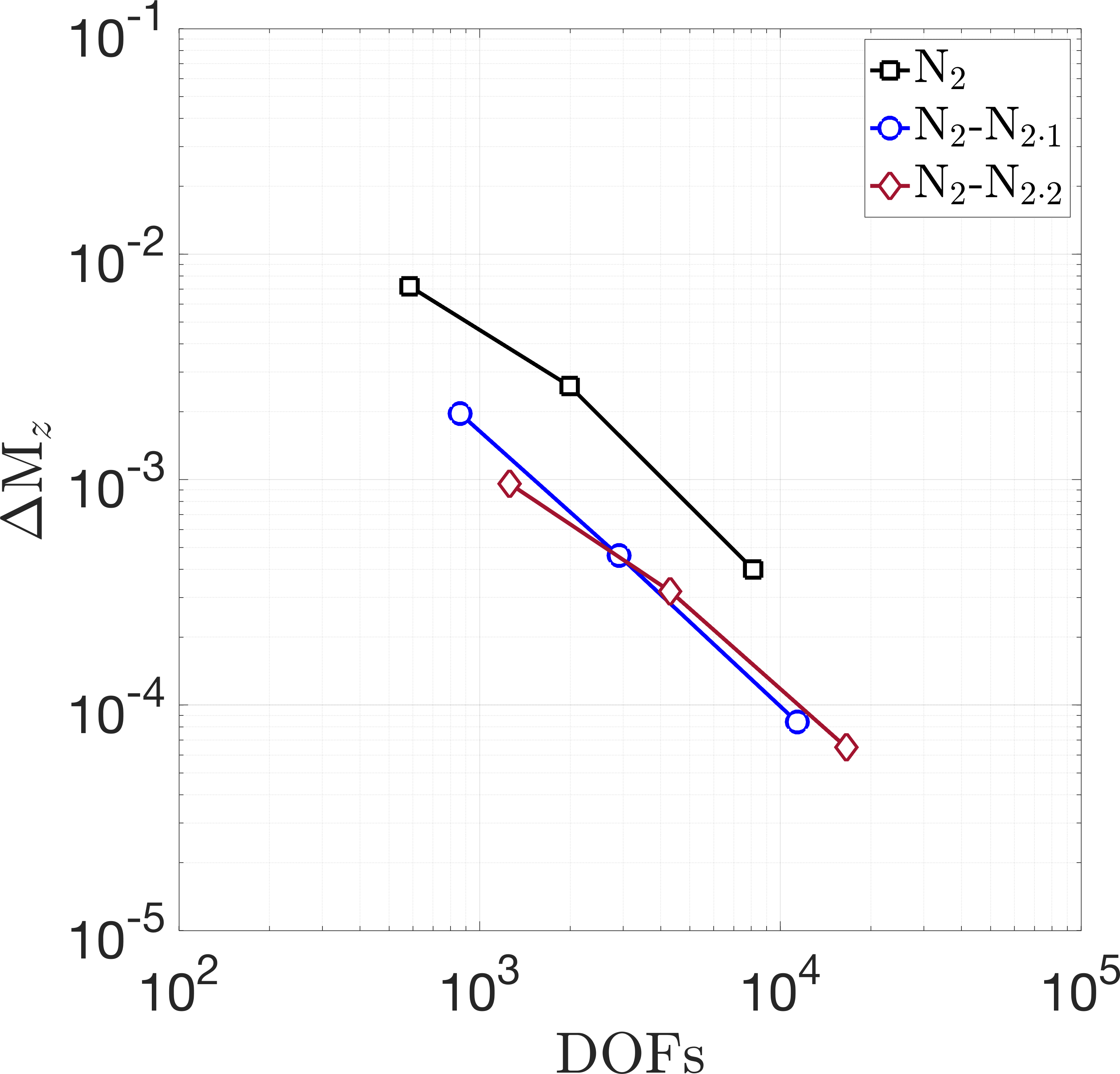}}
	\caption{Frictionaless twisting contact: Convergence of the amplitude error in the torque with different VO and standard NURBS discretization upon mesh refinement.} \label{fig:twisting_fless_conv}
\end{figure}
Finally, we investigate the convergence behavior of the reduction in the oscillation amplitude of the torque with different discretizations upon mesh refinement. Fig.~\ref{fig:twisting_fless_conv} shows the convergence plots for the VO and standard NURBS-based discretization. It can be observed that VO discretization clearly outperforms standard NURBS discretization at the same mesh level. Specifically, to attain the accuracy similar to N$_2$ at meshes m$_2$ and m$_3$, N$_2- $N$ _{2\cdot 1} $ uses approx. $2.32$ and $2.79$ times fewer number of DOFs. The VO discretization, thus, takes lower computational efforts in terms of the DOFs over the standard N$_2$ based NURBS discretization.

\subsubsection{Twisting with friction}
In this example, we use the same problem setup as in the previous one, except that frictional contact is additionally considered between the hemisphere and cube during the twisting. It is noted that during the downward vertical displacement of the hemisphere, frictionless contact is still assumed as in~\cite{Corbett2015}. Following~\cite{Corbett2015}, a high coefficient of friction $\mu = 0.5$ is applied at the interface, which leads to large shear deformation in both the contact surfaces and the interior of the bulk domains. As in the previous example, the normal penalty parameter $\epsilon_N$ is scaled in proportion to the average element length, and correspondingly, the tangential penalty parameter is set as $\epsilon_T = \epsilon_N$ at each mesh level. The same VO and standard NURBS discretization types, outlined in Table~\ref{table:twisting_DOFs} for frictionless twisting, are employed here. It is also noted that at the onset of the twisting of the hemisphere, when $\theta = 0^o$, all the contact points on the interface are considered to be in the sticking state, as per in~\cite{Corbett2015}. 
\begin{figure}[!ht]
\centering
\subfloat{\includegraphics[scale=0.18]{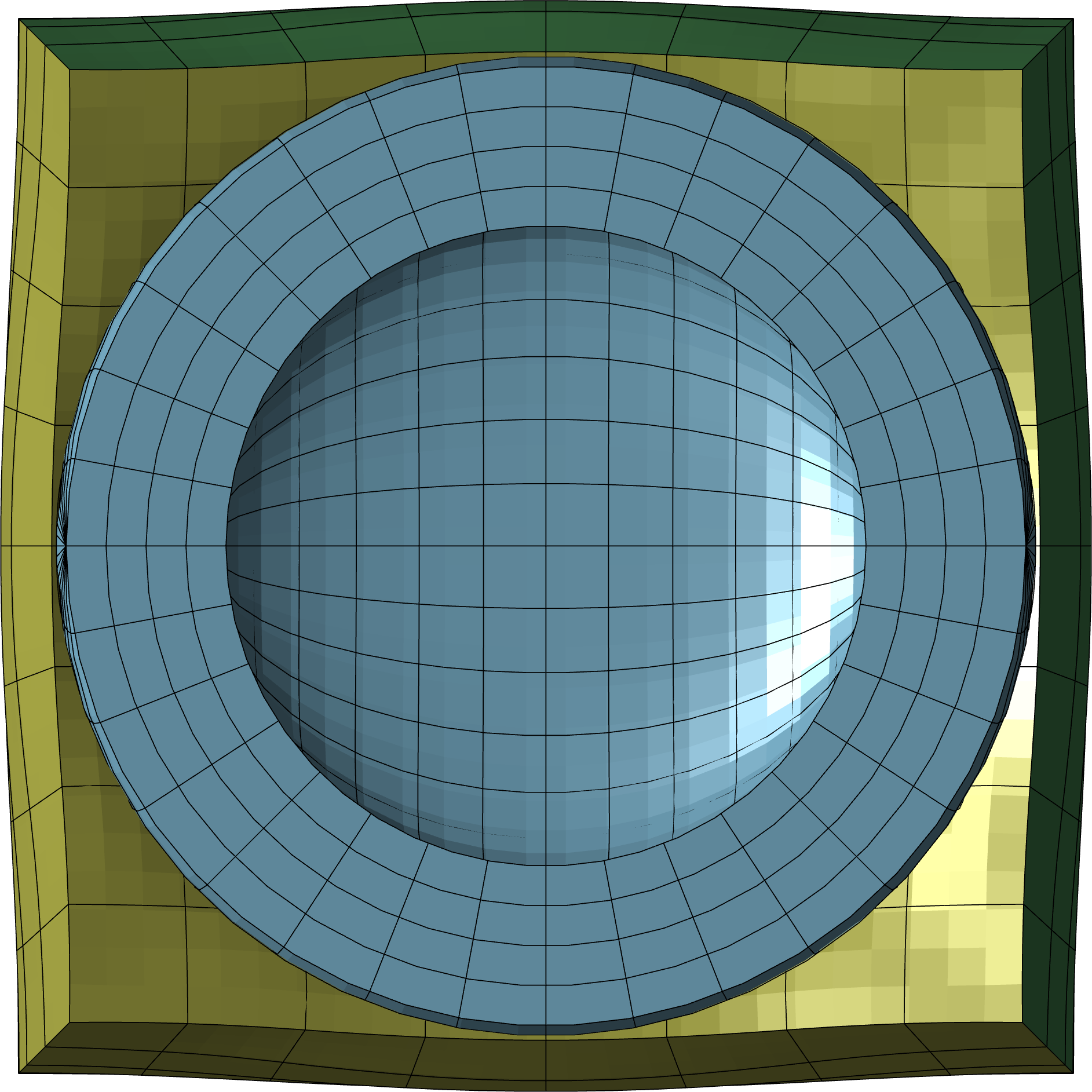}\label{fig:twist_iso0}}
\subfloat{\includegraphics[scale=0.18]{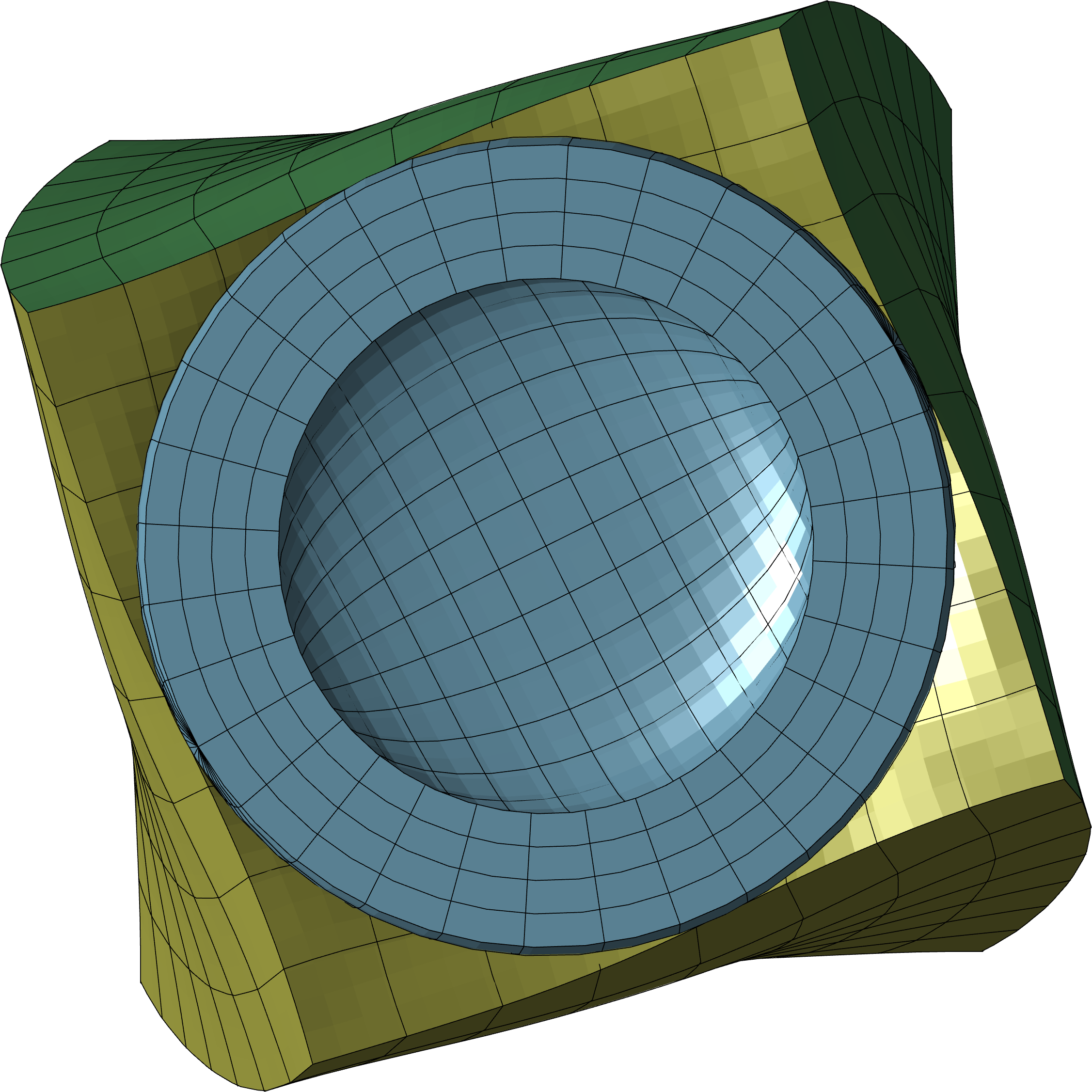}\label{fig:twist_iso30}}
\subfloat{\includegraphics[scale=0.18]{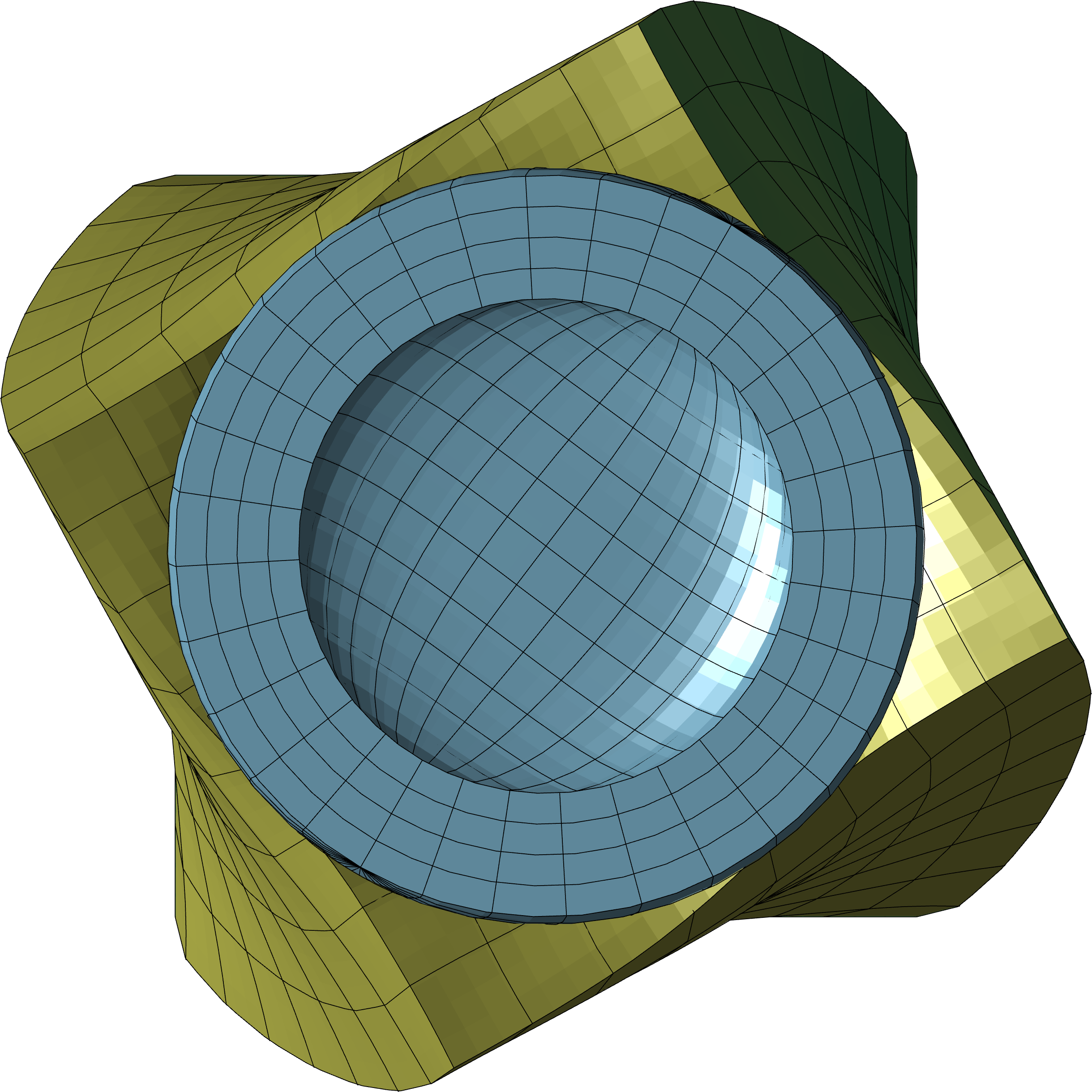}\label{fig:twist_iso60}}
\subfloat{\includegraphics[scale=0.18]{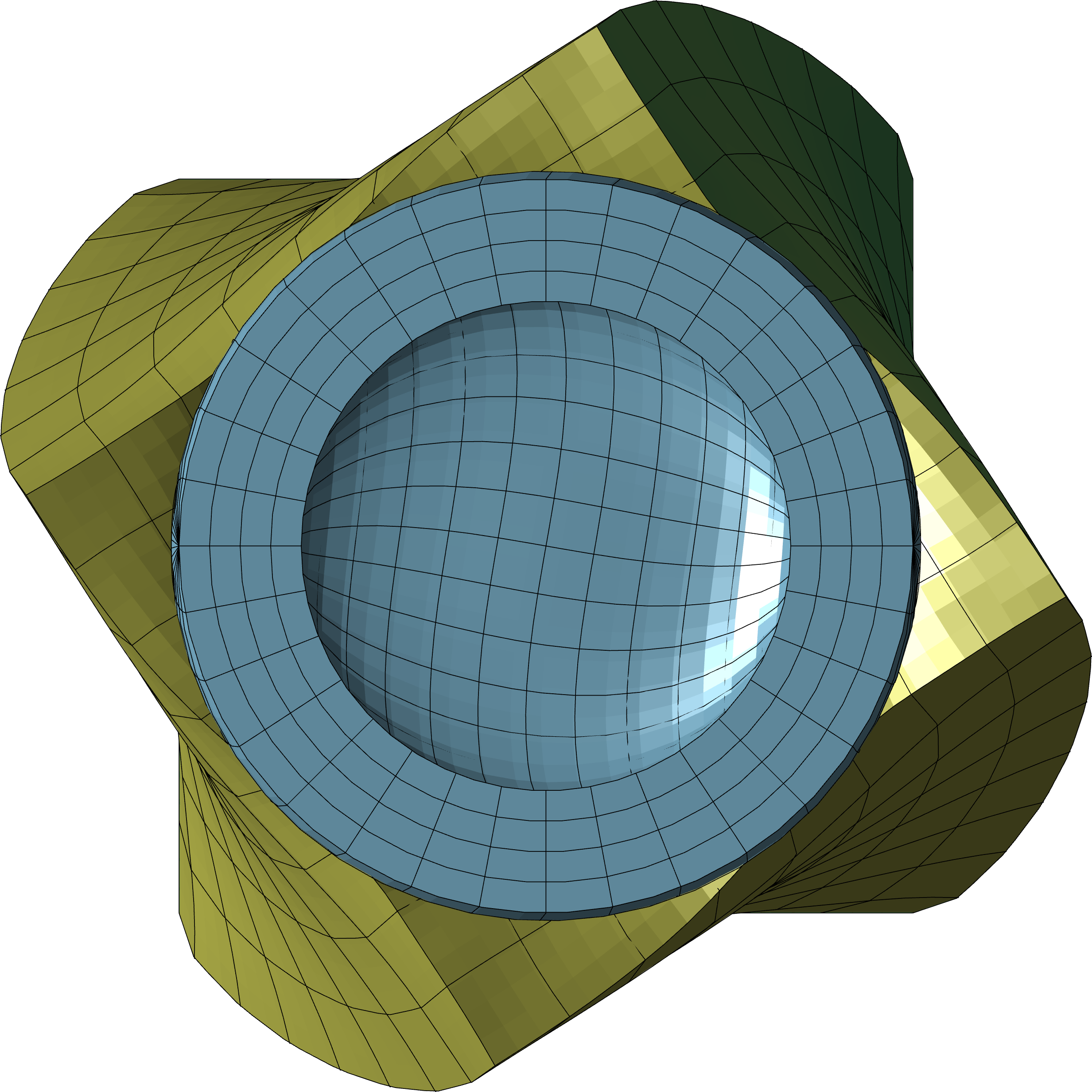}\label{fig:twist_iso180}} \\
\subfloat{\includegraphics[scale=0.2]{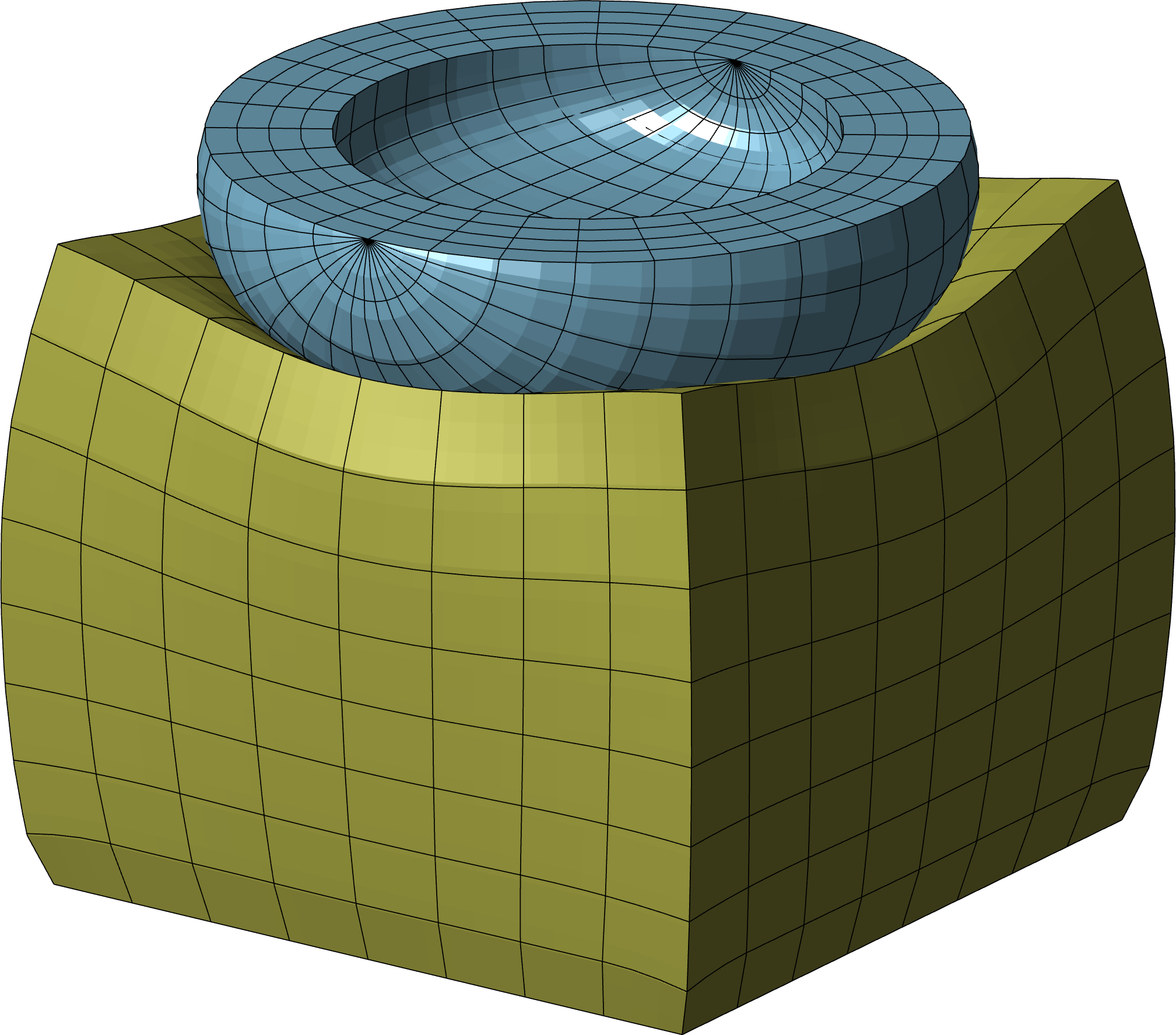}\label{fig:twist0}}
\subfloat{\includegraphics[scale=0.2]{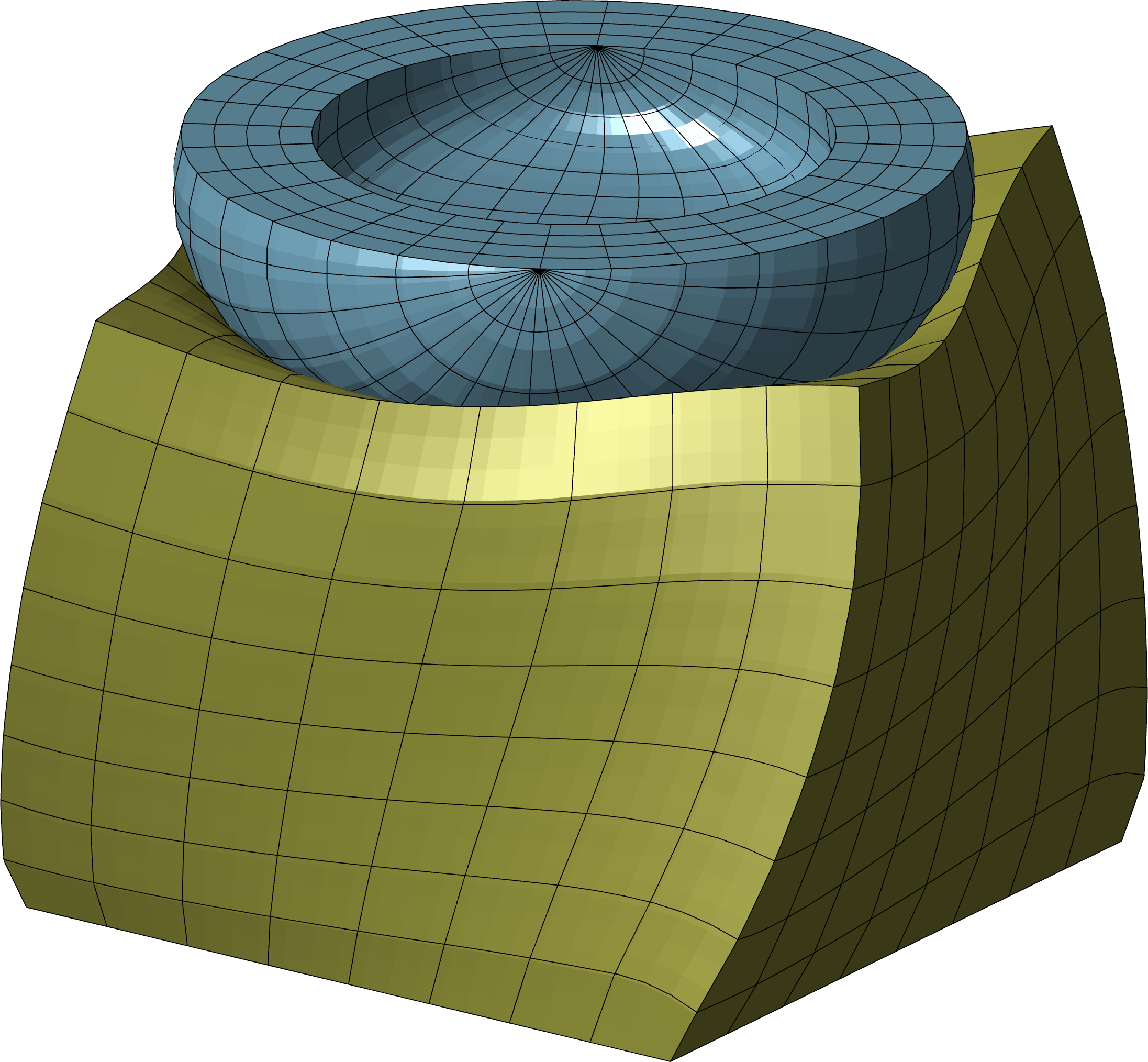}\label{fig:twist30}}
\subfloat{\includegraphics[scale=0.2]{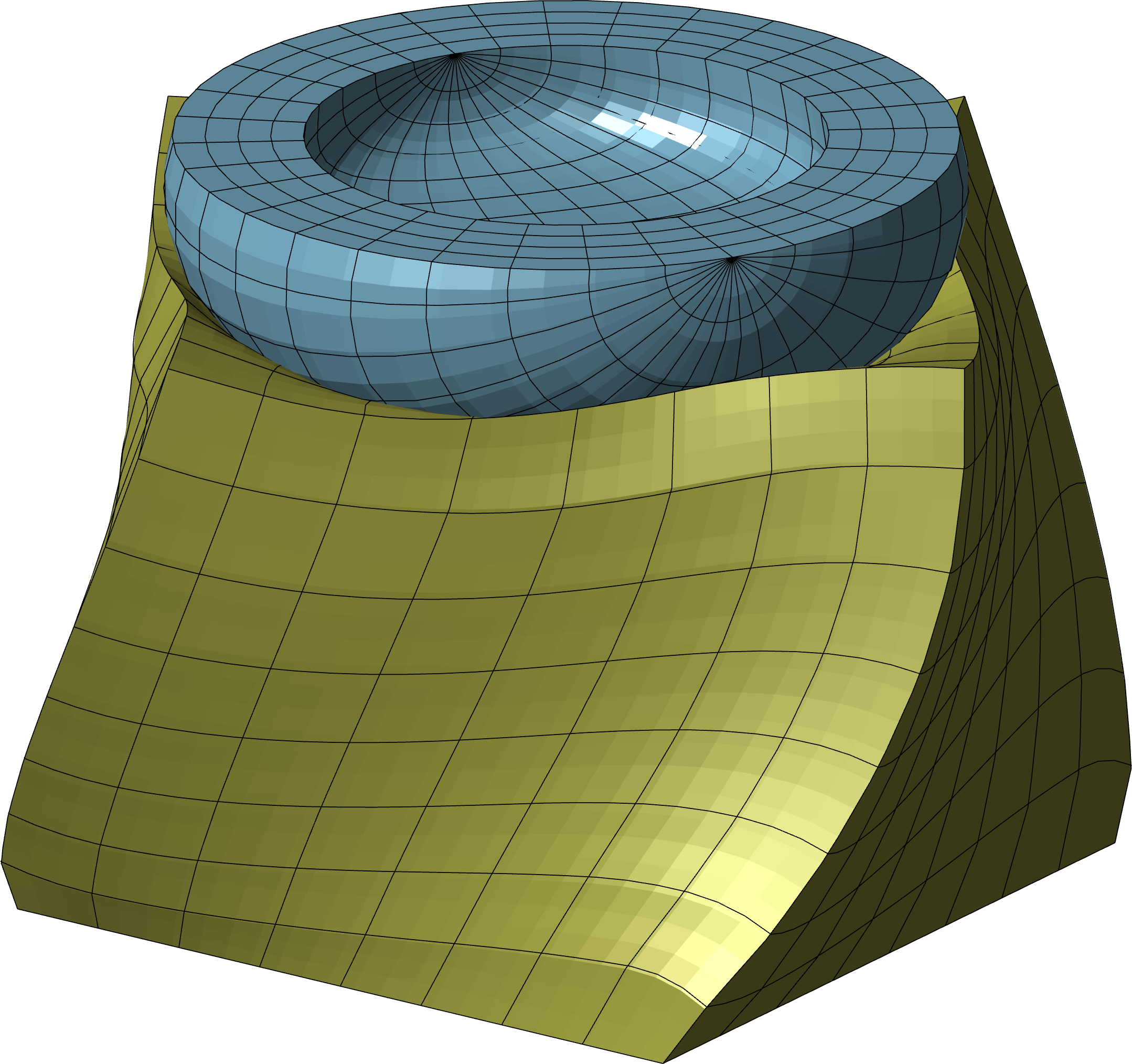}\label{fig:twist60}}
\subfloat{\includegraphics[scale=0.2]{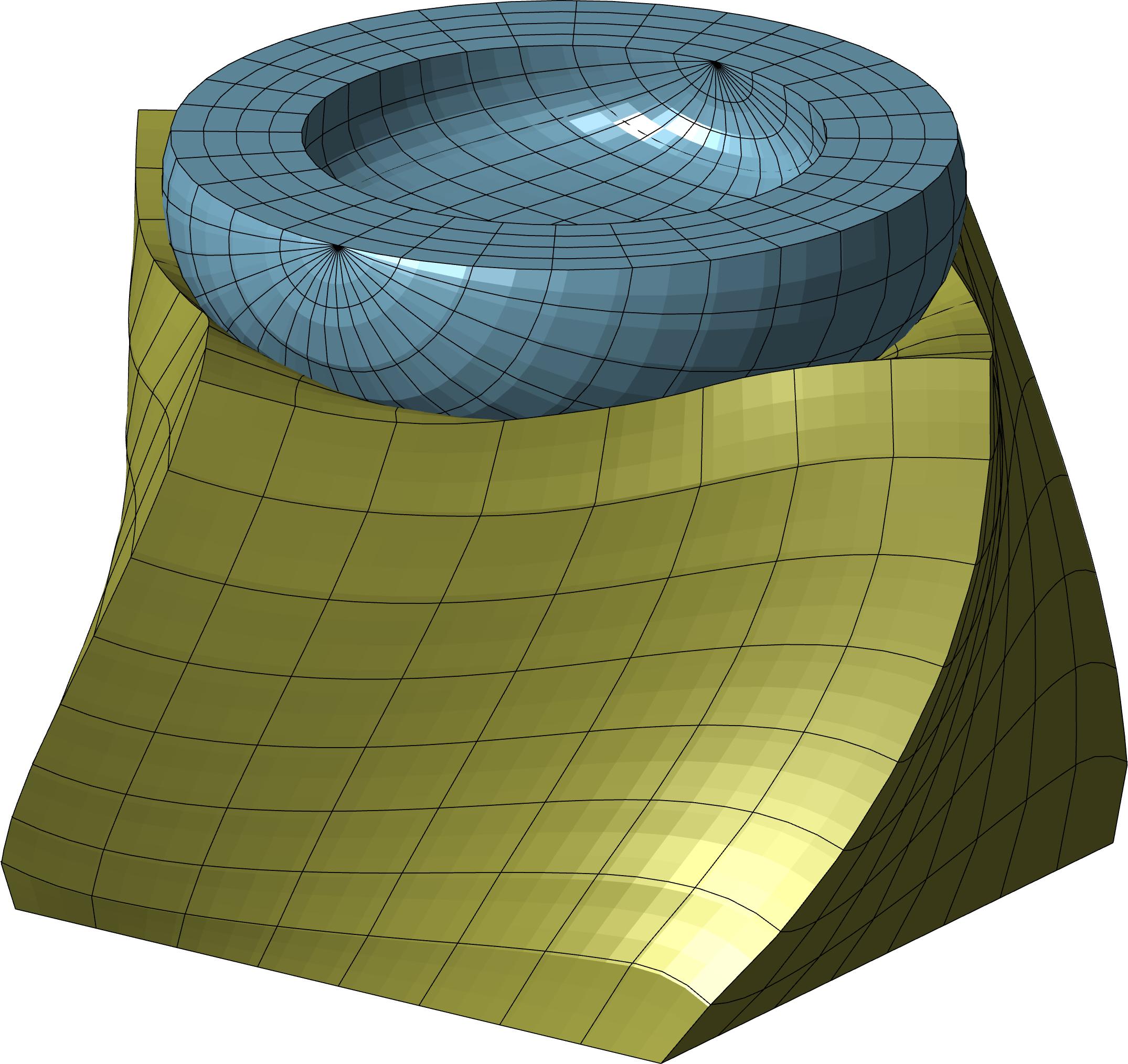}\label{fig:twist180}}
\caption{Frictional twisting contact: Deformed configurations of the setup (bottom row) and its top view (top row) at rotation angle $\theta = 0^o, 30^o, 60^o$ and $180^o$ (from left to right).} \label{fig:twisting_deformed_configs}
\end{figure}

Figure~\ref{fig:twisting_deformed_configs} illustrates the deformed configurations of the setup at various rotation angles. As evident, the high friction coefficient induces a significant shift in the frictional contact behavior. At the onset of twisting, when $\theta = 0^o$, the contact regions predominantly experience sticking due to the high resistance to the tangential motion. However, as the rotation progresses, the transition to slipping of the contact points occurs, and approx. after $\theta = 100^o$, slipping becomes the dominant behavior in the contact regions. This shift aligns with the observations in~\cite{Corbett2015}, where the same frictional effects were noted under such a high friction coefficient.

\begin{figure}[!htb]
	\centering
	\subfloat[]{\includegraphics[scale=0.25]{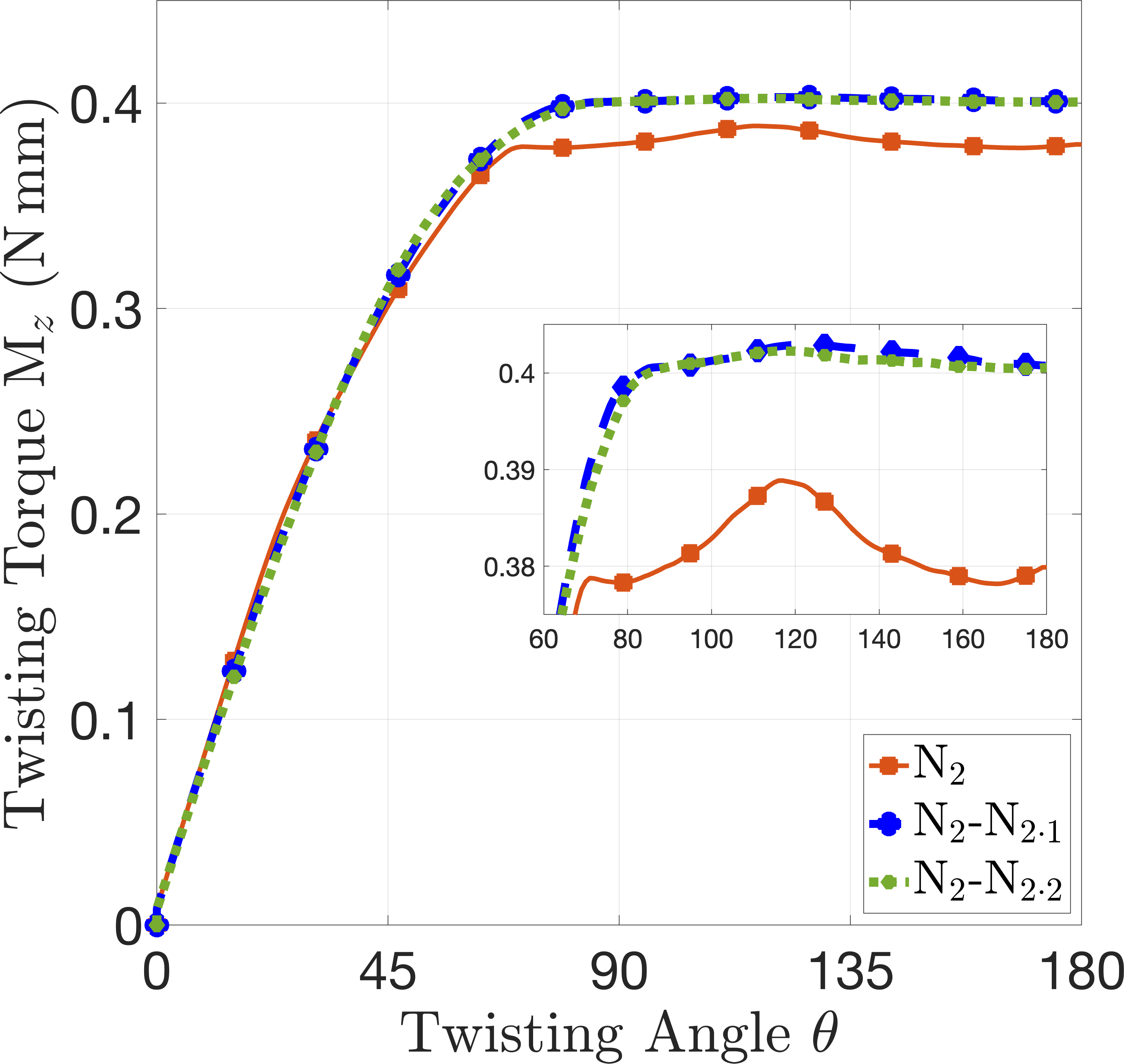}\label{fig:twist_m1}} ~~~~
	\subfloat[]{\includegraphics[scale=0.25]{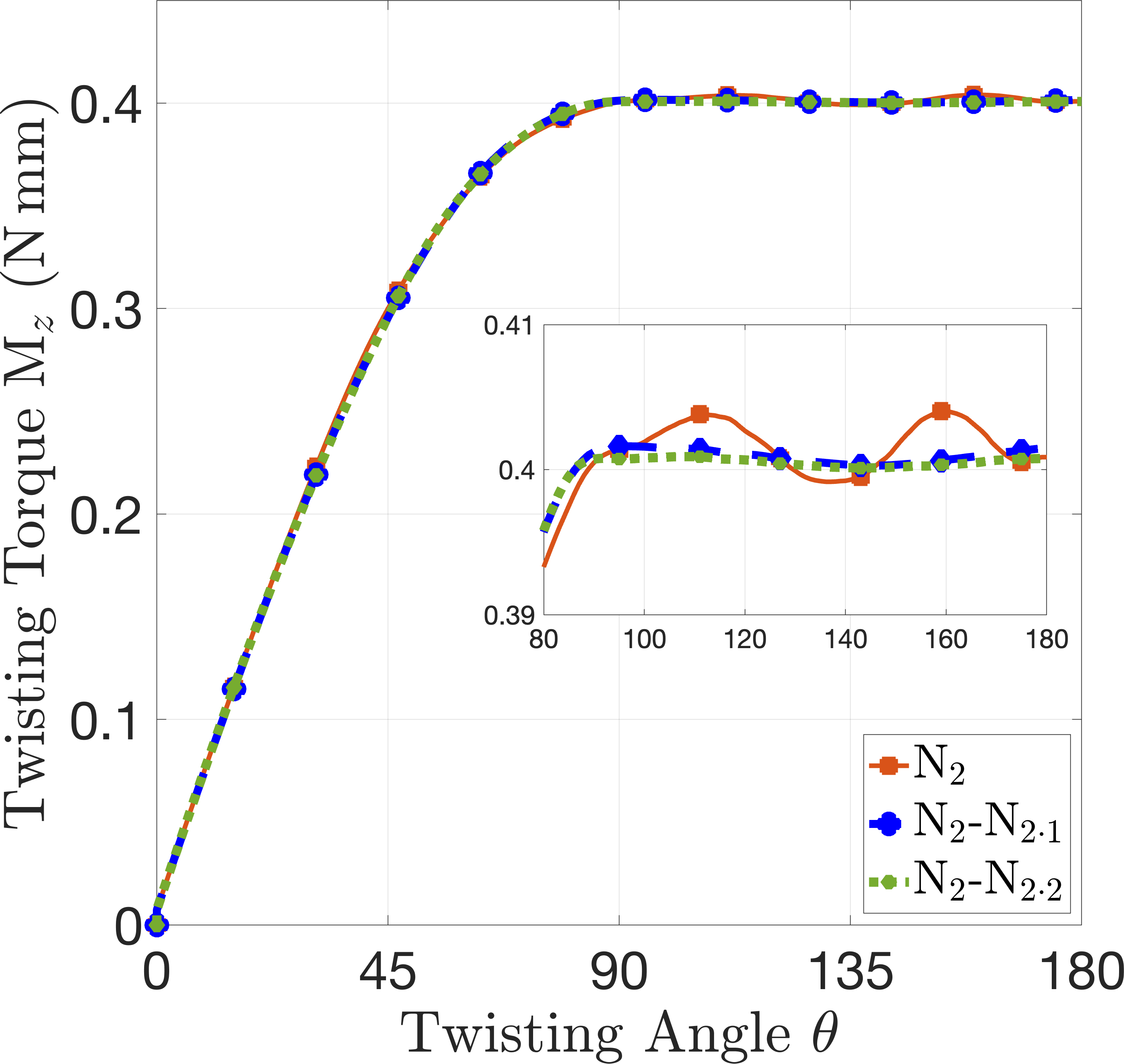}\label{fig:twist_m2}} 
	\caption{Frictional twisting contact: Evolution of total torque over the prescribed rotation angle with the VO and standard NURBS discretizations at mesh (a) m$_1$ and (b) m$_2$.} \label{fig:twisting_fresult}
\end{figure}

Figures.~\ref{fig:twist_m1} and~\ref{fig:twist_m2} show the evolution of the twisting torques over the applied rotation angles with different NURBS discretizations at mesh levels m$_1$ and m$_2$. The result obtained with N$_2$ discretization at m$_2$ agrees with that of Corbett and Sauer~\cite{Corbett2015}. It can be observed that up to the prescribed rotation angle, approx. $\theta = 100^o$, the torque around the vertical axis increases continuously due to the transition of the contact region from the sticking to the slipping. After this transition, once full slipping is achieved, the torque is expected to stabilize at a constant maximum value for the remainder of the twisting process.

However, as seen from Figs.~\ref{fig:twist_m1} and ~\ref{fig:twist_m2}, the standard N$_2$ based NURBS discretization yields a poor approximation of the torque, exhibiting both an underestimation of the average value of torque and non-physical oscillation. These inaccuracies arise from using a coarse mesh, which limits the ability of N$_2$ to accurately capture the contact responses. In contrast, the VO-based N$_2- $N$ _{2\cdot 1} $ discretization, although slightly over-predicts the average torque value at mesh m$_1$, it significantly improves the accuracy. It reduces the deviation of the average torque value by approx. $82\%$ and $79\%$  compared to N$_2$ discretization at mesh levels m$_1$ and m$_2$. Furthermore, applying the N$_2- $N$ _{2\cdot 2} $ at the same mesh level enhances accuracy only marginally compared to N$_2- $N$ _{2\cdot 1} $, consistent with the observations from the previous two examples. 

\begin{figure}[!ht]
	\centering
	{\includegraphics[scale=0.28]{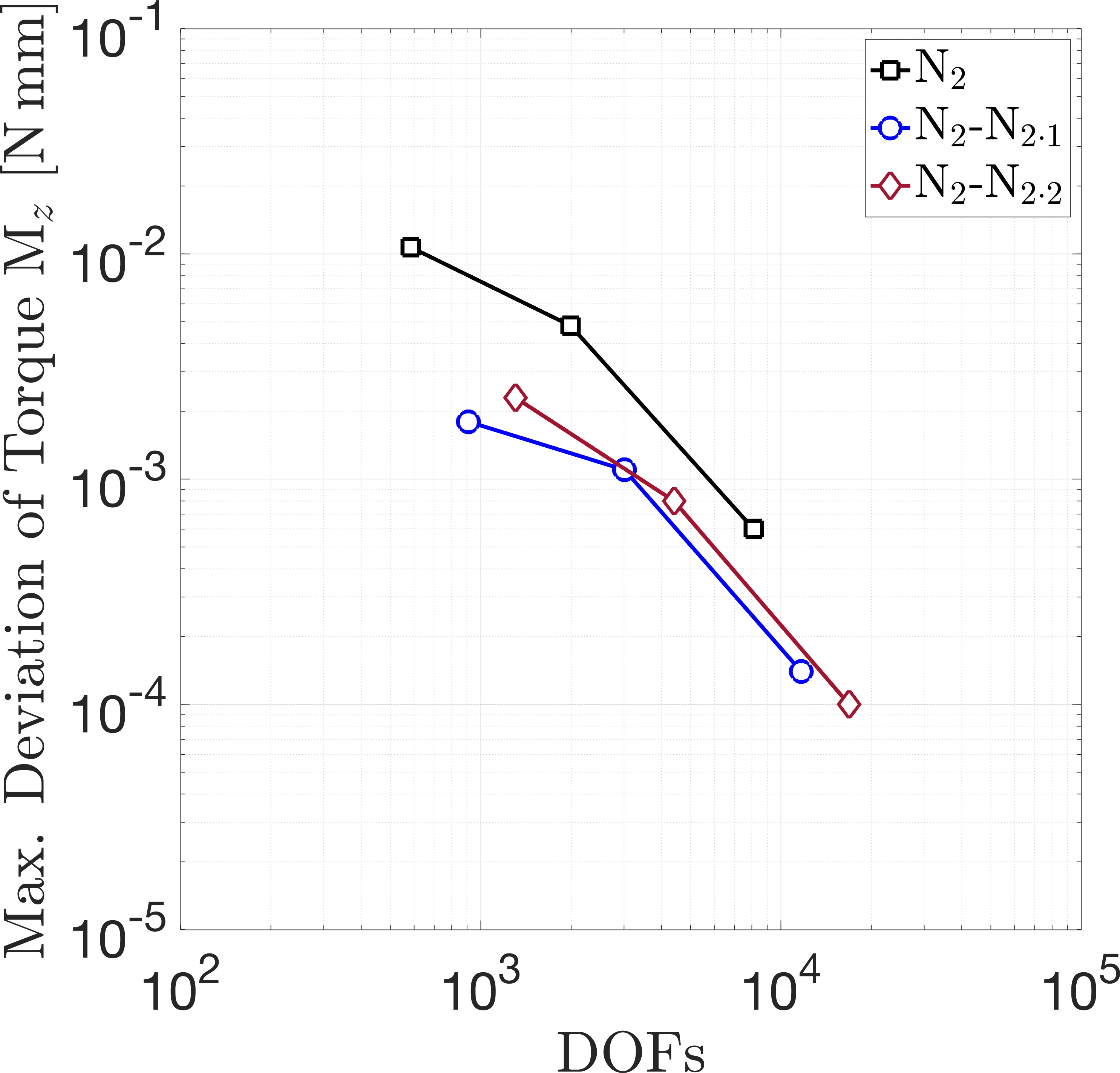}}
	\caption{Frictional twisting contact: Convergence of the maximum deviation of the average torque with mesh refinement for VO and standard NURBS discretizations.} \label{fig:twisting_conv}
\end{figure}

Next, the convergence of the maximum deviation of the torque is examined for various discretizations across three different meshes: m$_1$, m$_2$, and m$_3$, as summarized in Table~\ref{table:twisting_element}. Figure~\ref{fig:twisting_conv} illustrates the convergences of the deviation of the maximum average value of torque for the N$_2$, N$_2- $N$ _{2\cdot 1} $, as well as N$_2- $N$ _{2\cdot 2} $. As shown, all discretizations exhibit reduced torque deviation with increasing mesh resolution, and the convergence rates remain nearly the same for different discretizations. Notably, N$_2- $N$ _{2\cdot 1} $ accomplishes much higher accurate results than N$_2$ at a fixed mesh resolution. 

A closer look reveals that the result obtained with N$_2- $N$ _{2\cdot 1} $ at the first two mesh levels m$_1$ and m$_2$ closely matches with those of N$_2$ at subsequent finer mesh levels m$_2$ and m$_3$. To gain such an improvement in the accuracy, N$_2- $N$ _{2\cdot 1} $ takes approx. $2.2$ and $2.7$ times fewer DOFs than N$_2$, demonstrating a substantial gain in computational efficiency. Such a gain in the computational efficiency with the VO-based discretization is consistent with the trends observed in the previous two examples.

\section{Conclusions} \label{sec:conclusion}
For the accurate, robust, and efficient treatment of three-dimensional large deformation frictional contact using IGA, this work extends the VO-based NURBS discretization method from~\cite{Agrawal2020} to three dimensions. The central feature of the proposed method is that it enables controllable order elevation of a NURBS discretized tri-variate geometry such that the higher-order NURBS are used only for the evaluation of the contact integrals, while the minimum order of NURBS for the computation in the remaining of the 3D bulk domain, without the loss of the advantageous features of the IGA technique. 

In this work, we combine the Gauss-point-to-surface contact algorithm with the presented method to develop a simple yet computationally efficient technique for isogeometric contact analysis. Within the presented formulation, the penalty method is used to regularize the contact constraints, and the normal and sticking contact constraints are enforced independently at each quadrature point. The results demonstrate that the proposed VO NURBS discretization yields much more accurate results, even at a very coarse mesh, when compared to the widely popular standard N$_2$ based NURBS discretization. Specifically, the contact responses obtained with the VO discretization exhibit minimal oscillatory response at a fixed mesh with respect to the standard N$_2$ discretization. Moreover, we show that to attain an accuracy similar to N$_2$, VO uses much lower DOFs, hence considerably enhancing the computational efficiency of the IGA technique for the 3D small and large deformation contact problem with or without considering friction. This gain in efficiency is attributed to the use of higher-order NURBS in the contact region, which increases the conforming ability of the contact surface, enabling more accurate capture of both the normal and tangential contact forces even at a coarse mesh. Additionally, the obtained results show that further increasing the interpolation order of NURBS for the contact region with VO yields only marginal improvements in the accuracy, suggesting that VO-based N$_2$-N$_{2\cdot 1}$ discretization offers an optimal balance between the accuracy and computational efficiency. Furthermore, the simplicity of the proposed method facilitates its direct integration into the standard IGA framework with minimum modifications. Future work focuses on extending the capabilities of the presented method to more complex contact scenarios, such as self-contact and dynamic contact involving large deformation and frictional sliding.

\section*{Acknowledgments}
The author would like to gratefully acknowledge the anonymous reviewers for their critical comments and valuable suggestions, which have significantly improved the quality of the manuscript.

\section*{Declaration of competing interest}
The corresponding author states that there is no conflict of interest.

\section*{Data availability}
Data will be made available on request.

\section*{Appendix: Tangent matrices for frictional contact} \label{appendix:contact_matrices} 
In this section, we provide the contact tangent matrices that are obtained from the linearization of the contact force vectors $ \mathbf{f}^{ke}_{\mathrm{c}} $ given in Sec.~(\ref{sec:final_weak_form}). The final expressions of the contact tangent matrices in the framework of IGA are given by
\begin{equation}\label{eq:contact_tangent_matrices}
	\begin{aligned}
		\mathbf{k}_{\textrm{c}}^{\textrm{ss}e} &:= \frac{\partial \mathbf{f}_{\mathrm{c}}^{\mathrm{s}e}}{\partial \mathbf{u}^{\mathrm{s}e}}  = -\int_{\Gamma^{\textrm{s}e}_c} (\mathbf{R}^{\textrm{s}e})^{\mathrm{T}} \, \frac{\partial \bm{t}_c}{\partial \mathbf{u}^{\textrm{s}e}} ~\textrm{d}\Gamma - \int_{\Gamma^{\textrm{s}}_c} (\mathbf{R}^{\textrm{s}e})^{\mathrm{T}} \, \bm{t}_c \otimes {\bm{\uptau}}^{\alpha}_{\mathrm{s}} \, \mathbf{R}^{\textrm{s}e}_{,\alpha}  ~\textrm{d}\Gamma,~\,\alpha = \{1,2\} \\
		\mathbf{k}_{\textrm{c}}^{\textrm{sm}e} &:= \frac{\partial \mathbf{f}_{\mathrm{c}}^{\mathrm{s}e}}{\partial \mathbf{u}^{\mathrm{m}e}}  = \int_{\Gamma^{\textrm{s}e}_c} (\mathbf{R}^{\textrm{s}e})^{\mathrm{T}} \frac{\partial \bm{t}_c}{\partial \mathbf{u}^{\textrm{m}e}}  ~\textrm{d}\Gamma		\\
	\end{aligned}
\end{equation}
\begin{equation}
	\begin{aligned}
		\mathbf{k}_{\textrm{c}}^{\textrm{ms}e} &:= \frac{\partial \mathbf{f}_{\mathrm{c}}^{\mathrm{m}e}}{\partial \mathbf{u}^{\mathrm{s}e}} = \int_{\Gamma^{\textrm{s}e}_c} (\mathbf{R}^{\textrm{m}e}_{,\alpha})^{\mathrm{T}}\, \bm{t}^{\textrm{s}}\otimes\frac{\partial \bar{\xi}_{\textrm{m}}^{\alpha}}{\partial \mathbf{u}^{\textrm{s}e}}  ~\textrm{d}\Gamma \\
		&~~~~~~~~~~~~~~~~~~+ \int_{\Gamma^{\textrm{s}e}_c} (\mathbf{R}^{\textrm{m}e})^{\mathrm{T}} \, \frac{\partial \bm{t}_c}{\partial \mathbf{u}^{\textrm{s}e}} ~\textrm{d}\Gamma + \int_{\Gamma^{\textrm{s}}_c} (\mathbf{R}^{\textrm{m}e})^{\mathrm{T}} \, \bm{t}_c \otimes {\bm{\uptau}}^{\alpha}_{\mathrm{s}} \, \mathbf{R}^{\textrm{s}e}_{,\alpha} ~\textrm{d}\Gamma  \\
		\mathbf{k}_{\textrm{c}}^{\textrm{mm}e} &:= \frac{\partial \mathbf{f}_{\mathrm{c}}^{\mathrm{m}e}}{\partial \mathbf{u}^{\mathrm{m}e}}  = \int_{\Gamma^{\textrm{s}}_c} (\mathbf{R}^{\textrm{m}e}_{,\alpha})^{\mathrm{T}}\, \bm{t}_c \otimes \frac{\partial \bar{\xi}_{\textrm{m}}^{\alpha}}{\partial \mathbf{u}^{\textrm{m}e}} 
		~\textrm{d}\Gamma  +  \int_{\Gamma^{\textrm{s}}_c}  (\mathbf{R}^{\textrm{m}e})^{\mathrm{T}} \, \frac{\partial \bm{t}_c}{\partial \mathbf{u}^{\textrm{m}e}} ~\textrm{d}\Gamma\,. ~~~~~~~~~~~~~~~~~~\nonumber
	\end{aligned}
\end{equation}
Here, $ \mathbf{R}^{\textrm{s}e} = \mathbf{R}^{\mathrm{s}e}(\bm{\xi}^{\textrm{s}}) $ and $ \mathbf{R}^{\textrm{m}e} = \mathbf{R}^{\mathrm{m}e}(\bar{\bm{\xi}}_{\textrm{m}}) $, and following~\cite{unbiased2015} 
\begin{equation}
	\frac{\partial {\bar{\xi}}_{\mathrm{m}}^{\alpha}}{\partial \mathbf{u}^{\mathrm{s}e}} ~=~ \frac{\partial {\bar{\xi}}_{\mathrm{m}}^{\alpha}}{\partial \bm{x}^{\mathrm{s}e}} \,\frac{\partial \bm{x}^{\mathrm{s}e}}{\partial \mathbf{u}^{\mathrm{s}e}} ~=~ \bar{c}^{\alpha \beta} \, \bar{\bm{\uptau}}_{\beta} \, \mathbf{R}^{\mathrm{s}e}\,.
\end{equation}
In the above equation, the contravariant components $ \bar{c}^{\alpha \beta}_{\mathrm{m}} $ of a tensor $ \bm{c}_{\mathrm{m}} $ are computed using
\begin{equation}
	[\bar{{c}}^{\alpha \beta}] = \left[\bar{m}_{\alpha \beta} - \textrm{g}_{\textrm{N}}\,	\bar{k}_{\alpha \beta} \right]^{-1}\,.
\end{equation}

\noindent  The derivative $ {\partial {\bar{\xi}}_{\mathrm{m}}^{\alpha}}/{\partial \mathbf{u}^{\mathrm{m}e}} $ is given as~\cite{unbiased2013}
\begin{equation}
	\begin{aligned}		
		\frac{\partial {\bar{\xi}}_{\mathrm{m}}^{\alpha}}{\partial \mathbf{u}^{\mathrm{m}e}} = -\bar{c}^{\alpha \beta} \left(\bar{\bm{\uptau}}_{\beta} \, \mathbf{R}^{\mathrm{m}e} - \mathbf{g}_{\mathrm{N}} \, \bm{\bar{n}} \, \mathbf{R}^{\mathrm{m}e}_{\,,\beta}  \right)\,,
	\end{aligned}
\end{equation}

\noindent The partial derivatives of the total contact traction is given as
\begin{equation}
	\frac{\partial \bm{t}_c}{\partial \mathbf{u}^{ke}} = \frac{\partial \bm{t}_{\textrm{N}}}{\partial \mathbf{u}^{ke}} - \frac{\partial \bm{t}_{\textrm{T}}}{\partial \mathbf{u}^{ke}}, ~~~~~~~~\textrm{where } k = \{\mathrm{s,m}\}\,.
\end{equation}
\noindent  For normal contact, the derivatives of the contact traction are given as
\begin{equation}\label{eq:dtN_1}
	\begin{aligned}
		\frac{\partial \bm{t}_{\mathrm{N}}}{\partial \mathbf{u}^{\mathrm{s}e}} &= -\epsilon_{\text{N}}\bm{\bar{n}}\frac{\partial \mathrm{g}_{\mathrm{N}}}{\partial \mathbf{u}^{\mathrm{s}e}} -  \epsilon_{\text{N}}\, \mathrm{g}_{\mathrm{N}}\frac{\partial \bm{\bar{n}}}{\partial \mathbf{u}^{\mathrm{s}e}}\,, \\
		\frac{\partial \bm{t}_{\mathrm{N}}}{\partial \mathbf{u}^{\mathrm{m}e}} &= -\epsilon_{\text{N}}\bm{\bar{n}}\frac{\partial \mathrm{g}_{\mathrm{N}}}{\partial \mathbf{u}^{\mathrm{m}e}} - \epsilon_{\text{N}}\, \mathrm{g}_{\mathrm{N}}\frac{\partial \bm{\bar{n}}}{\partial \mathbf{u}^{\mathrm{m}e}}\,.
	\end{aligned}
\end{equation}

\noindent For tangential stick, the partial derivatives of tangential contact traction are given by
\begin{equation}\label{eq:dTt_duk_stick}
\begin{aligned}
	\frac{\partial \bm{t}_{\textrm{T}}}{\partial \mathbf{u}^{se}} &= \epsilon_{\text{T}} \, \bar{c}^{\alpha \beta} \bar{\bm{\uptau}}_{\alpha} \otimes \bar{\bm{\uptau}}_{\beta} \, \mathbf{R}^{\textrm{s}e}\,,~~~\textrm{and} \\
	\frac{\partial \bm{t}_{\textrm{T}}}{\partial \mathbf{u}^{me}} &= -\epsilon_{\text{T}} \, \bar{c}^{\alpha \beta} \bar{\bm{\uptau}}_{\alpha} \otimes \bar{\bm{\uptau}}_{\beta} \, \mathbf{R}^{\mathrm{m}e} + \epsilon_{\text{T}}\left( \mathbf{R}^{\mathrm{m}e} - \mathbf{R}^{\mathrm{m}e}({\bm{\xi}}_{\textrm{m}\,{sl}\,n}) \right) + \epsilon_{\text{T}} \, \mathrm{g}_{\mathrm{N}}\, \bar{c}^{\alpha \beta} \bar{\bm{\uptau}}_{\alpha} \otimes \bar{\bm{n}}\, \mathbf{R}^{\textrm{m}e}_{,\beta}\,\,.
\end{aligned}
\end{equation}

\noindent For a tangential slip step, the partial derivatives of the tangential contact traction are given as
\begin{equation}
	\begin{aligned}
		\frac{\partial \bm{t}_{\textrm{T}}}{\partial \mathbf{u}^{se}} = &\left( -\mu_f\,  \epsilon_{\text{N}} \, \bm{n}_{\textrm{T}} \otimes \bar{\bm{n}} - \epsilon_{\text{T}}\,\frac{\mu_f\, \epsilon_{\mathrm{N}}\,\mathrm{g}_{\mathrm{N}}}{|| \bm{t}^{\textrm{trial}}_{\textrm{T}}||} [\bm{I} - \bm{n}_{\textrm{T}}\otimes \bm{n}_{\textrm{T}}] \, \bar{c}^{\alpha \beta} \,\bar{\bm{\uptau}}_{\alpha} \otimes \bar{\bm{\uptau}}_{\beta}  \right) \mathbf{R}^{\textrm{s}e}
		\,, \\
		\frac{\partial \bm{t}_{\textrm{T}}}{\partial \mathbf{u}^{me}} =  &\left( \mu_f\, \epsilon_{\text{N}} \, \bm{n}_{\textrm{T}} \otimes \bar{\bm{n}} + \epsilon_{\text{T}}\,\frac{\mu_f \,\epsilon_{\mathrm{N}}\,\textrm{g}_{\textrm{N}}}{|| \bm{t}^{\textrm{trial}}_{\textrm{T}}||} [\bm{I} - \bm{n}_{\textrm{T}}\otimes \bm{n}_{\textrm{T}}]\, \bar{c}^{\alpha\beta} \,\bar{\bm{\uptau}}_{\alpha}\otimes \bar{\bm{\uptau}}_{\beta}  \right)\mathbf{R}^{\textrm{m}e} \\ 
		& - ~\epsilon_{\text{T}} \frac{\mu_f \,\epsilon_{\mathrm{N}}\,\textrm{g}_{\textrm{N}}}{|| \bm{t}^{\textrm{trial}}_{\textrm{T}}||} (\bm{I} - \bm{n}_{\textrm{T}}\otimes \bm{n}_{\textrm{T}})\left( \mathbf{R}^{\mathrm{m}e} - \mathbf{R}^{\mathrm{m}e}(\xi_{sl\, n}^{\textrm{m}}) \right)  \\
		& - ~\epsilon_{\text{T}} \, \frac{\mu_f\,\epsilon_{\textrm{T}}\,\mathrm{g}_{\mathrm{N}}}{||\bm{t}^{\textrm{trial}}_{\textrm{T}}||} (\bm{I} - \bm{n}_{\textrm{T}}\otimes \bm{n}_{\textrm{T}})\,\bar{c}^{\alpha \beta}\, \bar{\bm{\uptau}}_{\alpha} \otimes \bar{\bm{n}} \, \mathbf{R}^{\mathrm{m}e}_{,\beta}\,\,.
	\end{aligned}
\end{equation}

\begingroup
\fontsize{10pt}{10pt}\selectfont
\printbibliography
\endgroup
\end{document}